\documentclass[11pt, reqno]{amsart}
\usepackage{amsfonts,latexsym,enumerate}
\usepackage{amsmath}
\usepackage{amscd}
\usepackage{float,amsmath,amssymb,mathrsfs,bm,multirow,graphics}
\usepackage[dvips]{graphicx}
\usepackage[percent]{overpic}
\usepackage[pdftex]{color}
\usepackage{amsaddr}
\usepackage[numbers,sort&compress]{natbib}

\newcommand{\nequation}{\setcounter{equation}{0}}
\renewcommand{\theequation}{\mbox{\arabic{section}.\arabic{equation}}}
\newcommand{\R}{{\Bbb R}}

\newcommand{\C}{{\Bbb C}}

\newcommand{\proofbegin}{\noindent{\it Proof.\quad}}
\newcommand{\proofend}{\hfill$\Box$\bigskip}
\newcommand{\proofendcontinue}{\hfill \raisebox{.8mm}[0cm][0cm]{$\bigtriangledown$}\bigskip}

\newcommand{\tr}{\text{\upshape tr\,}}

\newcommand{\re}{\text{\upshape Re\,}}

\newcommand{\im}{\text{\upshape Im\,}}

\newcommand{\ntlim}{\lim^\angle}

\def\XXint#1#2#3{{\setbox0=\hbox{$#1{#2#3}{\int}$}
\vcenter{\hbox{$#2#3$}}\kern-.5\wd0}}

\allowdisplaybreaks

\newtheorem{theorem}{Theorem}[section]

\newtheorem{lemma}[theorem]{Lemma}

\newtheorem{remark}[theorem]{Remark}

\newtheorem{figuretext}{Figure}

\input epsf
\title[Nonlinear Fourier transforms]
{Nonlinear Fourier transforms and the mKdV equation in the quarter plane}

\author{Jonatan Lenells}
\address{Department of Mathematics, KTH Royal Institute of Technology, \\ 100 44 Stockholm, Sweden.}
\email{jlenells@kth.se}

\begin{document}
\begin{abstract} 
\noindent
The unified transform method introduced by Fokas can be used to analyze initial-boundary value problems for integrable evolution equations. %In analogy with the inverse scattering transform on the line, 
The method involves several steps, including the definition of spectral functions via nonlinear Fourier transforms and the formulation of a Riemann-Hilbert problem.
We provide a rigorous implementation of these steps in the case of the mKdV equation in the quarter plane under limited regularity and decay assumptions. We give detailed estimates for the relevant nonlinear Fourier transforms. 
Using the theory of $L^2$-RH problems, we consider the construction of quarter plane solutions which are $C^1$ in time and $C^3$ in space. 
\end{abstract}

\maketitle

\noindent
{\small{\sc AMS Subject Classification (2010)}: 35Q53, 35Q15, 41A60.}

\noindent
{\small{\sc Keywords}: Riemann-Hilbert problem, nonlinear Fourier transform, spectral analysis.}

%\tableofcontents

\section{Introduction}\nequation
Initial value problems for integrable evolution equations can be analyzed via the inverse scattering transform (IST) cf. \cite{AC1991, FT2007}. Starting with the initial data, certain spectral functions (often referred to as reflection and transmission coefficients) are defined via a nonlinear Fourier transform. Since the time evolution of these spectral functions is simple, the solution at time $t$ can be recovered via the solution of an inverse problem. The inverse problem is most conveniently formulated as a Riemann-Hilbert (RH) problem whose jump matrix involves the given spectral functions. 

Following the many successes of the inverse scattering approach, one of the main open problems in the area of integrable systems in the late twentieth century was the extension of the IST formalism to initial-boundary value (IBV) problems, see \cite{AC1991}. Such an extension was introduced by Fokas in \cite{F1997} (see also \cite{F2002, Fbook}) and has subsequently been developed and applied by several authors \cite{AK2014, BS2003, BFS2006, K2010, K2003, L2013, Ldnls, PP2010, MQ2014, IIK2011, IS2012, PP2009, XF2013, XF2014, P2007, TF2008, MF2011}.  
In analogy with the IST on the line, the unified transform of \cite{F1997} relies for the analysis of an IBV problem on the definition of several spectral functions via nonlinear Fourier transforms and on the formulation of a RH problem. 

In this paper, we provide a rigorous study of the nonlinear Fourier transforms and RH problems relevant for the analysis of the mKdV equation 
\begin{align}\label{mkdv}
u_t + 6 \lambda u^2u_x - u_{xxx} = 0,  \qquad \lambda = \pm 1,
\end{align}
in the quarter plane $\{x \geq 0, t \geq 0\}$. The unified transform method presents the solution of this problem in terms of the solution of a RH problem, which is defined in terms of four spectral functions $\{a(k), b(k), A(k), B(k)\}$, see \cite{BFS2004}. The functions $a(k)$Ê and $b(k)$ Êcan be viewed as half-line nonlinear Fourier transforms of the initial data, while the functions $A(k)$  and $B(k)$ Êcan be viewed as half-line nonlinear Fourier transforms of the boundary values. We present detailed estimates for these half-line transforms, which can be used to formulate an appropriate RH problem under limited regularity and decay assumptions. In particular, we derive uniform asymptotic expansions for large $k$ and give conditions under which these expansions can be differentiated termwise. We also show how smoothness and decay of the initial and boundary values translate into decay and smoothness of the spectral functions, respectively. Finally, using the theory of $L^2$-RH problems, we consider the construction of quarter plane solutions of (\ref{mkdv}) which are $C^1$ in time and $C^3$ in space. Our presentation can be viewed as an extension of \cite{BFS2004},  where equation (\ref{mkdv}) was analyzed on the half-line under less explicit regularity assumptions.
We present our results for the mKdV equation for definiteness, but similar arguments are applicable also to other integrable equations such as the nonlinear Schr\"odinger, KdV, sine-Gordon, and Camassa-Holm equations. 

The rigorous study of nonlinear Fourier transforms and RH problems is rather involved even in simple cases. For example, for the KdV equation on the line, the relevant nonlinear Fourier transform is that associated with the one-dimensional Schr\"odinger operator $-\partial_x^2 + u_0(x)$, and Deift and Trubowitz presented the rigorous analysis of this operator and of the associated transform in an elegant but long paper \cite{DT1979}. 
In many cases, to avoid technical details, results relying on inverse scattering techniques are presented under rather vague assumptions on the given data such as ``sufficient smoothness and decay''. This can sometimes be motivated by the fact that the qualitative outcomes of the theory are independent of the precise assumptions. However, there are situations where more precise formulations are vital also qualitatively. In the context of IBV problems, physically relevant examples of such situations include: 

\begin{enumerate}[1.]
\item The derivation of long-time asymptotics via the nonlinear steepest descent method. 

\item Problems with asymptotically time-periodic data. 

\item Problems whose initial and boundary data are not compatible to all orders at the points of the boundary for which $t = 0$. 
%The noncompatibility leads to a family of RH problems with subtle singular behavior as $(x,t) \to (0,0)$. 

\item Problems with step-like initial and/or boundary profiles. 
%The step-like discontinuities lead to slowly decaying spectral functions with possible singularities. 
\end{enumerate}
For the derivation of long-time asymptotics, the decay properties of the boundary values are particularly important---if the boundary values do not decay as $t \to \infty$, the asymptotic formulas will receive additional contributions from the boundary. For problems with $t$-periodic data, the spectral functions may have branch cuts, hence the formulation of a RH problem is intricate and a detailed understanding of the half-line Fourier transforms is crucial, see \cite{BIK2009, BKS2009, tperiodicI}. In addition to providing a detailed study of the mKdV equation in the quarter plane, the present paper intends to lay the foundation for future explorations of the above topics. 

The analysis of IBV problems is more involved than the analysis of pure initial value problems. Thus, although the unified transform method and the IST formalism share several characteristics, there are important differences. Let us comment on a few of these differences relevant for the present study: 

(a) For an initial value problem on the line, the IST formalism utilizes two eigenfunctions which are normalized at plus and minus infinity. For the corresponding quarter plane problem, the unified transform method utilizes one eigenfunction normalized at spatial infinity, one normalized at temporal infinity, and one normalized at the origin. The latter eigenfunction is entire, but has a more complicated large $k$ Êbehavior than the eigenfunctions normalized at infinity. Indeed, in order to be correctly normalized at the origin, this eigenfunction must be a linear combination of two solutions, one of which admits an expansion in $1/k$ whereas the other is exponentially small in each asymptotic sector. The contribution from the exponentially small solution can sometimes be ignored, but it becomes important as $k$ approaches the anti-Stokes lines that form the boundary of the asymptotic sector. In fact, as $k \to \infty$ along one of these lines, the two solutions contribute terms of comparable order asymptotically.

(b) Since the spectral functions for the IBV problem are defined via half-line transforms, they do not have rapid decay as $k \to \infty$ even in the case of smooth data. Indeed, the half-line transform of a function $f(x)$, $x \geq 0$, can be viewed as the transform on the whole line of $f_e(x)$ where $f_e = f$ for $x \geq 0$ and $f_e = 0$ for $x < 0$. If $f(0) \neq 0$, the dicontinuity of $f_e(x)$ at $x=0$ implies that the transform only decays as $1/k$ as $k \to \infty$. The rigorous formulation of a RH problem therefore involves a careful study of asymptotic expansions. 

(c) The initial data and the boundary values of a solution of an IBV problem are not independent. The relationship between the initial and boundary values is encoded in a relation among the spectral functions called the global relation. When formulating the main RH problem, we must assume that the global relation is fulfilled (see equation (\ref{GR}) below).

Section \ref{prelsec} contains some definitions and notational conventions. 
In Section \ref{xsec}, we consider the definition of $a(k)$ Êand $b(k)$.
In Section \ref{tsec}, we consider the definition of $A(k)$ Êand $B(k)$.
In Section \ref{mkdvsec}, we consider the construction of quarter plane solutions of (\ref{mkdv}) which are $C^1$ in time and $C^3$ in space.
A few results on $L^2$-RH problems are collected in the appendix. 

\section{Preliminaries} \nequation\label{prelsec}
\subsection{Lax pair}
The mKdV %mKdVII
equation (\ref{mkdv}) admits the Lax pair
\begin{align}\label{mulax}
\begin{cases}
\mu_x - ik[\sigma_3, \mu] = \mathsf{U} \mu, \\
\mu_t + 4ik^3[\sigma_3, \mu] = \mathsf{V} \mu,
\end{cases}
\end{align}
where $\mu(x,t,k)$ is a $2\times 2$-matrix valued eigenfunction, $k\in \C$ is the spectral parameter, and 
\begin{align}\nonumber
& \sigma_3 = \begin{pmatrix} 1 & 0 \\ 0 & -1 \end{pmatrix}, \qquad
\mathsf{U} = \begin{pmatrix} 0 & u \\ \lambda u & 0 \end{pmatrix}, 
	\\ \label{mathsfUVdef}
& \mathsf{V} = \begin{pmatrix} -2i \lambda k u^2  & -4uk^2 + 2iku_x - 2 \lambda u^3 + u_{xx} \\
-4\lambda k^2 u - 2i\lambda k u_x - 2 u^3 + \lambda u_{xx} & 2i \lambda k u^2 \end{pmatrix}.
\end{align}
The versions of (\ref{mkdv}) with $\lambda = 1$ and $\lambda = -1$ are referred to as the defocusing and focusing mKdV equations, respectively.

\subsection{Notation}
\begin{figure}
\begin{center}
\bigskip \bigskip
\begin{overpic}[width=.45\textwidth]{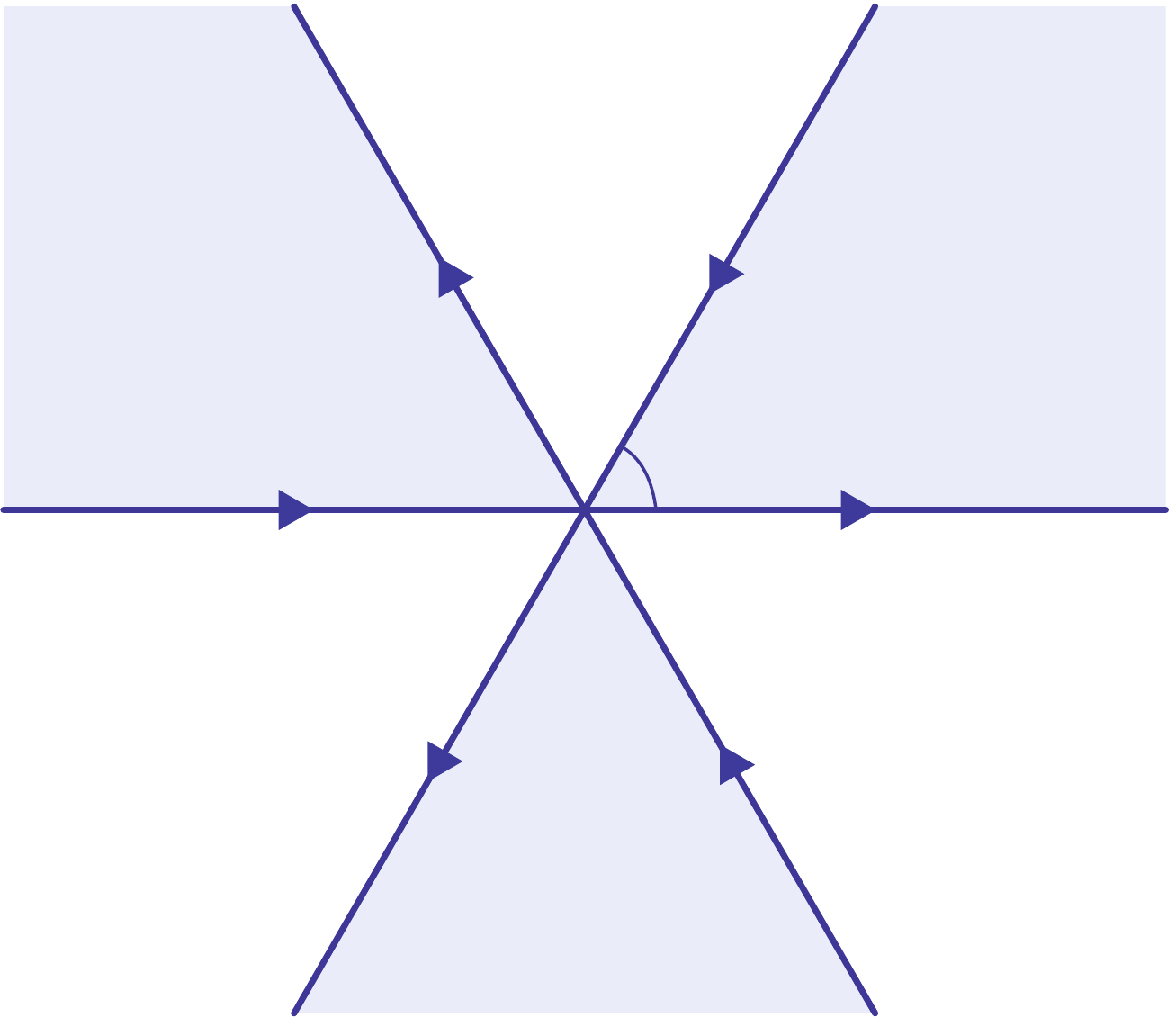}
      \put(77,57){$D_3$}
      \put(46,70){$D_4$}
      \put(15,57){$D_3$}
      \put(15,25){$D_2$}
      \put(46,12){$D_1$}
      \put(77,25){$D_2$}
      \put(57,47){$\pi/3$}
      \put(103,41){$\Gamma$}
      \end{overpic}
     \begin{figuretext}\label{DjsmKdV.pdf}
       The contour $\Gamma$ and the domains $\{D_j\}_1^4$ in the complex $k$-plane.
     \end{figuretext}
     \end{center}
\end{figure}

For a $2 \times 2$ matrix $A$, we let $A^{(d)}$ and $A^{(o)}$ denote the diagonal and off-diagonal parts, respectively. If $A$ is an $n \times m$ matrix, we define $|A|$ by 
$$|A| = \sqrt{\sum_{i,j} |A_{ij}|^2} = \sqrt{\tr \bar{A}^T A}.$$
Then $|A + B| \leq |A| + |B|$ and $|AB| \leq |A| |B|$. 
For a contour $\gamma \subset \C$ and $1 \leq p \leq \infty$, we write $A \in L^p(\gamma)$  if $|A|$ belongs to $L^p(X)$. We define $\|A\|_{L^p(\gamma)} := \| |A|\|_{L^p(\gamma)}$. Note that $A \in L^p(\gamma)$ if and only if each entry $A_{ij}$ belongs to $L^p(\gamma)$. 
We let $\{\sigma_j\}_1^3$ denote the three Pauli matrices. 
We let $\hat{\sigma}_3$ act on a $2\times 2$ matrix $A$ by $\hat{\sigma}_3A = [\sigma_3, A]$, i.e. $e^{\hat{\sigma}_3} A = e^{\sigma_3} A e^{-\sigma_3}$. 
For a $2 \times 2$ matrix $A$, we let $[A]_1$ and $[A]_2$ denote the first and second columns of $A$.
We let $\C_+ = \{\im k > 0\}$ and $\C_- = \{\im k < 0\}$ denote the open upper and lower half-planes; $\bar{\C}_\pm = \C_\pm \cup \R$ will denote their closures. 
The notation $k \in (\C_+, \C_-)$ indicates that the first and second columns are valid for $k \in \C_+$ and $k \in \C_-$, respectively. Given $x \in \R$, $[x]$ will denote the integer part of $x$.
The open domains $\{D_j\}_1^4$ of the complex $k$-plane are defined by (see Figure \ref{DjsmKdV.pdf})
\begin{align}\nonumber
D_1 = \{\im k < 0\} \cap \{\im k^3 > 0\},  \qquad
D_2 = \{\im k < 0\} \cap \{\im k^3 < 0\}, 
	\\ \nonumber
D_3 = \{\im k > 0\} \cap \{\im k^3 > 0\},  \qquad
D_4 = \{\im k > 0\} \cap \{\im k^3 < 0\}.
\end{align}
We let $D_+ = D_1 \cup D_3$ and $D_- = D_2 \cup D_4$. 
We let $\Gamma = \R \cup e^{i\pi/3}\R \cup e^{2i\pi/3}\R$ denote the contour separating the $D_j$'s  oriented so that $D_+$ lies to the left of $\Gamma$.
Throughout the paper $C$ denotes a generic constant. 

\section{Spectral analysis of the $x$-part}\nequation\label{xsec}
Let $u(x)$ be a real-valued function defined for $x\geq 0$ and let
$$\mathsf{U}(x) = \begin{pmatrix} 0 & u(x) \\ \lambda u(x) & 0 \end{pmatrix}.$$
Consider the linear differential equation
\begin{align}\label{xpartF}
X_x - ik[\sigma_3, X] = \mathsf{U}X,
\end{align}
where $X(x,k)$ is a $2 \times 2$-matrix valued eigenfunction and $k \in \C$ is a spectral parameter.
We define two $2 \times 2$-matrix valued solutions of (\ref{xpartF}) as the solutions of the linear Volterra integral equations
\begin{subequations}\label{Fpmdef}
\begin{align}  \label{Fpmdefa}
 & X(x,k) = I + \int_{\infty}^x e^{-i k(x'-x)\hat{\sigma}_3} (\mathsf{U}X)(x',k) dx',
  	\\
&  Y(x,k) = I + \int_{0}^x e^{-i k(x'-x)\hat{\sigma}_3} (\mathsf{U}Y)(x',k) dx'.
\end{align}
\end{subequations}
The proof of the following theorem is given in Section \ref{xproof1subsec}.

\begin{theorem}\label{xth1}
Let $m \geq 1$ and $n \geq 1$ be integers. 
Suppose 
\begin{align}\label{uassump}
\begin{cases}
u \in C^{m+1}([0,\infty)), &
	\\
(1+x)^{n}\partial^i u(x) \in L^1([0,\infty)), & i = 0,1, \dots, m+1.
\end{cases}
\end{align}
Then the equations (\ref{Fpmdef}) uniquely define two $2 \times 2$-matrix valued solutions $X$ and $Y$ of (\ref{xpartF}) with the following properties:
\begin{enumerate}[$(a)$]
\item The function $X(x, k)$ is defined for $x \geq 0$ and $k \in (\bar{\C}_+, \bar{\C}_-)$. For each $k \in (\bar{\C}_+, \bar{\C}_-)$, the function $X(\cdot, k) \in C^1([0,\infty))$ satisfies (\ref{xpartF}).

\item The function $Y(x, k)$ is defined for $x \geq 0$ and $k \in \C$. For each $k \in \C$,  the function $Y(\cdot, k) \in C^1([0,\infty))$ satisfies (\ref{xpartF}). 

\item For each $x \geq 0$, the function $X(x,\cdot)$ is bounded and continuous for $k \in (\bar{\C}_+, \bar{\C}_-)$ and analytic for $k \in (\C_+, \C_-)$.

\item For each $x \geq 0$, the function $Y(x,\cdot)$ is an entire function of $k \in \C$ which is bounded for $k \in (\bar{\C}_-, \bar{\C}_+)$.

\item For each $x \geq 0$ and each $j = 1, \dots, n$, the partial derivative $\frac{\partial^j X}{\partial k^j}(x, \cdot)$ has a continuous extension to $(\bar{\C}_+, \bar{\C}_-)$.

\item $X$ and $Y$ satisfy the following estimates: 
\begin{subequations}\label{Fest}
\begin{align}\label{Festa}
& \bigg|\frac{\partial^j}{\partial k^j}\big(X(x,k) - I\big) \bigg| \leq 
\frac{C}{(1+x)^{n-j}}, \qquad x \geq 0, \quad  k \in (\bar{\C}_+, \bar{\C}_-),
	\\\label{Festb}
& \bigg|\frac{\partial^j}{\partial k^j}\Big(Y(x,k) - I\Big) \bigg| \leq 
C\min(x,1)(1 + x)^j, \qquad x \geq 0, \quad k \in (\bar{\C}_-, \bar{\C}_+),
\end{align}
\end{subequations}
for $ j = 0, 1, \dots, n$.
\end{enumerate}
\end{theorem}

\subsection{Behavior as $k \to \infty$}
In addition to the properties listed in Theorem \ref{xth1}, we also need to know the behavior of the eigenfunctions $X$ and $Y$ as $k \to \infty$. To this end, we note that equation (\ref{xpartF}) admits formal power series solutions $X_{formal}$ and $Y_{formal}$, normalized at $x = \infty$ and $x = 0$ respectively, such that
\begin{align}\label{varphiformalx}
& X_{formal}(x,k) = I + \frac{X_1(x)}{k} + \frac{X_2(x)}{k^2} + \cdots,
	\\
& Y_{formal}(x,k) =  I + \frac{Z_1(x)}{k} + \frac{Z_2(x)}{k^2} + \cdots + \bigg(\frac{W_1(x)}{k} + \frac{W_2(x)}{k^2} + \cdots \bigg) e^{-2ikx\sigma_3},
\end{align}
where the coefficients $\{X_j(x), Z_j(x), W_j(x)\}_1^\infty$ satisfy
\begin{align}\label{Fjnormalization}
\lim_{x\to \infty} X_j(x) = 0, \qquad Z_j(0) + W_j(0) = 0, \qquad j \geq 1.
\end{align}
Indeed, substituting
$$X = I + \frac{X_1(x)}{k} + \frac{X_2(x)}{k^2} + \cdots$$
into (\ref{xpartF}), the off-diagonal terms of $O(k^{-j})$ and the diagonal terms of $O(k^{-j-1})$ yield the relations
\begin{align}\label{xrecursive}
\begin{cases}
X_{j+1}^{(o)} = -\frac{i}{2}\sigma_3\big(X_{jx}^{(o)} - \mathsf{U}X_j^{(d)}\big),
	\\
\partial_xX_{j+1}^{(d)} = \frac{i}{2}\sigma_3\mathsf{U}\big(X_{jx}^{(o)} - \mathsf{U}X_j^{(d)}\big).
\end{cases}
\end{align}
Similarly, substituting
$$X = \bigg(\frac{W_1(x)}{k} + \frac{W_2(x)}{k^2} + \cdots\bigg) e^{-2ikx\sigma_3}$$
into (\ref{xpartF}), the diagonal terms of $O(k^{-j})$ and the off-diagonal terms of $O(k^{-j-1})$ yield the relations
%$$H_x - 2ik\sigma_3 H^{(d)} = \mathsf{U}H,$$
\begin{align}\label{xrecursive2}
\begin{cases}
W_{j+1}^{(d)} = - \frac{i}{2}\sigma_3\big(W_{jx}^{(d)} - \mathsf{U}W_j^{(o)}\big),
	\\
\partial_x W_{j+1}^{(o)} = \frac{i}{2}\sigma_3\mathsf{U}\big(W_{jx}^{(d)} - \mathsf{U}W_j^{(o)}\big).
\end{cases}
\end{align}

The coefficients $\{X_j(x), Z_j(x), W_j(x)\}$ are determined recursively from (\ref{Fjnormalization})-(\ref{xrecursive2}), the equations obtained from (\ref{xrecursive}) by replacing $\{X_j\}$ with $\{Z_j\}$, and the initial assignments 
$$X_{-1} = 0, \qquad X_0 = I, \qquad Z_{-1} = 0, \qquad Z_0 = I, \qquad W_{-1} = 0, \qquad W_0 = 0.$$ 
Let 
$$\sigma_\lambda = \begin{pmatrix} 0 & 1 \\ - \lambda & 0 \end{pmatrix}.$$
Then the first few coefficients are given by
\begin{align}\nonumber
X_1(x) = &\; \frac{iu(x)}{2} \sigma_\lambda 
- \frac{i\lambda}{2} \sigma_3 \int_{\infty}^x u^2(x') dx',
	\\ \nonumber
X_2(x) = &\; \frac{1}{4}\big[2iu (X_1)_{22} + u_x\big]\sigma_3\sigma_\lambda
 + I\int_{\infty}^x \frac{\lambda u}{4}\big[2iu (X_1)_{22} + u_x\big] dx',
	\\ \nonumber
X_3(x) = &\; \frac{1}{8}\big[2(X_1)_{22} u_x + i\lambda u^3 + 4iu (X_2)_{22} - iu_{xx}\big]\sigma_\lambda
	\\\nonumber
& - \sigma_3\int_{\infty}^x \frac{iu}{8}\big[u^3 + 4\lambda u(X_2)_{22} - 2i\lambda (X_1)_{22} u_x - \lambda u_{xx})\big]  dx',
	\\ \nonumber
X_4(x) = &\; \frac{1}{16}\big[5\lambda u^2u_x + 4(X_2)_{22}u_x + 2i(\lambda u^3 - u_{xx}) (X_1)_{22} + 8i u (X_3)_{22} - u_{xxx}\big]\sigma_3\sigma_\lambda
	\\ \nonumber
& + I\int_{\infty}^x \frac{\lambda u}{16}\big[5\lambda u^2u_x + 4(X_2)_{22}u_x + 2i(\lambda u^3 - u_{xx}) (X_1)_{22} 
	\\ \label{Fjexplicit}
& + 8i u (X_3)_{22} - u_{xxx}\big] dx',
\end{align}
\begin{align}\nonumber
Z_1(x) = &\; \frac{iu(x)}{2} \sigma_\lambda 
- \frac{i\lambda}{2} \sigma_3 \int_{0}^x u^2 dx',
	\\ \nonumber
Z_2(x) = &\; \frac{1}{4}\big[2iu (Z_1)_{22} + u_x\big]\sigma_3\sigma_\lambda
 + I\bigg\{\int_0^x \frac{\lambda u}{4}\big[2iu (Z_1)_{22} + u_x\big] dx' + \frac{\lambda}{4}u^2(0)\bigg\},
	\\ \nonumber
Z_3(x) = &\; \frac{1}{8}\big[2(Z_1)_{22} u_x + i\lambda u^3 + 4iu (Z_2)_{22} - iu_{xx}\big]\sigma_\lambda
	\\\nonumber
& - \sigma_3\int_{0}^x \frac{iu}{8}\big[u^3 + 4\lambda u(Z_2)_{22} - 2i\lambda (Z_1)_{22} u_x - \lambda u_{xx})\big]  dx',
	\\ \nonumber
Z_4(x) = &\; \frac{1}{16}\big[5\lambda u^2u_x + 4(Z_2)_{22}u_x + 2i(\lambda u^3 - u_{xx}) (Z_1)_{22} + 8i u (Z_3)_{22} - u_{xxx}\big]\sigma_3\sigma_\lambda
	\\ \nonumber
& + I\bigg\{\int_0^x \frac{\lambda u}{16}\big[5\lambda u^2u_x + 4(Z_2)_{22}u_x + 2i(\lambda u^3 - u_{xx}) (Z_1)_{22} 
	\\ \label{Fjmexplicit}
& + 8i u (Z_3)_{22} - u_{xxx}\big] dx' + \frac{3u^4(0)}{16} - \frac{\lambda}{8}u(0)u_{xx}(0) + \frac{\lambda u_x^2(0)}{16}\bigg\},
\end{align}
and
\begin{align}\nonumber
W_1(x) = &\; -\frac{iu(0)}{2}\sigma_\lambda,
	\\ \nonumber
W_2(x) = &\; -\frac{\lambda i u (W_1)_{12}}{2}I - \sigma_3\sigma_\lambda \bigg\{\frac{\lambda i}{2}(W_1)_{12}\int_0^x u^2  dx' + \frac{u_x(0)}{4}\bigg\},
	\\ \nonumber
W_3(x) = &\; -\frac{\lambda}{4}(-2iu(W_2)_{12} + (W_1)_{12}u_x)\sigma_3 
	\\\nonumber
&- \frac{\lambda}{4}\sigma_\lambda \bigg\{\int_0^x u\big[2iu (W_2)_{12} - (W_1)_{12} u_x\big]dx' + iu^3(0) - \frac{i\lambda}{2}u_{xx}(0) \bigg\},
	\\ \nonumber
W_4(x) = &\; -\frac{i\lambda}{8}(\lambda u^3(W_1)_{12} + 4u (W_3)_{12} + 2i(W_2)_{12} u_x - (W_1)_{12} u_{xx})I
	\\\nonumber
& - \frac{i}{8}\sigma_3\sigma_\lambda\bigg\{ \int_0^x u\big[u^3(W_1)_{12} + \lambda(4u (W_3)_{12} + 2i(W_2)_{12} u_x - (W_1)_{12} u_{xx})\big]dx' 
	\\ \label{calFjmexplicit}
& - 3i\lambda u^2(0)u_x(0) + \frac{i}{2}u_{xxx}(0)\bigg\}.
\end{align}
If $u(x)$ has a finite degree of regularity and decay, only finitely many coefficients $\{X_j, Z_j, W_j\}$ are well-defined. The following result, whose proof is given in Section  \ref{xproof2subsec}, describes the behavior of $X$ and $Y$ as $k \to \infty$.

\begin{theorem}\label{xth2}
Suppose $u(x)$ satisfies (\ref{uassump}) for some integers $m \geq 1$ and $n \geq 1$.
Then, as $k \to \infty$, $X$ and $Y$ coincide to order $m$ with $X_{formal}$ and $Y_{formal}$ respectively in the following sense: The functions
\begin{align}\label{Fhatdef}
&\hat{X}(x,k) = I + \frac{X_1(x)}{k} + \cdots + \frac{X_{m+1}(x)}{k^{m+1}},
	\\
&\hat{Y}(x,k) = I + \frac{Z_1(x)}{k} + \cdots + \frac{Z_{m+1}(x)}{k^{m+1}}
 + \bigg(\frac{W_1(k)}{k} + \cdots + \frac{W_{m+1}(x)}{k^{m+1}}\bigg) e^{-2ikx\sigma_3},
\end{align}
are well-defined and there exists a $K > 0$ such that
\begin{subequations}\label{varphiasymptotics}
\begin{align}\label{varphiasymptoticsa}
\bigg|\frac{\partial^j}{\partial k^j}\big(X - \hat{X}\big) \bigg| \leq 
\frac{C}{|k|^{m+1}(1+x)^{n-j}}, \qquad x \geq 0, \quad  j = 0, 1, \dots, n,
\end{align}
for all $k \in (\bar{\C}_+, \bar{\C}_-)$ with  $|k| > K$, and
\begin{align}\label{varphiasymptoticsb}
\bigg|\frac{\partial^j}{\partial k^j}\big(Y - \hat{Y}\big) \bigg| \leq 
\frac{C(1 + x)^{j+2}e^{\frac{Cx}{|k|^{m+1}}}}{|k|^{m+1}}, \qquad x \geq 0, \quad j = 0, 1, \dots, n,
\end{align}
for all $k \in (\bar{\C}_-, \bar{\C}_+)$ with  $|k| > K$, and
\begin{align}\label{varphiasymptoticsc}
\bigg|\frac{\partial^j}{\partial k^j}\big((Y - \hat{Y})e^{2ikx\sigma_3}\big) \bigg| \leq 
\frac{C(1 + x)^{j+2}e^{\frac{Cx}{|k|^{m+1}}}}{|k|^{m+1}}, \qquad x \geq 0, \quad j = 0, 1, \dots, n,
\end{align}
\end{subequations}
for all $k \in (\bar{\C}_+, \bar{\C}_-)$ with  $|k| > K$.
\end{theorem}

\begin{remark}\upshape
It is relatively straightforward to derive from equations (\ref{Fpmdef}) that $X$ and $Y$ have the properties listed in Theorem \ref{xth1}. For the proof of Theorem \ref{xth2}, which is much more involved, the equations (\ref{Fpmdef}) are not suitable; instead we will consider the equations satisfied by the `errors' $\hat{X}^{-1} X$, $\hat{Z}^{-1} Y$, and $\hat{W}^{-1} Y e^{2ikx\sigma_3}$, where $\hat{Y} = \hat{Z} + \hat{W} e^{-2ikx\sigma_3}$.
\end{remark}

\subsection{Proof of Theorem \ref{xth1}}\label{xproof1subsec} 
In view of the symmetries 
\begin{align}\label{Fpmsymm}
F(x,k) = \begin{cases} 
\sigma_1\overline{F(x, \bar{k})} \sigma_1, & \lambda = 1, \\
\sigma_2\overline{F(x, \bar{k})} \sigma_2, & \lambda = -1,
\end{cases}
\end{align}
which are valid for $F = X$Ê and $F = Y$, it is enough to prove the theorem for $[X]_2$ and $[Y]_2$. We first consider the construction of $[X]_2$.

Let $\Psi(x,k)$ denote the second column of $X(x,k)$. Then, by (\ref{Fpmdefa}), 
\begin{align}\label{xPsiVolterra}
\Psi(x,k) = \begin{pmatrix} 0 \\ 1 \end{pmatrix} - \int_x^\infty E(x,x',k)\mathsf{U}(x')\Psi(x',k) dx', \qquad x \geq 0,
\end{align}
where
$$E(x,x',k) =  \begin{pmatrix} e^{-2ik(x' - x)} & 0 \\ 0 & 1 \end{pmatrix}.$$
We  use successive approximations to show that the Volterra integral equation (\ref{xPsiVolterra}) has a unique solution $\Psi(x,k)$ for each $k \in \bar{\C}_-$.
Let $\Psi_0 = \big(\begin{smallmatrix} 0 \\ 1 \end{smallmatrix}\big)$ and define $\Psi_l$ for $l \geq 1$ inductively by 
\begin{align*}
\Psi_{l+1}(x,k) = - \int_x^\infty E(x,x',k) \mathsf{U}(x') \Psi_l(x', k) dx', \qquad x \geq 0, \quad k \in \bar{\C}_-.
\end{align*}
Then
\begin{align}\label{xPhiliterated}
\Psi_l(x,k) = & (-1)^l \int_{x = x_{l+1} \leq x_l \leq \cdots \leq x_1 < \infty} \prod_{i = 1}^l E(x_{i+1}, x_i, k) \mathsf{U}(x_i) \Psi_0 dx_1 \cdots dx_l. 
\end{align}
Using the estimates 
\begin{align*}
  |E(x, x', k)| < C, \qquad 0 \leq x \leq x' < \infty, \quad k \in \bar{\C}_-,
\end{align*}
and
$$\|\mathsf{U}\|_{L^1([x,\infty))} < \frac{C}{(1+x)^n}, \qquad x \geq 0,$$
we find
\begin{align}\nonumber
|\Psi_l(x,k)| \leq \; & C \int_{x \leq x_l \leq \cdots \leq x_1 < \infty} \prod_{i = 1}^l  |\mathsf{U}(x_i)| |\Psi_0| dx_1 \cdots dx_l 
	\\ \label{xPhilestimate}
\leq & \; \frac{C}{l!} \|\mathsf{U}\|_{L^1([x,\infty))}^l 
\leq \frac{C}{l!} \bigg(\frac{C}{(1+x)^n}\bigg)^l, \qquad x \geq 0, \quad k \in \bar{\C}_-.
\end{align}
Hence the series
$$\Psi(x,k) = \sum_{l=0}^\infty \Psi_l(x,k)$$
converges absolutely and uniformly for $x \geq 0$ and $k \in \bar{\C}_-$ to a continuous solution $\Psi(x,k)$ of (\ref{xPsiVolterra}). 
Moreover,
\begin{align}\label{xPsiminusf}
|\Psi(x,k) - \Psi_0| \leq \sum_{l=1}^\infty | \Psi_l(x,k)| \leq \frac{C}{(1+x)^n}, \qquad x \geq 0, \quad k \in \bar{\C}_-,
\end{align}
which proves (\ref{Festa}) for $j = 0$. Using the estimate
\begin{align}\label{xE1bound}
|\partial_k^jE(x, x', k)| < C(1+|x' - x|)^j, \qquad 0 \leq x \leq x' < \infty, \quad k \in \bar{\C}_-, \quad j = 0, 1, \dots, n,
\end{align}
with $j = 1$ to differentiate under the integral sign in (\ref{xPhiliterated}), we see that $\Psi_l(x, \cdot)$ is analytic in $\C_-$ for each $l$; the uniform convergence then proves that $\Psi$  is analytic in $\C_-$. 

It remains to show that $[X]_2 = \Psi$ satisfies $(e)$ and $(f)$ for $j = 1, \dots, n$.
Let 
\begin{align}\label{xLambda0}
\Lambda_0(x,k) = \;& - \int_x^\infty (\partial_kE)(x,x',k) \mathsf{U}(x') \Psi(x', k) dx'.
\end{align}
Differentiating the integral equation (\ref{xPsiVolterra}) with respect to $k$, we find that $\Lambda := \partial_k\Psi$ satisfies
\begin{align}\label{xintegraleq3}
\Lambda(x,k) = \Lambda_0(x,k) - \int_x^\infty E(x,x',k) \mathsf{U}(x') \Lambda(x', k) dx'
\end{align}
for each $k$ in the interior of $\bar{\C}_-$; the differentiation can be justified by dominated convergence using (\ref{xE1bound}) and a Cauchy estimate for $\partial_k\Psi$.
We seek a solution of (\ref{xintegraleq3}) of the form $\Lambda = \sum_{l=0}^\infty \Lambda_l$, where the $\Lambda_l$'s are defined by replacing $\{\Psi_l\}$ by $\{\Lambda_l\}$ in (\ref{xPhiliterated}). Proceeding as in (\ref{xPhilestimate}), we find
\begin{align}\label{EE0}
|\Lambda_l(x,k)| \leq \frac{C}{l!}\|\Lambda_0(\cdot, k)\|_{L^\infty([x,\infty))}\bigg(\frac{C}{(1+x)^n}\bigg)^l, \qquad x \geq 0, \quad k \in \bar{\C}_-.
\end{align}
Using (\ref{xPsiminusf}) and (\ref{xE1bound}) in (\ref{xLambda0}), we obtain
\begin{align}\label{EE1}
|\Lambda_0(x,k)| \leq \frac{C}{(1+x)^{n-1}}, \qquad x \geq 0, \quad k \in \bar{\C}_-.
\end{align}
In particular $\|\Lambda_0(\cdot, k)\|_{L^\infty([x,\infty))}$ is bounded for $x \geq 0$ and $k \in \bar{\C}_-$.
Thus, $\sum_{l=0}^\infty \Lambda_l$ converges uniformly on $[0,\infty) \times \bar{\C}_-$ to a continuous solution $\Lambda$ of (\ref{xintegraleq3}), which satisfies the following analog of (\ref{xPsiminusf}):
\begin{align}\label{EE2}
|\Lambda(x,k) - \Lambda_0(x,k)| \leq 
\frac{C}{(1+x)^{n}}, \qquad x \geq 0,  \quad k \in \bar{\C}_-.
\end{align}
In view of equations (\ref{EE1}) and (\ref{EE2}), we conclude that $[X]_2 = \Psi$ satisfies $(e)$ and $(f)$ for $j = 1$.

Proceeding inductively, we find that $\Lambda^{(j)} := \partial_k^j\Psi$ satisfies an integral equation of the form
\begin{align*}
\Lambda^{(j)}(x,k) = \Lambda^{(j)}_0(x,k)
- \int_x^\infty E(x,x',k) \mathsf{U}(x') \Lambda^{(j)}(x', k) dx', 
\end{align*}
where
$$\big|\Lambda^{(j)}_0(x, k)\big| \leq 
\frac{C}{(1+x)^{n-j}}, \qquad x \geq 0, \quad k \in \bar{\C}_-.$$
If $1 \leq j \leq n$, then $\|\Lambda^{(j)}_0(\cdot, k)\|_{L^\infty([x,\infty))}$ is bounded for $x \geq 0$ and $k \in \bar{\C}_-$; hence the associated series $ \Lambda^{(j)} = \sum_{l=0}^\infty \Lambda^{(j)}_l$ converges uniformly on $[0,\infty) \times \bar{\C}_-$ to a continuous solution with the desired properties.

We now consider the the construction of $[Y]_2$. In this case, we still let $\Psi_0 = \big(\begin{smallmatrix} 0 \\ 1 \end{smallmatrix}\big)$, but instead of (\ref{xPhiliterated}) we now introduce $\{\Psi_l\}_1^\infty$ by 
\begin{align*}
\Psi_l(x,k) = & (-1)^l \int_{0 \leq x_1 \leq  \cdots \leq x_l \leq x < \infty} \prod_{i = 1}^l E(x_{i+1}, x_i, k) \mathsf{U}(x_i) \Psi_0 dx_1 \cdots dx_l,
\end{align*}
where $x \geq 0$ and $k \in \C$. This leads to an entire function $\Psi(x,k) = \sum_{l=0}^\infty \Psi_l(x,k)$ satisfying (\ref{xpartF}). Moreover, as in (\ref{xPhilestimate}), we find
\begin{align}\nonumber
|\Psi_l(x,k)| \leq \frac{C}{l!} \|\mathsf{U}\|_{L^1([0,x])}^l 
\leq \frac{C\min(x,1)^l}{l!}, \qquad x \geq 0, \quad k \in \bar{\C}_+,
\end{align}
which leads to the following analog of (\ref{xPsiminusf}):
\begin{align*}
|\Psi(x,k) - \Psi_0| \leq \sum_{l=1}^\infty | \Psi_l(x,k)| \leq C\min(x,1), \qquad x \geq 0, \quad k \in \bar{\C}_+.
\end{align*}
This proves (\ref{Festb}) for $j = 0$. Letting
\begin{align*}
\Lambda_0(x,k) = \;& \int_0^x (\partial_kE)(x,x',k) \mathsf{U}(x') \Psi(x', k) dx',
\end{align*}
we find that $\Lambda := \partial_k\Psi$ satisfies
\begin{align*}
\Lambda(x,k) = \Lambda_0(x,k) + \int_0^x E(x,x',k) \mathsf{U}(x') \Lambda(x', k) dx'
\end{align*}
The proof of (\ref{Festb}) for $j = 1$ follows from the following analogs of (\ref{EE0})-(\ref{EE2}):
\begin{align*}
& |\Lambda_l(x,k)| \leq \frac{C}{l!}\|\Lambda_0(\cdot, k)\|_{L^\infty([0,x])}, \qquad 
	\\ 
&  |\Lambda_0(x,k)| \leq C \int_0^x (1 + |x-x'|) |u(x')| dx'
\leq Cx + C x \int_0^x |u(x')| dx'
 \leq Cx, 
	\\ \nonumber
& |\Lambda(x,k) - \Lambda_0(x,k)| \leq Cx,
\end{align*}
which are valid for $x \geq 0$ and $k \in \bar{\C}_+$; the proof for $j \geq 2$ is similar.
\proofend

\subsection{Proof of Theorem \ref{xth2}}\label{xproof2subsec}
We first consider $[X]_2$. Our first goal is to show that $\hat{X}$ is well-defined and invertible for $k$ large enough. 

\medskip\noindent
{\bf Claim 1.}  $\{X_j(x)\}_1^{m+1}$ are $C^1$ functions of $x \geq 0$ satisfying
\begin{align}\label{Fjbounds}
\begin{cases}
(1 + x)^n X_j(x) \in L^1([0,\infty)) \cap L^\infty([0,\infty)), \\
(1 + x)^n X_j'(x) \in L^1([0,\infty)),
\end{cases} \qquad j = 1, \dots, m+1.
\end{align}

{\it Proof of Claim 1.}
The assumption (\ref{uassump}) implies 
\begin{align}\label{uibound}
  |u^{(i)}(x)| \leq \frac{C}{(1 + x)^n}, \qquad x \geq 0, \quad i = 0, 1, \dots, m.
\end{align}
Indeed, if $i = 0, 1, \dots, m$, and $x \geq 0$, then
\begin{align}\label{1xuiest}
\big|(1 + x)^n u^{(i)}(x)\big| \leq |u^{(i)}(0)| + \int_0^x \big|n(1 + x')^{n-1}u^{(i)}(x') + (1+x')^nu^{(i+1)}(x')\big|dx' < C.
\end{align}

Let $S_j$ refer to the statement
$$
\begin{cases}
X_j^{(o)} \in C^{m+2-j}([0,\infty)), \qquad
X_j^{(d)} \in C^{m+3-j}([0,\infty)), 
	\\
\text{$(1 + x)^n \partial^i X_j^{(o)} \in L^1([0,\infty))$ for $i = 0,1, \dots, m+2-j$,}
	\\
\text{$(1 + x)^n \partial^i X_j^{(d)} \in L^1([0,\infty))$ for $i = 0,1, \dots, m+3-j$.}
\end{cases}$$
By (\ref{uibound}),
$$(1 + x)^{n+1} \int_x^\infty u^2(x') dx'
\leq \int_x^\infty (1 + x')^{n+1} u^2(x') dx' 
\leq C \int_x^\infty (1 + x') |u(x')| dx' \to 0$$
as $x \to \infty$. Consequently, an integration by parts yields
\begin{align}\nonumber
& \|(1 + x)^n X_1^{(d)}\|_{L^1([0,\infty))}
 \leq C \int_0^\infty (1+x)^n \int_x^\infty u^2(x') dx' dx
	\\ \label{xF1dest}
& \leq C + C\int_0^\infty \frac{(1+x)^{n+1}}{n+1} u^2(x) dx
 \leq  C + C \int_0^\infty (1+x) |u(x)| dx < \infty.
 \end{align}
Using the estimate (\ref{xF1dest}) and the expression (\ref{Fjexplicit}) for $X_1$, we conclude that $S_1$ holds.
Similar estimates together with the relations (\ref{xrecursive}) imply that if $1 \leq j \leq m$ and $S_j$ holds, then $S_{j+1}$ also holds. Thus, by induction, $\{S_j\}_{j=1}^{m+1}$ hold. This shows that $\{X_j\}_1^{m+1}$ are $C^1$ functions satisfying $(1+x)^n\partial^iX_j(x) \in L^1([0,\infty))$ for $i =0,1$ and $j = 1, \dots, m+1$. The boundedness of $(1+x)^nX_j(x)$ follows by an estimate analogous to (\ref{1xuiest}). 
\proofendcontinue

\medskip\noindent
{\bf Claim 2.} There exists a $K> 0$ such that $\hat{X}(x,k)^{-1}$ exists for all $k \in \C$ with $|k| \geq K$. 
Moreover, letting $A = ik\sigma_3 + \mathsf{U}$ and 
\begin{align}\label{Ahatdef}
\hat{A}(x,k) = \big(\hat{X}_x(x,k) + ik \hat{X}(x,k) \sigma_3\big)\hat{X}(x,k)^{-1}, \qquad x \geq 0, \quad |k| \geq K,
\end{align}
the difference $\Delta(x,k) = A(x,k) - \hat{A}(x,k)$ satisfies
\begin{align}\label{preABsmall}
|\Delta(x,k)| \leq \frac{Cf(x)}{|k|^{m+1}(1 + x)^n}, \qquad x \geq 0, \quad |k| \geq K,
\end{align}
where $f$ is a function in $L^1([0,\infty)) \cap C([0,\infty))$. 
In particular,
\begin{align}\label{ABsmall}
\|\Delta(\cdot, k)\|_{L^1([x,\infty))} \leq \frac{C}{|k|^{m+1}(1+x)^n}, \qquad x \geq 0, \quad |k| \geq K.
\end{align}

{\it Proof of Claim 2.}
By Claim 1, there exists a bounded continuous function $g \in L^1([0,\infty))$ such that
\begin{align}\label{Fjxbound}
|X_j(x)| \leq \frac{g(x)}{(m+1)(1+x)^n}, \qquad x \geq 0, \quad j = 1, \dots, m+1.
\end{align}
In particular,
$$\bigg| \sum_{j=1}^{m+1} \frac{X_j(x)}{k^j}\bigg| \leq \frac{g(x)}{|k|(1+x)^n}, \qquad x \geq 0, \quad |k| \geq 1.$$
Choose $K > \max(1, \|g\|_{L^\infty([0,\infty))})$. Then $\hat{X}(x,k)^{-1}$ exists whenever $|k| \geq K$ and is given by the absolutely and uniformly convergent Neumann series
$$\hat{X}(x,k)^{-1} = \sum_{l =0}^\infty \bigg(-\sum_{j=1}^{m+1} \frac{X_j(x)}{k^j}\bigg)^l, \qquad x \geq 0, \quad |k| \geq K.$$
Furthermore,
\begin{align}\label{tailestimate}
\bigg|\sum_{l = m+2}^\infty \bigg(-\sum_{j=1}^{m+1} \frac{X_j(x)}{k^j}\bigg)^l\bigg|
\leq \sum_{l = m+2}^\infty \bigg(\frac{g(x)}{|k|(1+x)^n}\bigg)^l
\leq \frac{Cg(x)}{|k|^{m+2}(1+x)^n},
\end{align}
for $x \geq 0$ and $|k| \geq K$. Now let $Q_0(x) + \frac{Q_1(x)}{k} + \frac{Q_2(x)}{k^2} + \cdots$ be the formal power series expansion of $\hat{X}(x,k)^{-1}$ as $k \to \infty$, i.e.
\begin{align*}
& Q_0(x) = I, \quad Q_1(x) = -X_1(x), \quad Q_2(x) = X_1(x)^2 - X_2(x), 
	\\
& Q_3(x) = X_1(x)X_2(x) + X_2(x) X_1(x) - X_1(x)^3 - X_3(x), \quad \dots.
\end{align*}
Equation (\ref{Fjxbound}) and the inequality (\ref{tailestimate}) imply that the function $\mathcal{E}(x,k)$ defined by
$$\mathcal{E}(x,k) = \hat{X}(x,k)^{-1} - \sum_{j=0}^{m+1} \frac{Q_j(x)}{k^j}$$
satisfies
\begin{align}\label{Eestimate}
|\mathcal{E}(x,k)| \leq \frac{Cg(x)}{|k|^{m+2}(1+x)^n}, \qquad x \geq 0, \quad |k| \geq K.
\end{align}

Let $\hat{A}(x,k)$ be given by (\ref{Ahatdef}).
Since $X_{formal}$ is a formal solution of (\ref{xpartF}), the coefficient of $k^{-j}$ in the formal expansion of $\Delta = A - \hat{A}$ as $k \to \infty$  vanishes for $j \leq m$; hence, in view of Claim 1 and (\ref{Eestimate}),
\begin{align}\nonumber
|\Delta| 
& = \bigg|A - (\hat{X}_x + ik \hat{X} \sigma_3)\bigg(\sum_{j=0}^{m+1} \frac{Q_j}{k^j} + \mathcal{E}\bigg)\bigg|
	\\\nonumber
& \leq \frac{Cg(x)}{|k|^{m+1}(1+x)^n} + |(\hat{X}_x + ik \hat{X} \sigma_3)\mathcal{E}|
	\\ \nonumber
& \leq \frac{Cf(x)}{|k|^{m+1}(1+x)^n}, \qquad x \geq 0, \quad |k| \geq K,
\end{align}
where $f$ is a continuous (not necessarily bounded) function in $L^1([0,\infty))$. This proves (\ref{preABsmall}).
\proofendcontinue

Given $K > 0$, we let $\bar{\C}_\pm^K = \bar{\C}_\pm \cap \{|k| \geq K\}$.

\medskip\noindent
{\bf Claim 3.} The Volterra integral equation
\begin{align} \label{integraleq}
\Psi(x,k) = &\; \Psi_0(x,k) - \int_x^\infty E(x,x',k) \Delta(x', k) \Psi(x', k) dx',
\end{align}
where $\Psi_0(x,k) =  [\hat{X}(x,k)]_2$ and
\begin{align}\label{E1E1E1def}
& E(x,x',k) =  \hat{X}(x,k) \begin{pmatrix} e^{-2ik(x' - x)} & 0 \\ 0 & 1 \end{pmatrix} \hat{X}(x',k)^{-1},
\end{align}
has a unique solution $\Psi(x,k)$ for each $k \in \bar{\C}_-^K$. This solution satisfies $\Psi = [X]_2$.

{\it Proof of Claim 3.}
Define $\Psi_l$ for $l \geq 1$ inductively by 
\begin{align*}
\Psi_{l+1}(x,k) = - \int_x^\infty E(x,x',k) \Delta(x', k) \Psi_l(x', k) dx', \qquad x \geq 0, \quad k \in \bar{\C}_-^K.
\end{align*}
Then
\begin{align}\label{Philiterated}
\Psi_l(x,k) = & (-1)^l \int_{x = x_{l+1} \leq x_l \leq \cdots \leq x_1 < \infty} \prod_{i = 1}^l E(x_{i+1}, x_i, k) \Delta(x_i, k) \Psi_0(x_1, k) dx_1 \cdots dx_l. 
\end{align}
Using the estimate 
\begin{align*}
  |E(x, x', k)| < C, \qquad 0 \leq x \leq x' < \infty, \quad k \in \bar{\C}_-^K,
\end{align*}
as well as (\ref{ABsmall}), we find, for $x \geq 0$ and $k \in \bar{\C}_-^K$,
\begin{align}\nonumber
|\Psi_l(x,k)| \leq \; & C \int_{x \leq x_l \leq \cdots \leq x_1 < \infty} \prod_{i = 1}^l  |\Delta(x_i, k)| |\Psi_0(x_1, k)| dx_1 \cdots dx_l 
	\\ \nonumber
\leq & \; \frac{C}{l!}\|\Psi_0(\cdot, k)\|_{L^\infty([x,\infty))} \|\Delta(\cdot, k)\|_{L^1([x,\infty))}^l
	\\ \label{Philestimate}
\leq & \; \frac{C}{l!} \|\Psi_0(\cdot, k)\|_{L^\infty([x,\infty))}\bigg(\frac{C}{|k|^{m+1}(1+x)^n}\bigg)^l.
\end{align}
Since
$$\sup_{k \in \bar{\C}_-^K} \|\Psi_0(\cdot, k)\|_{L^\infty([0,\infty))} < C,$$ 
we find
$$|\Psi_l(x,k)| \leq \frac{C}{l!}\bigg(\frac{C}{|k|^{m+1}(1+x)^n}\bigg)^l, \qquad x \geq 0, \quad k \in \bar{\C}_-^K.$$
Hence the series $\Psi(x,k) = \sum_{l=0}^\infty \Psi_l(x,k)$
converges absolutely and uniformly for $x \geq 0$ and $k \in \bar{\C}_-^K$ to a continuous solution $\Psi(x,k)$ of (\ref{integraleq}).
Moreover,
\begin{align}\label{Psiminusf}
|\Psi(x,k) - \Psi_0(x,k)| \leq \sum_{l=1}^\infty | \Psi_l(x,k)| \leq \frac{C}{|k|^{m+1}(1+x)^n}, \qquad x \geq 0, \quad k \in \bar{\C}_-^K.
\end{align}
It follows from the integral equation (\ref{integraleq}) that $\Psi$ satisfies the second column of (\ref{xpartF}). By (\ref{Psiminusf}), $\Psi(x,k) \sim \begin{pmatrix} 0 \\ 1 \end{pmatrix}$ as $x \to \infty$. Hence, by uniqueness of solution, $\Psi = [X]_2$ for $k \in \bar{\C}_-^K$.
\proofendcontinue

\medskip\noindent
{\bf Claim 4.} $[X]_2$ satisfies (\ref{varphiasymptoticsa}).

{\it Proof of Claim 4.}
Equation (\ref{Psiminusf}) implies that $[X]_2 = \Psi$ satisfies (\ref{varphiasymptoticsa}) for $j = 0$.

A glance at the proof of Claim 2 shows that the inequality (\ref{preABsmall}) can be extended to derivatives of $\Delta$ with respect to $k$:
\begin{subequations}\label{DeltaEbound}
\begin{align}\label{DeltaEbounda}
\big|\partial_k^j\Delta(x,k)\big| \leq \frac{Cf(x)}{|k|^{m+1}(1+x)^n}, \qquad x \geq 0, \quad |k| \geq K, \quad j = 0,1, \dots, n,
\end{align}
where $f \in L^1([0,\infty))$ is a continuous function of $x \geq 0$. We will also need the following estimate for $j = 0,1, \dots, n$:
\begin{align}\label{DeltaEboundb}
|\partial_k^jE(x, x', k)| < C(1+|x' - x|)^j, \qquad 0 \leq x \leq x' < \infty, \quad k \in \bar{\C}_-^K,
\end{align}
\end{subequations}
where $E$ Êis defined in (\ref{E1E1E1def}).
Let 
\begin{align}\label{Lambda0}
\Lambda_0(x,k) = \;& [\partial_k\hat{X}(x,k)]_2 
- \int_x^\infty \frac{\partial}{\partial k}\big[E(x,x',k) \Delta(x', k)\big] \Psi(x', k) dx'.
\end{align}
Differentiating the integral equation (\ref{integraleq}) with respect to $k$, we find that $\Lambda := \partial_k\Psi$ satisfies
\begin{align}\label{integraleq3}
\Lambda(x,k) = \Lambda_0(x,k)
- \int_x^\infty E(x,x',k) \Delta(x', k) \Lambda(x', k) dx'
\end{align}
for each $k$ in the interior of $\bar{\C}_-^K$; the differentiation can be justified by dominated convergence using (\ref{DeltaEbound}) and a Cauchy estimate for $\partial_k\Psi$.
We seek a solution of (\ref{integraleq3}) of the form $\Lambda = \sum_{l=0}^\infty \Lambda_l$,  where the $\Lambda_l$'s are defined by replacing $\{\Psi_l\}$ by $\{\Lambda_l\}$ in (\ref{Philiterated}). Proceeding as in (\ref{Philestimate}), we find
$$|\Lambda_l(x,k)| \leq \frac{C}{l!}\|\Lambda_0(\cdot, k)\|_{L^\infty([x,\infty))}\bigg(\frac{C}{|k|^{m+1}(1+x)^n}\bigg)^l, \qquad x \geq 0, \quad k \in \bar{\C}_-^K.$$
Using (\ref{Psiminusf}) and (\ref{DeltaEbound}) in (\ref{Lambda0}), we obtain
\begin{align}\label{FF1}
|\Lambda_0(x,k) - [\partial_k\hat{X}(x,k)]_2 | \leq \frac{C}{|k|^{m+1}(1+x)^{n-1}}, \qquad x \geq 0, \quad k \in \bar{\C}_-^K.
\end{align}
In particular $\|\Lambda_0(\cdot, k)\|_{L^\infty([x,\infty))}$ is bounded for $k \in \bar{\C}_-^K$ and $x \geq 0$.
Thus, $\sum_{l=0}^\infty \Lambda_l$ converges uniformly on $[0,\infty) \times \bar{\C}_-^K$ to a continuous solution $\Lambda$ of (\ref{integraleq3}), which satisfies the following analog of (\ref{Psiminusf}):
\begin{align}\label{FF2}
|\Lambda(x,k) - \Lambda_0(x,k)| \leq 
\frac{C}{|k|^{m+1}(1+x)^{n}}, \qquad x \geq 0,  \quad k \in \bar{\C}_-^K.
\end{align}
Equations (\ref{FF1}) and (\ref{FF2}) show that $[X]_2 = \Psi$ satisfies (\ref{varphiasymptoticsa}) for $j = 1$.

Proceeding inductively, we find that $\Lambda^{(j)} := \partial_k^j\Psi$ satisfies an integral equation of the form
\begin{align*}
\Lambda^{(j)}(x,k) = \Lambda^{(j)}_0(x,k)
- \int_x^\infty E(x,x',k) \Delta(x', k) \Lambda^{(j)}(x', k) dx', 
\end{align*}
where
$$\big|\Lambda^{(j)}_0(x, k) - [\partial_k^{j} \hat{X}(x, k)]_2 \big| \leq 
\frac{C}{|k|^{m+1}(1+x)^{n-j}}, \qquad x \geq 0, \quad k \in \bar{\C}_-^K.$$
If $1 \leq j \leq n$, then $\|\Lambda^{(j)}_0(\cdot, k)\|_{L^\infty([x,\infty))}$ is bounded for $k \in \bar{\C}_-^K$ and $x \geq 0$; hence the associated series $ \Lambda^{(j)} = \sum_{l=0}^\infty \Lambda^{(j)}_l$ converges uniformly on $[0,\infty) \times \bar{\C}_-^K$ to a continuous solution with the desired properties. 
\proofendcontinue

The above claims prove the theorem for $X$. We now consider $[Y]_2$.

\medskip\noindent
{\bf Claim 5.}  $\{Z_j(x), W_j(x)\}_1^{m+1}$ are $C^1$ functions of $x \geq 0$ satisfying
\begin{align*}
\begin{cases}
Z_j, W_j \in L^\infty([0,\infty)), \\
Z_j', W_j' \in L^1([0,\infty)),
\end{cases} \qquad j = 1, \dots, m+1.
\end{align*}

{\it Proof of Claim 5.}
Let $S_j$ and $\mathcal{S}_j$ refer to the statements
$$
\begin{cases}
Z_j^{(o)} \in C^{m+2-j}([0,\infty)), \qquad
Z_j^{(d)} \in C^{m+3-j}([0,\infty)), 
	\\
\text{$\partial^i Z_j^{(o)} \in L^1([0,\infty))$ for $i = 1, \dots, m+2-j$,}
	\\
\text{$\partial^i Z_j^{(d)} \in L^1([0,\infty))$ for $i = 1, \dots, m+3-j$,}
\end{cases}$$
and
$$\begin{cases}
W_j^{(d)} \in C^{m+2-j}([0,\infty)), \qquad
W_j^{(o)} \in C^{m+3-j}([0,\infty)), 
	\\
\text{$\partial^i W_j^{(d)} \in L^1([0,\infty))$ for $i = 1, \dots, m+2-j$,}
	\\
\text{$\partial^i W_j^{(o)} \in L^1([0,\infty))$ for $i = 1, \dots, m+3-j$,}
\end{cases}$$
respectively.
Using the expressions (\ref{Fjmexplicit}) and (\ref{calFjmexplicit}) for $Z_1$ and $W_1$, we conclude that $S_1$ and $\mathcal{S}_1$ hold. The relations (\ref{xrecursive}) and (\ref{xrecursive2}) imply by induction that $\{S_j, \mathcal{S}_j\}_{j=1}^{m+1}$ hold. This shows that $\{Z_j, W_j\}_1^{m+1}$ are $C^1$ functions satisfying $Z_j', W_j' \in L^1([0,\infty))$ for $j = 1, \dots, m+1$. Integration shows that $Z_j, W_j \in L^\infty([0,\infty))$ for $j = 1, \dots, m+1$. \proofendcontinue

We write
$$\hat{Y}(x,k) = \hat{Z}(x,k) + \hat{W}(x,k) e^{-2ikx\sigma_3},$$
where $\hat{Z}$Ê and $\hat{W}$ are defined by
$$\hat{Z}(x,k) = I + \frac{Z_1(x)}{k} + \cdots + \frac{Z_{m+1}(x)}{k^{m+1}}, \qquad
\hat{W}(x,k) = \frac{W_1(k)}{k} + \cdots + \frac{W_{m+1}(x)}{k^{m+1}}.$$

\medskip\noindent
{\bf Claim 6.}  
There exists a $K> 0$ such that $\hat{Z}(x,k)^{-1}$ and $\hat{W}(x,k)^{-1}$ exist for all $k \in \C$ with $|k| \geq K$. 
Moreover, letting $A = ik\sigma_3 + \mathsf{U}$ and 
\begin{align*}
\begin{cases}
\hat{A}_1(x,k) = \big(\hat{Z}_x(x,k) + ik \hat{Z}(x,k) \sigma_3\big)\hat{Z}(x,k)^{-1}, 
	\\
\hat{A}_2(x,k) = \big(\hat{W}_x(x,k) - ik \hat{W}(x,k) \sigma_3\big)\hat{W}(x,k)^{-1},
\end{cases}
\qquad x \geq 0, \quad |k| \geq K,
\end{align*}
the differences 
$$\Delta_l(x,k) = A(x,k) - \hat{A}_l(x,k), \qquad l = 1,2,$$ 
satisfy
\begin{align*}
|\partial_k^j\Delta_l(x,k)| \leq \frac{C + f(x)}{|k|^{m+1}}, \qquad x \geq 0, \quad |k| \geq K, \quad j = 0, 1, \dots, n, \quad l = 1,2,
\end{align*}
where $f$ is a function in $L^1([0,\infty)) \cap C([0,\infty))$. 
In particular,
\begin{align}\label{EEDeltaest}
\|\Delta_l(\cdot, k)\|_{L^1([0,x])} \leq \frac{Cx}{|k|^{m+1}}, \qquad x \geq 0, \quad |k| \geq K, \quad l = 1,2.
\end{align}

{\it Proof of Claim 6.}
The proof uses Claim 5 and is similar to that of Claim 2.
\proofendcontinue

\medskip\noindent
{\bf Claim 7.} We have
\begin{subequations}
\begin{align}\label{ZYWa}
& \bigg| \frac{\partial^j}{\partial k^j}\Big[\hat{Z}(x,k) e^{ikx\hat{\sigma}_3} \hat{Z}^{-1}(0,k) - \hat{Y}(x,k)\Big]_2 \bigg| \leq C(1+x)^j\frac{1 + |e^{2ikx}|}{|k|^{m+2}}, 
	\\ \label{ZYWb}
& \bigg| \frac{\partial^j}{\partial k^j}\Big[\hat{W}(x,k) e^{-ikx\hat{\sigma}_3} \hat{W}^{-1}(0,k) - \hat{Y}(x,k)e^{2ikx\sigma_3}\Big]_2 \bigg| \leq C(1+x)^j\frac{1 + |e^{-2ikx}|}{|k|^{m+2}},
\end{align}
for all $x \geq 0$, $|k| \geq K$, and $j = 0, 1, \dots, n$.
\end{subequations}

{\it Proof of Claim 7.}
We write
$$Y_{formal}(x,k) = Z_{formal}(x,k) + W_{formal}(x,k)e^{-2ikx\sigma_3},$$
where
$$Z_{formal}(x,k) = I + \frac{Z_1(x)}{k} + \frac{Z_2(x)}{k^2} + \cdots, \qquad
W_{formal}(x,k) = \frac{W_1(k)}{k} + \frac{W_2(x)}{k^2} + \cdots.$$
Then
\begin{align}\label{ZZY}
Z_{formal}(x,k) e^{ikx\hat{\sigma}_3} Z_{formal}^{-1}(0,k) = Y_{formal}(x,k)
\end{align}
formally to all orders in $k$. Indeed, both sides of (\ref{ZZY}) are formal solutions of (\ref{xpartF}) satisfying the same initial condition at $x = 0$. Truncating (\ref{ZZY}) at order $k^{-m-1}$, it follows that
$$\hat{Z}(x,k) e^{ikx\hat{\sigma}_3} \hat{Z}^{-1}(0,k) = \hat{Y}(x,k) + O(k^{-m-2}) + O(k^{-m-2})e^{-2ikx\sigma_3}.$$
Using Claim 5 and estimating the inverse $\hat{Z}^{-1}(0,k)$ as in the proof of Claim 2, we find (\ref{ZYWa}) for $j = 0$. Using the estimate $|\partial_k e^{\pm 2ikx}| \leq C |x e^{\pm 2ikx}|$, we find (\ref{ZYWa}) also for $j \geq 1$. 

Similarly, we have
$$W_{formal}(x,k) e^{-ikx\hat{\sigma}_3} W_{formal}^{-1}(0,k) = Y_{formal}(x,k)e^{2ikx\sigma_3}$$
to all orders in $k$ and truncation leads to (\ref{ZYWb}).
\proofendcontinue

\medskip\noindent
{\bf Claim 8.} $[Y]_2$ satisfies (\ref{varphiasymptoticsb}).

{\it Proof of Claim 8.}
Using that $Y(x,k)$ satisfies (\ref{xpartF}), we compute
\begin{align*}
  (\hat{Z}^{-1}Y)_x & = - \hat{Z}^{-1}\hat{Z}_x \hat{Z}^{-1} Y + \hat{Z}^{-1} Y_x
  	\\
&  = - \hat{Z}^{-1}(\hat{A}_1\hat{Z} - ik\hat{Z}\sigma_3) \hat{Z}^{-1} Y + \hat{Z}^{-1} (AY - ikY\sigma_3)
	\\
&  = \hat{Z}^{-1}\Delta_1 Y + ik[\sigma_3, \hat{Z}^{-1}Y].
\end{align*}
Hence
$$\Big(e^{-ikx\hat{\sigma}_3}\hat{Z}^{-1}Y\Big)_x = e^{-ikx\hat{\sigma}_3}\hat{Z}^{-1} \Delta_1 Y.$$
Integrating and using the initial condition $Y(0,k) = I$, we conclude that $Y$ satisfies the Volterra integral equation
\begin{align}\label{EEYVolterra}
Y(x,k) = \hat{Z}(x,k) e^{ikx\hat{\sigma}_3}\hat{Z}^{-1}(0,k) + \int_0^x \hat{Z}(x,k) e^{ik(x-x')\hat{\sigma}_3} (\hat{Z}^{-1}\Delta_1 Y)(x',k) dx'.
\end{align}
Letting $\Psi = [Y]_2$ and $\Psi_0(x,k) = [\hat{Z}(x,k) e^{ikx\hat{\sigma}_3}\hat{Z}^{-1}(0,k) ]_2$, we can write the second column of (\ref{EEYVolterra}) as
\begin{align}\label{EEYVolterra2}
\Psi(x,k) = \Psi_0(x,k) + \int_0^x E(x,x',k) (\Delta_1 \Psi)(x',k) dx',
\end{align}
where
$$E(x,x',k) = \hat{Z}(x,k) \begin{pmatrix} e^{2ik(x-x')} & 0 \\ 0 & 1 \end{pmatrix} \hat{Z}^{-1}(x',k).$$
We seek a solution $\Psi(x,k) = \sum_{l=0}^\infty \Psi_l(x,k)$ where
$$\Psi_l(x,k) = (-1)^l \int_{0 \leq x_1 \leq \cdots \leq x_l \leq x < \infty} \prod_{i = 1}^l E(x_{i+1}, x_i, k) \Delta_1(x_i, k) \Psi_0(x_1, k) dx_1 \cdots dx_l.$$ 
The estimates
\begin{align}\label{EEDeltaEboundb}
|\partial_k^jE(x, x', k)| < C(1+|x' - x|)^j, \qquad 0 \leq x' \leq x < \infty, \quad k \in \bar{\C}_+^K, \quad j = 0, 1, \dots, n,
\end{align}
and
$$|\Psi_0(x,k)| \leq C, \qquad x \geq 0, \quad k \in \bar{\C}_+^K,$$
together with  (\ref{EEDeltaest}) yield
\begin{align*}\nonumber
|\Psi_l(x,k)| \leq \; & C \int_{0 \leq x_1 \leq \cdots \leq x_l \leq x < \infty} \prod_{i = 1}^l  |\Delta_1(x_i, k)| |\Psi_0(x_1, k)| dx_1 \cdots dx_l 
	\\ \nonumber
\leq & \; \frac{C}{l!}\|\Psi_0(\cdot, k)\|_{L^\infty([0,x])} \|\Delta_1(\cdot, k)\|_{L^1([0,x])}^l
	\\ 
\leq &\; \frac{C}{l!} \bigg(\frac{Cx}{|k|^{m+1}}\bigg)^l, \qquad x \geq 0, \quad k \in \bar{\C}_+^K.
\end{align*}
Hence
\begin{align}\label{EEPsiPsi0}
|\Psi(x,k) - \Psi_0(x,k)| \leq \sum_{l=1}^\infty | \Psi_l(x,k)| \leq \frac{Cx e^{\frac{Cx}{|k|^{m+1}}}}{|k|^{m+1}}, \qquad x \geq 0, \quad k \in \bar{\C}_+^K.
\end{align}
Equations (\ref{ZYWa}) and (\ref{EEPsiPsi0}) prove the second column of (\ref{varphiasymptoticsb}) for $j = 0$.

Differentiating the integral equation (\ref{EEYVolterra2}) with respect to $k$, we find that $\Lambda := \partial_k\Psi$ satisfies
\begin{align}\label{EEintegraleq3}
\Lambda(x,k) = \Lambda_0(x,k)
+ \int_0^x E(x,x',k) \Delta_1(x', k) \Lambda(x', k) dx'
\end{align}
for each $k$ in the interior of $\bar{\C}_+^K$, where
\begin{align}\label{EELambda0}
\Lambda_0(x,k) = \;& [\partial_k\Psi_0(x,k)]_2 
+ \int_0^x \frac{\partial}{\partial k}\big[E(x,x',k) \Delta_1(x', k)\big] \Psi(x', k) dx'.
\end{align}
We seek a solution of (\ref{EEintegraleq3}) of the form $\Lambda = \sum_{l=0}^\infty \Lambda_l$. Proceeding as above, we find
$$|\Lambda_l(x,k)| \leq \frac{C}{l!}\|\Lambda_0(\cdot, k)\|_{L^\infty([0, x])}\bigg(\frac{Cx}{|k|^{m+1}}\bigg)^l, \qquad x \geq 0, \quad k \in \bar{\C}_+^K.$$
Using (\ref{Festb}) and (\ref{EEDeltaEboundb}) in (\ref{EELambda0}), we obtain
\begin{align}\label{EEFF1}
\big|\Lambda_0(x,k) - [\partial_k\Psi_0(x,k)]_2 \big| \leq \frac{C(1+x)^2}{|k|^{m+1}}, \qquad x \geq 0, \quad k \in \bar{\C}_+^K.
\end{align}
Thus $\sum_{l=0}^\infty \Lambda_l$ converges uniformly on compact subsets of $[0,\infty) \times \bar{\C}_+^K$ to a continuous solution $\Lambda$ of (\ref{EEintegraleq3}), which satisfies
\begin{align}\label{EEFF2}
|\Lambda(x,k) - \Lambda_0(x,k)| 
\leq \frac{C(1+x)^2xe^{\frac{Cx}{|k|^{m+1}}}}{|k|^{m+1}}, \qquad x \geq 0,  \quad k \in \bar{\C}_+^K.
\end{align}
Equations (\ref{ZYWa}), (\ref{EEFF1}), and (\ref{EEFF2}) show that $[Y]_2 = \Psi$ satisfies (\ref{varphiasymptoticsb}) for $j = 1$. 
Extending the above argument, we find that (\ref{varphiasymptoticsb}) holds also for $j = 2, \dots, n$.
\proofendcontinue

\medskip\noindent
{\bf Claim 9.}  $[Y]_2$ satisfies (\ref{varphiasymptoticsc}). 

{\it Proof of Claim 9.}
Let $y(x,k) = Y(x,k) e^{2ikx\sigma_3}$. Then $y$ satisfies $y_x = Ay + iky\sigma_3$. 
Thus
\begin{align*}
  (\hat{W}^{-1}y)_x & = - \hat{W}^{-1}\hat{W}_x \hat{W}^{-1} y + \hat{W}^{-1} y_x
  	\\
&  = - \hat{W}^{-1}(\hat{A}_2\hat{W} + ik\hat{W}\sigma_3) \hat{W}^{-1} y + \hat{W}^{-1} (Ay + iky\sigma_3)
	\\
&  = \hat{W}^{-1}\Delta_2 y - ik[\sigma_3, \hat{W}^{-1}y].
\end{align*}
Hence
$$\Big(e^{ikx\hat{\sigma}_3}\hat{W}^{-1}y\Big)_x = e^{ikx\hat{\sigma}_3}\hat{W}^{-1} \Delta_2 y.$$
Integrating and using the initial condition $y(0,k) = I$, we conclude that $y$ satisfies the following Volterra integral equation:
\begin{align}\label{yWhatVolterra}
y(x,k) = \hat{W}(x,k) e^{-ikx\hat{\sigma}_3}\hat{W}^{-1}(0,k) + \int_0^x \hat{W}(x,k) e^{ik(x'-x)\hat{\sigma}_3} (\hat{W}^{-1}\Delta_2 y)(x',k) dx'.
\end{align}
Letting $\Psi = [y]_2$ and $\Psi_0(x,k) = [\hat{W}(x,k) e^{-ikx\hat{\sigma}_3}\hat{W}^{-1}(0,k) ]_2$, we can write the second column of (\ref{yWhatVolterra}) as
$$\Psi(x,k) = \Psi_0(x,k) + \int_0^x E(x,x',k) (\Delta_2 \Psi)(x',k) dx',$$
where
$$E(x,x',k) = \hat{W}(x,k) \begin{pmatrix} e^{2ik(x'-x)} & 0 \\ 0 & 1 \end{pmatrix} \hat{W}^{-1}(x',k).$$
As in the proof of Claim 8, the estimates
$$|\partial_k^jE(x, x', k)| < C(1+|x' - x|)^j, \qquad 0 \leq x' \leq x < \infty, \quad k \in \bar{\C}_-, \quad j = 0, 1, \dots, n,$$
and
$$|\Psi_0(x,k)| \leq C, \qquad x \geq 0, \quad k \in \bar{\C}_-^K,$$
together with  (\ref{EEDeltaest}) yield
\begin{align}\label{HHPsiPsi0}
|\Psi(x,k) - \Psi_0(x,k)| \leq \frac{Cx e^{\frac{Cx}{|k|^{m+1}}}}{|k|^{m+1}}, \qquad x \geq 0, \quad k \in \bar{\C}_+^K.
\end{align}
Equations (\ref{ZYWb}) and (\ref{HHPsiPsi0}) prove the second column of (\ref{varphiasymptoticsc}) for $j = 0$. Proceeding as in the proof of Claim 8, equation (\ref{varphiasymptoticsc}) follows also for $j = 1, \dots, n$. 

\proofend

\begin{remark}\upshape
The proof of Theorem \ref{xth2} was inspired by Chapter 6 of \cite{CL1955}, where asymptotic results are derived for differential equations on a finite interval. 
\end{remark}

\subsection{The spectral functions $\{a(k), b(k)\}$}
Let $s(k) = X(0,k)$. 
Since $X$ obeys the symmetries (\ref{Fpmsymm}),
we may define the spectral functions $a(k)$ and $b(k)$ for $\im k \leq 0$ by 
\begin{align}\label{abdef}
s(k) = \begin{pmatrix} 
\overline{a(\bar{k})} 	&	b(k)	\\
\lambda \overline{b(\bar{k})}	&	a(k)
\end{pmatrix}, \qquad k \in (\bar{\C}_+, \bar{\C}_-).
\end{align}

\begin{theorem}\label{abth}
Suppose $u(x)$ satisfies (\ref{uassump}) for some integers $m \geq 1$ and $n \geq 1$.
Then the spectral functions  $a(k)$ and $b(k)$ have the following properties:
\begin{enumerate}[$(a)$]
\item $a(k)$ and $b(k)$ are continuous for $\im k \leq 0$ and analytic for $\im k < 0$.

\item There exist complex constants $\{a_j, b_j\}_1^m$ such that
\begin{subequations}\label{abexpansions}
\begin{align}\nonumber
& a(k) = 1 + \frac{a_1}{k} + \cdots + \frac{a_m}{k^m} + O\bigg(\frac{1}{k^{m+1}}\bigg), 
	\\ 
& b(k) = \frac{b_1}{k} + \cdots + \frac{b_m}{k^m} + O\bigg(\frac{1}{k^{m+1}}\bigg),	
\end{align}
uniformly as $k \to \infty$ with $\im k \leq 0$.

\item For $j = 1, \dots n$, the derivatives $a^{(j)}(k)$ and $b^{(j)}(k)$ have continuous extensions to $\im k \leq 0$ and 
\begin{align}\nonumber
& a^{(j)}(k) = \frac{d^j}{dk^k}\bigg(1 + \frac{a_1}{k} + \cdots + \frac{a_m}{k^m}\bigg) + O\bigg(\frac{1}{k^{m+1}}\bigg), 
	\\
& b^{(j)}(k) = \frac{d^j}{dk^j}\bigg(\frac{b_1}{k} + \cdots + \frac{b_m}{k^m}\bigg) + O\bigg(\frac{1}{k^{m+1}}\bigg), 
\end{align}
\end{subequations}
uniformly as $k \to \infty$ with $\im k \leq 0$.

\item $a$ and $b$ obey the symmetries
\begin{align}\label{absymm}
\begin{cases}
a(k) = \overline{a(-\bar{k})}, \\ 
b(k) = \overline{b(-\bar{k})},
\end{cases} \qquad \im k \leq 0.
\end{align}

\item $|a(k)|^2 - \lambda |b(k)|^2 = 1$ for $k \in \R$.

\item If $\lambda = 1$, then $a(k) \neq 0$  for $\im k \leq 0$. 

\end{enumerate}
\end{theorem}
\proofbegin 
Letting $a_j = (X_j(0))_{22}$ and $b_j = (X_j(0))_{12}$, properties $(a)$-$(c)$ follow immediately from Theorems \ref{xth1} and \ref{xth2}. 
Property $(d)$ is a consequence of the symmetry $X(x,k) = \overline{X(x,-\bar{k})}$.
Property $(e)$ follows since $\det X = 1$. 

Suppose $\lambda = 1$. Then $(e)$ implies  that $|a(k)| \geq 1$  for $k \in \R$. It only remains to prove that $a(k) \neq 0$ for $\im k < 0$.

Suppose $a(k_0) = 0$ for some $k_0 \in \C$ with  $\im k_0 < 0$. Consider the space $L^2(\R, \C^2)$ of vector valued functions $f = (f_1, f_2)$ equipped with the inner product
$$\langle f, g\rangle = \int_\R (\bar{f}_1 g_1 + \bar{f}_2 g_2) dx.$$
Let 
$$u_e(x) = \begin{cases} u(x), & x \geq 0, \\ 0, & x < 0, \end{cases} \quad \text{and} \quad \mathsf{U}_e = \begin{pmatrix} 0 & u_e \\ u_e & 0 \end{pmatrix}.$$
Then the operator $L = i\sigma_3 \partial_x - i\sigma_3\mathsf{U}_e$ satisfies
$$\langle Lf, g \rangle = \langle f, Lg\rangle \quad \text{whenever} \quad f,g \in H^1(\R, \C^2) \subset L^2(\R, \C^2).$$
Define $h \in L^2(\R, \C^2)$ by
$$h(x) = \begin{cases}
[X(x,k_0)]_2 e^{-ik_0 x}, & x \geq 0, \\
\begin{pmatrix} b(k_0) e^{ik_0 x} \\ 0 \end{pmatrix}, & x < 0.
\end{cases}$$
The condition $a(k_0) = 0$ implies that $h$ is continuous at $x = 0$. Moreover, since $\im k_0 < 0$, $h$ has exponential decay as $x \to \pm \infty$. It follows that $h \in H^1(\R, \C^2)$. But since $Lh = k_0 h$ this leads to the contradiction that the eigenvalue $k_0$ must be real:
$$\bar{k}_0 \langle h, h \rangle = \langle Lh, h \rangle =\langle h, Lh \rangle = k_0 \langle h, h \rangle.$$
This proves $(f)$. 
\proofend

%Proof that a function in H^1 goes to zero at infinity!! (The proof in master's thesis p. 26 is incomplete.)
%Theorem If $f \in H^1(\R, \C)$, then $\lim_{|x| \to \infty} f(x) = 0$. 
%Proof. By taking real and imaginary parts, we may assume that $f$ is real valued. 
%Suppose $\limsup_{x \to \infty} |f(x)| = 2c > 0$. Then there exists a sequence $x_k$ such that $x_k \to \infty$ and $|f(x_k)| > c$ for each $k$. Choose $k$ so large that $\int_{x_k}^\infty (f^2 + f_x^2) dx < \frac{c^2}{2}$. Then integration of the inequality
%$$-2f f_x\leq f^2 + f_x^2$$
%from $x_k$ to $y$ yields
%$$f^2(x_k) - f^2(y) \leq \int_{x_k}^y (f^2 + f_x^2) dx < \frac{c^2}{2}$$
%for all $y > x_k$. But then
%$$f^2(y) > f^2(x_k) - \frac{c^2}{2} > c^2 - \frac{c^2}{2} = \frac{c^2}{2}$$
%for all $y > x_k$. But this contradicts the fact that $f \in L^2(\R)$. 
%\proofend

\section{Spectral analysis of the $t$-part}\nequation\label{tsec}
Let $\{g_j(t)\}_{j=0}^2$ be real-valued functions defined for $t \geq 0$ and let
\begin{align*}
\mathsf{V}(t,k) & = -4k^2g_0\sigma_3\sigma_\lambda - 2i \lambda kg_0^2\sigma_3 + 2i kg_1\sigma_\lambda + (g_2 - 2\lambda g_0^3)\sigma_3\sigma_\lambda
	\\
& = \begin{pmatrix} -2i\lambda kg_0^2 & -4k^2 g_0 + 2ik g_1 - 2 \lambda g_0^3 + g_2 \\
-4\lambda k^2g_0 - 2i\lambda k g_1 - 2 g_0^3 + \lambda g_2 & 2i\lambda kg_0^2 \end{pmatrix}.
\end{align*}
Consider the linear differential equation
\begin{align}\label{tpartG}
T_t + 4ik^3[\sigma_3, T] = \mathsf{V} T,
\end{align}
where $T(t,k)$ is a $2 \times 2$-matrix valued eigenfunction and $k \in \C$ is a spectral parameter. We define two $2 \times 2$-matrix valued solutions of (\ref{tpartG}) as the solutions of the linear Volterra integral equations
\begin{subequations}\label{Gpmdef}
\begin{align}  \label{Gpmdefa}
 & T(t,k) = I + \int_{\infty}^t e^{4i k^3(t'-t)\hat{\sigma}_3} (\mathsf{V}T)(t',k) dt',
  	\\
&  U(t,k) = I + \int_{0}^t e^{4i k^3(t'-t)\hat{\sigma}_3} (\mathsf{V}U)(t',k) dt'.
\end{align}
\end{subequations}
The proof of the following theorem is given in Section \ref{tproof1subsec}.

\begin{theorem}\label{tth1}
Let $m \geq 1$ and $n \geq 1$ be integers. 
Suppose 
\begin{align}\label{gjassump}
\begin{cases}
g_j \in C^{[\frac{m+5 - j}{3}]}([0,\infty)), & j = 0,1,2,
	\\
(1+t)^{n}g_j^{(i)}(t) \in L^1([0,\infty)), & j = 0, 1, 2, \quad i = 0,1, \dots, [\frac{m + 5 - j}{3}].
\end{cases}
\end{align}

Then equation (\ref{tpartG}) admits two $2 \times 2$-matrix valued solutions $T$ and $U$ with the following properties:
\begin{enumerate}[$(a)$]
\item The function $T(t, k)$ is defined for $t \geq 0$ and $k \in (\bar{D}_-, \bar{D}_+)$. For each $k \in (\bar{D}_-, \bar{D}_+)$, the function $T(\cdot, k) \in C^1([0,\infty))$ satisfies (\ref{tpartG}).

\item The function $U(t, k)$ is defined for $t \geq 0$ and $k \in \C$. For each $k \in \C$,  the function $U(\cdot, k) \in C^1([0,\infty))$ satisfies (\ref{tpartG}). 

\item For each $t \geq 0$, the function $T(t,\cdot)$ is bounded and continuous for $k \in (\bar{D}_-, \bar{D}_+)$ and analytic for $k \in (D_-, D_+)$.

\item For each $t \geq 0$, the function $U(t,\cdot)$ is an entire function of $k \in \C$ which is bounded for $k \in (\bar{D}_+, \bar{D}_-)$.

\item For each $t \geq 0$ and each $j = 1, \dots, n$, the partial derivative $\frac{\partial^j T}{\partial k^j}(t, \cdot)$ has a continuous extension to $(\bar{D}_-, \bar{D}_+)$.

\item $T$ and $U$ satisfy the following estimates: 
\begin{subequations}\label{Gest}
\begin{align}\label{Gesta}
& \bigg|\frac{\partial^j}{\partial k^j}\big(T(t,k) - I\big) \bigg| \leq 
\frac{C(1+|k|)^{2 + 4j} e^{\frac{C(1+|k|)^2}{(1+t)^n}}}{(1+t)^{n-j}}, \qquad t \geq 0, \quad  k \in (\bar{D}_-, \bar{D}_+),
	\\\label{Gestb}
& \bigg|\frac{\partial^j}{\partial k^j}\Big(U(t,k) - I\Big) \bigg| \leq 
C\min(t,1)(1 + t)^j e^{C(1+|k|)^2}, \qquad t \geq 0, \quad k \in (\bar{D}_+, \bar{D}_-),
\end{align}
\end{subequations}
for $ j = 0, 1, \dots, n$.
\end{enumerate}
\end{theorem}

\subsection{Behavior as $k \to \infty$}
Equation (\ref{tpartG}) admits formal power series solutions $T_{formal}$ and $U_{formal}$, normalized at $t = \infty$ and $t = 0$ respectively, such that
\begin{align}\label{psiformalt}
& T_{formal}(t,k) = I + \frac{T_1(t)}{k} + \frac{T_2(t)}{k^2} + \cdots,
	\\
& U_{formal}(t,k) =  I + \frac{V_1(t)}{k} + \frac{V_2(t)}{k^2} + \cdots + \bigg(\frac{W_1(t)}{k} + \frac{W_2(t)}{k^2} + \cdots \bigg) e^{8ik^3t\sigma_3},	
\end{align}
and
$$\lim_{t \to \infty} T_j(t) = 0, \qquad V_j(0) + W_j(0) = 0, \qquad j \geq 1.$$ 
Indeed, substituting (\ref{psiformalt}) into (\ref{tpartG}) the diagonal terms of $O(k^{-j})$ yield
\begin{align}\label{Gjdiag1}
\partial_t T_j^{(d)} = -4g_0 \sigma_3\sigma_\lambda T_{j+2}^{(o)} - 2i\lambda g_0^2 \sigma_3 T_{j+1}^{(d)} 
+ 2i g_1 \sigma_\lambda T_{j+1}^{(o)} + (g_2 - 2 \lambda g_0^3) \sigma_3\sigma_\lambda T_j^{(o)},
\end{align}
while the off-diagonal terms of $O(k^{-j+3})$ yield
\begin{subequations}\label{Gjrecursive}
\begin{align}\nonumber
T_{j}^{(o)} = &\; \frac{i}{8}\sigma_3\Big(\partial_t T_{j-3}^{(o)} + 4g_0\sigma_3\sigma_\lambda T_{j-1}^{(d)} + 2i\lambda g_0^2 \sigma_3 T_{j-2}^{(o)} - 2ig_1 \sigma_\lambda T_{j-2}^{(d)} 
	\\ \label{Gjoff}
& - (g_2 - 2\lambda g_0^3)\sigma_3\sigma_\lambda T_{j-3}^{(d)}\Big).
\end{align} 
Using (\ref{Gjoff}) to eliminate $T_{j+2}^{(o)}$, $T_{j+1}^{(o)}$, and $T_j^{(o)}$ from (\ref{Gjdiag1}), we find after simplification
\begin{align}\nonumber
 \partial_t T_j^{(d)} = & \frac{ig_0}{2}\sigma_\lambda \partial_t T_{j-1}^{(o)} 
 + \frac{1}{4} g_1 \sigma_3\sigma_\lambda  \partial_t T_{j-2}^{(o)} 
 - \frac{i}{8} (g_2 - \lambda g_0^3)\sigma_\lambda  \partial_t T_{j-3}^{(o)}
 - \frac{i\lambda}{2}g_1 g_0^2 \sigma_\lambda T_{j-1}^{(o)}
 	\\\nonumber
& 
+ \frac{i}{2}(3g_0^4 - 2\lambda g_2g_0 + \lambda g_1^2)\sigma_3 T_{j-1}^{(d)}
- \frac{1}{4}(\lambda g_2 - g_0^3)g_0^2\sigma_3\sigma_\lambda T_{j-2}^{(o)}
	\\\label{Gjdiag2}
&+ \frac{1}{4}g_0^3g_1 T_{j-2}^{(d)} 
+ \frac{i}{8}(g_2 - \lambda g_0^3)(\lambda g_2 - 2g_0^3)\sigma_3 T_{j-3}^{(d)}.
\end{align}
\end{subequations}
Equations (\ref{Gjrecursive}) provide the recursive equations necessary to generate the $T_j$'s. The coefficients $\{V_j\}$ satisfy the equations obtained by replacing $\{T_j\}$ with $\{V_j\}$ in (\ref{Gjrecursive}). 
Moreover, substituting
$$T = \bigg(\frac{W_1(t)}{k} + \frac{W_2(t)}{k^2} + \cdots \bigg) e^{8ik^3t\sigma_3}$$
into  (\ref{tpartG}) we find that the coefficients $\{W_j\}$ satisfy the equations obtained from (\ref{Gjrecursive}) by replacing $\{T_j^{(d)}\}$ and $\{T_j^{(o)}\}$ with $\{W_j^{(o)}\}$ and $\{W_j^{(d)}\}$, respectively. 
%the off-diagonal terms of $O(k^{-j})$ yield
%\begin{align}
%\partial_t W_j^{(o)} = -4g_0 \sigma_3\sigma_\lambda W_{j+2}^{(d)} - 2i\lambda g_0^2 \sigma_3 W_{j+1}^{(o)}  + 2i g_1 \sigma_\lambda W_{j+1}^{(d)} + (g_2 - 2 \lambda g_0^3) \sigma_3\sigma_\lambda W_j^{(d)}.
%\end{align}
%while the diagonal terms of $O(k^{-j+3})$ yield
%\begin{align}\nonumber
%W_{j}^{(d)} = &\; \frac{i}{8}\sigma_3\Big(\partial_t W_{j-3}^{(d)} + 4g_0\sigma_3\sigma_\lambda W_{j-1}^{(o)} + 2i\lambda g_0^2 \sigma_3 W_{j-2}^{(d)} - 2ig_1 \sigma_\lambda W_{j-2}^{(o)} 
%	\\ 
%& - (g_2 - 2\lambda g_0^3)\sigma_3\sigma_\lambda W_{j-3}^{(o)}\Big).
%\end{align} 
%\begin{align}\nonumber
% \partial_t W_j^{(o)} = & \frac{ig_0}{2}\sigma_\lambda \partial_t W_{j-1}^{(d)}  + \frac{1}{4} g_1 \sigma_3\sigma_\lambda  \partial_t W_{j-2}^{(d)}  - \frac{i}{8} (g_2 - \lambda g_0^3)\sigma_\lambda  \partial_t W_{j-3}^{(d)} - \frac{i\lambda}{2}g_1 g_0^2 \sigma_\lambda W_{j-1}^{(d)}
% 	\\\nonumber
%& + \frac{i}{2}(3g_0^4 - 2\lambda g_2g_0 + \lambda g_1^2)\sigma_3 W_{j-1}^{(o)}
%- \frac{1}{4}(\lambda g_2 - g_0^3)g_0^2\sigma_3\sigma_\lambda W_{j-2}^{(d)}
%	\\
%&+ \frac{1}{4}g_0^3g_1 W_{j-2}^{(o)} + \frac{i}{8}(g_2 - \lambda g_0^3)(\lambda g_2 - 2g_0^3)\sigma_3 W_{j-3}^{(o)}.
%\end{align}
The coefficients $\{T_j(t), V_j(t), W_j(t)\}$ are determined recursively from the above equations and the initial assignments 
$$T_{-2} = T_{-1} = 0, \quad T_0 = I, \quad V_{-2} = V_{-1} = 0, \quad V_0 = I,
\quad W_{-2} = W_{-1} = W_0 = 0.$$ 
The first few coefficients are given by
\begin{align}\nonumber
T_1(t) = &\; \frac{i}{2} \begin{pmatrix} 0 & g_0(t) \\ -\lambda g_0(t) & 0 \end{pmatrix} 
+ \sigma_3 \frac{i\lambda}{2}\int_{\infty}^t (3\lambda g_0^4 +  g_1^2 - 2 g_0g_2) d\tau,
	\\ \nonumber
T_2(t) = &\;\frac{1}{4} (g_1+2 i g_0 (T_1)_{22}) \sigma_3\sigma_\lambda
	\\ \nonumber
& -I \int_{\infty}^t \frac{i\lambda}{4} \left[g_0 \left(-4 g_2 (T_1)_{22}+i
   g_{0t}\right)+6 \lambda g_0^4 (T_1)_{22}+2 g_1^2 (T_1)_{22} \right] d\tau,
	\\ \nonumber
T_3(t) = &\; -\frac{i}{8} \left(-4 g_0 (T_2)_{22}-\lambda g_0^3+2 i g_1
   (T_1)_{22}+g_2\right) \sigma_\lambda
	\\ \nonumber
& + \sigma_3 \int_{\infty}^t \frac{i\lambda}{8} \big[4 g_0^6 
-4\lambda  g_0^3 g_2
+g_2^2
+12 \lambda g_0^4 (T_2)_{22}
+ g_1 \left(4 g_1 (T_2)_{22}-g_{0t}\right)
	\\ \label{Gjexplicit}
&+  g_0 \left(2 i (T_1)_{22} g_{0t}+g_{1t}-8 g_2 (T_2)_{22}\right)\big] d\tau,
%    	\\ \nonumber
%T_4(t) = &\; \frac{1}{16} \left(-\lambda g_0^2 g_1 +4 g_1 (T_2)_{22} +2 i \left((\lambda g_0^3 -g_2) (T_1)_{22} +4 g_0 (T_3)_{22}\right) -g_{0t}\right) \sigma_3 \sigma_\lambda
%	\\ \nonumber
%& + I \int_{\infty}^t \frac{\lambda}{16} \big[5 \lambda g_0^3 g_{0t} -g_2 g_{0t} +4 g_0 (T_2)_{22}g_{0t} +g_1 g_{1t} -g_0 g_{2t}
%	\\\nonumber
%& -2 i \big(4 g_0^6 (T_1)_{22} -4 \lambda g_0^3  g_2 (T_1)_{22} +g_2^2 (T_1)_{22} +12 \lambda g_0^4 (T_3)_{22} + 4 g_1^2 (T_3)_{22}
%	\\ \label{Gjexplicit}
%& - g_1 g_{0t} (T_1)_{22}  - 8 g_2 g_0 (T_3)_{22} + g_{1t} g_0 (T_1)_{22} \big)  \big] d\tau,
\end{align}
\begin{align}\nonumber
V_1(t) = &\; \frac{i}{2} \begin{pmatrix} 0 & g_0(t) \\ -\lambda g_0(t) & 0 \end{pmatrix} 
+ \sigma_3 \frac{i\lambda}{2}\int_0^t (3\lambda g_0^4 +  g_1^2 - 2 g_0g_2) d\tau,
	\\ \nonumber
V_2(t) = &\;\frac{1}{4} (g_1+2 i g_0 (V_1)_{22}) \sigma_3\sigma_\lambda
+I\bigg\{-\int_0^t \frac{i\lambda}{4} \big[g_0 \left(-4 g_2 (V_1)_{22}+i
   g_{0t}\right)
	\\ \nonumber
& +6 \lambda g_0^4 (V_1)_{22}+2 g_1^2 (V_1)_{22} \big] d\tau + \frac{\lambda g_0^2(0)}{4}\bigg\},
	\\ \nonumber
V_3(t) = &\; -\frac{i}{8} \left(-4 g_0 (V_2)_{22}-\lambda g_0^3+2 i g_1
   (V_1)_{22}+g_2\right) \sigma_\lambda
	\\ \nonumber
& + \sigma_3 \int_0^t \frac{i\lambda}{8} \big[4 g_0^6 
-4\lambda  g_0^3 g_2
+g_2^2
+12 \lambda g_0^4 (V_2)_{22}
+ g_1 \left(4 g_1 (V_2)_{22}-g_{0t}\right)
	\\ \label{Vjexplicit}
&+  g_0 \left(2 i (V_1)_{22} g_{0t}+g_{1t}-8 g_2 (V_2)_{22}\right)\big] d\tau,
\end{align}
and
\begin{align}\nonumber
W_1(t) = &\; -\frac{ig_0(0)}{2}\sigma_\lambda,
	\\ \nonumber
W_2(t) = &\; - \frac{\lambda i}{2}g_0 (W_1)_{12} I	
+ \sigma_3 \sigma_\lambda \bigg\{\int_0^t \frac{i}{2} \big[3g_0^4 + \lambda g_1^2 - 2\lambda g_0g_2\big](W_1)_{12} d\tau - \frac{g_1(0)}{4}\bigg\},
	\\ \nonumber
W_3(t) = &\; -\frac{\lambda}{4}(g_1 (W_1)_{12} - 2ig_0(W_2)_{12})\sigma_3 
+ \sigma_\lambda \bigg\{\int_0^t \frac{1}{4} \big[2i(3g_0^4 + \lambda g_1^2 - 2\lambda g_0g_2)(W_2)_{12}
	\\ \nonumber
&  + \lambda g_0(W_1)_{12} g_{0t} \big] d\tau + \frac{i}{8}(g_2(0) - 2\lambda g_0^3(0))\bigg\}.
\end{align}
If $\{g_j(t)\}_0^2$ have a finite degree of regularity and decay, only finitely many coefficients $\{T_j, V_j, W_j\}$ are well-defined.
The following result, whose proof is given in Section  \ref{tproof2subsec}, describes the behavior of $T$ and $U$ as $k \to \infty$.

\begin{theorem}\label{tth2}
Suppose $\{g_j(t)\}_0^2$ satisfy (\ref{gjassump}) for some integers $m \geq 1$ and $n \geq 1$.
Then, as $k \to \infty$, $T$ and $U$ coincide to order $m$ with $T_{formal}$ and $U_{formal}$ respectively in the following sense: The functions
\begin{align}\label{Ghatdef}
& \hat{T}(t,k) = I + \frac{T_1(t)}{k} + \cdots + \frac{T_{m+3}(t)}{k^{m+3}}
	\\
& \hat{U}(t,k) =  I + \frac{V_1(t)}{k} + \cdots + \frac{V_{m+3}(t)}{k^{m+3}} + \bigg(\frac{W_1(t)}{k} + \cdots + \frac{W_{m+3}(t)}{k^{m+3}}\bigg) e^{8ik^3t\sigma_3},	
\end{align}
are well-defined and there exists a $K > 0$ such that
\begin{subequations}\label{Gasymptotics}
\begin{align}\label{Gasymptoticsa}
\bigg|\frac{\partial^j}{\partial k^j}\big(T - \hat{T}\big) \bigg| \leq 
\frac{C}{|k|^{m+1 - 2j}(1+ t)^{n-j}}, \qquad t \geq 0, \quad  j = 0, 1, \dots, n,
\end{align}
for all $k \in (\bar{D}_-, \bar{D}_+)$ with  $|k| > K$, and
\begin{align}\label{Gasymptoticsb}
\bigg|\frac{\partial^j}{\partial k^j}\big(U - \hat{U}\big) \bigg| \leq 
\frac{C(1 + t)^{j+2}e^{\frac{Ct}{|k|^{m+1}}}}{|k|^{m+1 -2j}}, \qquad t \geq 0, \quad j = 0, 1, \dots, n,
\end{align}
for all $k \in (\bar{D}_+, \bar{D}_-)$ with  $|k| > K$,
\begin{align}\label{Gasymptoticsc}
\bigg|\frac{\partial^j}{\partial k^j}\big((U - \hat{U})e^{-8ik^3t\sigma_3}\big) \bigg| \leq 
\frac{C(1 + t)^{j+2}e^{\frac{Ct}{|k|^{m+1}}}}{|k|^{m+1 -2j}}, \qquad t \geq 0, \quad j = 0, 1, \dots, n,
\end{align}
for all $k \in (\bar{D}_-, \bar{D}_+)$ with  $|k| > K$.
\end{subequations}
\end{theorem}

\subsection{Proof of Theorem \ref{tth1}}\label{tproof1subsec} 
In view of the symmetries 
\begin{align}\label{Gpmsymm}
F(t,k) = \begin{cases} 
\sigma_1\overline{F(t, \bar{k})} \sigma_1, & \lambda = 1, \\
\sigma_2\overline{F(t, \bar{k})} \sigma_2, & \lambda = -1,
\end{cases}
\end{align}
which are valid for $F = T$ and $F = U$, it is enough to prove the theorem for $[T]_2$ and $[U]_2$. We first consider the construction of $[T]_2$.

Let $\Psi(t,k)$ denote the second column of $T(t,k)$. As suggested by (\ref{Gpmdefa}), we define $\Psi$ by the integral equation
\begin{align}\label{ttPsiVolterra}
\Psi(t,k) = \begin{pmatrix} 0 \\ 1 \end{pmatrix} - \int_t^\infty E(t,t',k)\mathsf{V}(t', k)\Psi(t',k) dt', \qquad t \geq 0,
\end{align}
where
$$E(t,t',k) =  \begin{pmatrix} e^{8ik^3(t' - t)} & 0 \\ 0 & 1 \end{pmatrix}.$$
We  use successive approximations to show that the Volterra integral equation (\ref{ttPsiVolterra}) has a unique solution $\Psi(t,k)$ for each $k \in \bar{D}_+$.
Let $\Psi_0 = \begin{pmatrix} 0 \\ 1 \end{pmatrix}$ and define $\Psi_l$ for $l \geq 1$ inductively by 
\begin{align*}
\Psi_{l+1}(t,k) = - \int_t^\infty E(t,t',k) \mathsf{V}(t', k) \Psi_l(t', k) dt', \qquad t \geq 0, \quad k \in \bar{D}_+.
\end{align*}
Then
\begin{align}\label{ttPhiliterated}
\Psi_l(t,k) = & (-1)^l \int_{t = t_{l+1} \leq t_l \leq \cdots \leq t_1 < \infty} \prod_{i = 1}^l E(t_{i+1}, t_i, k) \mathsf{V}(t_i, k) \Psi_0 dt_1 \cdots dt_l. 
\end{align}
Using the estimates 
\begin{align*}
  |E(t, t', k)| < C, \qquad 0 \leq t \leq t' < \infty, \quad k \in \bar{D}_+,
\end{align*}
and
$$\|\mathsf{V}(\cdot, k)\|_{L^1([t,\infty))} < \frac{C(1+|k|)^2}{(1+t)^n}, \qquad t \geq 0,\quad k \in \C,$$
we find
\begin{align}\nonumber
|\Psi_l(t,k)| \leq \; & C \int_{t \leq t_l \leq \cdots \leq t_1 < \infty} \prod_{i = 1}^l  |\mathsf{V}(t_i, k)| |\Psi_0| dt_1 \cdots dt_l 
	\\ \label{ttPhilestimate}
\leq & \; \frac{C}{l!} \|\mathsf{V}(\cdot, k)\|_{L^1([t,\infty))}^l 
\leq \frac{C}{l!} \bigg(\frac{C(1+|k|)^2}{(1+t)^n}\bigg)^l, \qquad t \geq 0, \quad k \in \bar{D}_+.
\end{align}
Hence the series
$$\Psi(t,k) = \sum_{l=0}^\infty \Psi_l(t,k)$$
converges absolutely and uniformly for $t \geq 0$ and $k$ in compact subsets of $\bar{D}_+$ to a continuous solution $\Psi(t,k)$ of (\ref{ttPsiVolterra}). 
Moreover,
\begin{align}\label{ttPsiminusf}
|\Psi(t,k) - \Psi_0| \leq \sum_{l=1}^\infty | \Psi_l(t,k)| \leq \frac{C(1+|k|)^2e^{\frac{C(1+|k|)^2}{(1+t)^n}}}{(1+t)^n}, \qquad t \geq 0, \quad k \in \bar{D}_+,
\end{align}
which proves (\ref{Gesta}) for $j = 0$. Using the estimate
\begin{align}\nonumber
& |\partial_k^jE(t, t', k)| < C(1 + |k|)^{2j} (1+|t' - t|)^j, 
	\\ \label{ttE1bound}
& \hspace{5cm} 0 \leq t \leq t' < \infty, \quad k \in \bar{D}_+, \quad j = 0, 1, \dots, n,
\end{align}
with $j = 1$ to differentiate under the integral sign in (\ref{ttPhiliterated}), we see that $\Psi_l(t, \cdot)$ is analytic in $D_+$ for each $l$; the uniform convergence then proves that $\Psi$  is analytic in $D_+$. 

It remains to show that $[T]_2 = \Psi$ satisfies $(e)$ and $(f)$ for $j = 1, \dots, n$.
Let 
\begin{align}\label{ttLambda0}
\Lambda_0(t,k) = \;& - \int_t^\infty (\partial_kE)(t,t',k) \mathsf{V}(t',k) \Psi(t', k) dt'.
\end{align}
Differentiating the integral equation (\ref{ttPsiVolterra}) with respect to $k$, we find that $\Lambda := \partial_k\Psi$ satisfies
\begin{align}\label{ttintegraleq3}
\Lambda(t,k) = \Lambda_0(t,k) - \int_t^\infty E(t,t',k) \mathsf{V}(t',k) \Lambda(t', k) dt'
\end{align}
for each $k$ in the interior of $\bar{D}_+$; the differentiation can be justified by dominated convergence using (\ref{ttE1bound}) and a Cauchy estimate for $\partial_k\Psi$.
We seek a solution of (\ref{ttintegraleq3}) of the form $\Lambda = \sum_{l=0}^\infty \Lambda_l$. Proceeding as in (\ref{ttPhilestimate}), we find
\begin{align}\label{GG0}
|\Lambda_l(t,k)| \leq \frac{C}{l!}\|\Lambda_0(\cdot, k)\|_{L^\infty([t,\infty))}\bigg(\frac{C(1 + |k|)^2}{(1+t)^n}\bigg)^l, \qquad t \geq 0, \quad k \in \bar{D}_+.
\end{align}
Using (\ref{ttPsiminusf}) and (\ref{ttE1bound}) in (\ref{ttLambda0}), we obtain
\begin{align}\label{GG1}
|\Lambda_0(t,k)| \leq \frac{C(1+|k|)^6 e^{\frac{C(1+|k|)^2}{(1+t)^n}}}{(1+t)^{n-1}}, \qquad t \geq 0, \quad k \in \bar{D}_+.
\end{align}
In particular $\|\Lambda_0(\cdot, k)\|_{L^\infty([t,\infty))}$ is bounded for $t \geq 0$ and $k$ in compact subsets of $\bar{D}_+$.
Thus, $\sum_{l=0}^\infty \Lambda_l$ converges uniformly for $t \geq 0$ and $k$ in compact subsets of $\bar{D}_+$ to a continuous solution $\Lambda$ of (\ref{ttintegraleq3}), which satisfies the following analog of (\ref{ttPsiminusf}):
\begin{align}\label{GG2}
|\Lambda(t,k) - \Lambda_0(t,k)| \leq 
\frac{C(1+|k|)^6e^{\frac{C(1+|k|)^2}{(1+t)^n}}}{(1+t)^{n-1}}, \qquad t \geq 0,  \quad k \in \bar{D}_+.
\end{align}
In view of equations (\ref{GG1}) and (\ref{GG2}), we conclude that $[X]_2 = \Psi$ satisfies $(e)$ and $(f)$ for $j = 1$.

Proceeding inductively, we find that $\Lambda^{(j)} := \partial_k^j\Psi$ satisfies an integral equation of the form
\begin{align*}
\Lambda^{(j)}(t,k) = \Lambda^{(j)}_0(t,k)
- \int_t^\infty E(t,t',k) \mathsf{V}(t', k) \Lambda^{(j)}(t', k) dt', 
\end{align*}
where
$$\big|\Lambda^{(j)}_0(t, k)\big| \leq 
\frac{C(1+|k|)^{2+4j} e^{\frac{C(1+|k|)^2}{(1+t)^n}}}{(1+t)^{n-j}}, \qquad t \geq 0, \quad k \in \bar{D}_+.$$
If $1 \leq j \leq n$, then $\|\Lambda^{(j)}_0(\cdot, k)\|_{L^\infty([t,\infty))}$ is bounded for $t \geq 0$ and $k$ in compact subsets of $\bar{D}_+$; hence the associated series $ \Lambda^{(j)} = \sum_{l=0}^\infty \Lambda^{(j)}_l$ converges uniformly on compact subsets of $[0,\infty) \times \bar{D}_+$ to a continuous solution with the desired properties.

We now consider the the construction of $[U]_2$. In this case, we introduce $\{\Psi_l\}_1^\infty$ by the following analog of (\ref{ttPhiliterated}):
\begin{align*}
\Psi_l(t,k) = & (-1)^l \int_{0 \leq t_1 \leq  \cdots \leq t_l \leq t < \infty} \prod_{i = 1}^l E(t_{i+1}, t_i, k) \mathsf{V}(t_i, k) \Psi_0 dt_1 \cdots dt_l,
\end{align*}
where $t \geq 0$ and $k \in \C$. This leads to an entire function $\Psi(t,k) = \sum_{l=0}^\infty \Psi_l(t,k)$ satisfying (\ref{tpartG}). Moreover, as in (\ref{ttPhilestimate}), we find
\begin{align}\nonumber
|\Psi_l(t,k)| \leq \frac{C}{l!} \|\mathsf{V}(\cdot, k)\|_{L^1([0,t])}^l 
\leq \frac{C\min(t,1)^l(1+|k|)^{2l}}{l!}, \qquad t \geq 0, \quad k \in \bar{D}_-,
\end{align}
which leads to the following analog of (\ref{tPsiminusf}):
\begin{align*}
|\Psi(t,k) - \Psi_0| \leq \sum_{l=1}^\infty | \Psi_l(t,k)| \leq C\min(t,1)e^{C(1+|k|)^{2}}, \qquad t \geq 0, \quad k \in \bar{D}_-.
\end{align*}
This proves (\ref{Gestb}) for $j = 0$. Letting
\begin{align*}
\Lambda_0(t,k) = \;& \int_0^t (\partial_kE)(t,t',k) \mathsf{V}(t',k) \Psi(t', k) dt',
\end{align*}
we find that $\Lambda := \partial_k\Psi$ satisfies
\begin{align*}
\Lambda(t,k) = \Lambda_0(t,k) + \int_0^t E(t,t',k) \mathsf{V}(t',k) \Lambda(t', k) dt'
\end{align*}
The proof of (\ref{Gestb}) for $j = 1$ follows from the following analogs of (\ref{GG0})-(\ref{GG2}):
\begin{align*}
& |\Lambda_l(t,k)| \leq \frac{C}{l!}\|\Lambda_0(\cdot, k)\|_{L^\infty([0,t])}(1 + |k|)^{2l}, \qquad 
	\\ 
&  |\Lambda_0(t,k)| \leq Ct(1+|k|)^6 e^{C(1+|k|)^2}, 
	\\ \nonumber
& |\Lambda(t,k) - \Lambda_0(t,k)| \leq Ct(1+|k|)^6e^{C(1+|k|)^2},
\end{align*}
which are valid for $t \geq 0$ and $k \in \bar{D}_-$; the proof for $j = 2, \dots, n$ is similar.
\proofend

\subsection{Proof of Theorem \ref{tth2}}\label{tproof2subsec}
We first consider $[T]_2$. Our first goal is to show that $\hat{T}$ is well-defined and invertible for $k$ large enough. 

\medskip\noindent
{\bf Claim 1.}  $\{T_j(t)\}_1^{m+3}$ are $C^1$ functions of $t \geq 0$ satisfying 
\begin{align}
\begin{cases}
(1 + t)^n T_j(t) \in L^1([0,\infty)) \cap L^\infty([0,\infty)), 
	\\
(1 + t)^n T_j'(t) \in L^1([0,\infty)), 
\end{cases} \qquad j = 1, \dots, m+3.
\end{align}

{\it Proof of Claim 1.}
The assumption (\ref{gjassump}) implies that
\begin{align}\label{gjibound}
g_j^{(i)}(t) \leq \frac{C}{(1 + t)^n}, \qquad t \geq 0, \quad j = 0,1,2, \quad i = 0, 1, \dots, \bigg[\frac{m + 5 - j}{3}\bigg] - 1.
\end{align}
Indeed, if $i = 0, 1, \dots, [\frac{m + 5 - j}{3}] - 1$ and $t \geq 0$, then
\begin{align}\label{1tgiest}
\big|(1 + t)^n g_j^{(i)}(t)\big| \leq |g_j^{(i)}(0)| + \int_0^t \big|n(1 + t')^{n-1}g_j^{(i)}(t') + (1+t')^n g_j^{(i+1)}(t')\big|dt' < C.
\end{align}

Let $S_j$ refer to the statement
$$
\begin{cases}
T_j^{(o)} \in C^{[\frac{m+6-j}{3}]}([0,\infty)), \qquad
T_j^{(d)} \in C^{[\frac{m+7-j}{3}]}([0,\infty)), 
	\\
\text{$(1 + t)^n \partial^i T_j^{(o)}(t) \in L^1([0,\infty))$ for $i = 0, 1, \dots, [\frac{m+6-j}{3}]$,}
	\\
\text{$(1 + t)^n \partial^i T_j^{(d)}(t) \in L^1([0,\infty))$ for $i = 0, 1, \dots, [\frac{m+7-j}{3}]$.}
\end{cases}$$
The same kind of argument that led to (\ref{xF1dest}) now shows that
\begin{align}\label{tG1dest}
& \|(1 + t)^n T_j^{(d)}\|_{L^1([0,\infty))} < \infty, \qquad j = 1, 2, 3.
\end{align}
Using (\ref{tG1dest}) and the explicit expressions in (\ref{Gjexplicit}), we conclude that $S_1$, $S_2$, and $S_3$ hold.
Estimates similar to (\ref{xF1dest}) together with the relations (\ref{Gjrecursive}) imply that if $4 \leq j \leq m+3$ and $S_{j-1}, S_{j-2}, S_{j-3}$ hold, then $S_{j}$ also holds. Thus, by induction, $\{S_j\}_{j=1}^{m+3}$ hold. This shows that $\{T_j\}_1^{m+3}$ are $C^1$ functions satisfying $(1+t)^n\partial^iT_j(t) \in L^1([0,\infty))$ for $i =0,1$ and $j = 1, \dots, m+3$. The boundedness of $(1+t)^nT_j(t)$ follows by an estimate analogous to (\ref{1tgiest}). 
\proofendcontinue

\medskip\noindent
{\bf Claim 2.}  There exists a $K> 0$ such that $\hat{T}(t,k)^{-1}$ exists for all $k \in \C$ with $|k| \geq K$. 
Moreover, letting $A = -4ik^3\sigma_3 + \mathsf{V}$ and 
\begin{align}\label{tAhatdef}
\hat{A}(t,k) = \big(\hat{T}_t(t,k) - 4ik^3 \hat{T}(t,k) \sigma_3\big)\hat{T}(t,k)^{-1}, \qquad t \geq 0, \quad |k| \geq K,
\end{align}
the difference $\Delta(t,k) = A(t,k) - \hat{A}(t,k)$ satisfies
\begin{align}\label{tpreABsmall}
|\Delta(t,k)| \leq \frac{Cf(t)}{|k|^{m+1}(1+t)^n}, \qquad t \geq 0, \quad |k| \geq K,
\end{align}
where $f$ is a function in $L^1([0,\infty)) \cap C([0,\infty))$. 
In particular,
\begin{align}\label{tABsmall}
\|\Delta(\cdot, k)\|_{L^1([t,\infty))} \leq \frac{C}{|k|^{m+1}(1+t)^n}, \qquad t \geq 0, \quad |k| \geq K.
\end{align}

{\it Proof of Claim 2.}
By Claim 1, there exists a bounded continuous function $g \in L^1([0,\infty))$ such that
\begin{align}\label{Gjtbound}
|T_j(t)| \leq \frac{g(t)}{(m+3)(1+t)^n}, \qquad t \geq 0, \quad j = 1, \dots, m+3.
\end{align}
In particular,
$$\bigg| \sum_{j=1}^{m+3} \frac{T_j(t)}{k^j}\bigg| \leq \frac{g(t)}{|k|(1+t)^n}, \qquad t \geq 0, \quad |k| \geq 1.$$
Choose $K > \max(1, \|g\|_{L^\infty([0,\infty))})$. Then $\hat{T}(t,k)^{-1}$ exists whenever $|k| \geq K$ and is given by the absolutely and uniformly convergent Neumann series
$$\hat{T}(t,k)^{-1} = \sum_{l =0}^\infty \bigg(-\sum_{j=1}^{m+3} \frac{T_j(t)}{k^j}\bigg)^l, \qquad t \geq 0, \quad |k| \geq K.$$
Furthermore,
\begin{align}\label{ttailestimate}
\bigg|\sum_{l = m+4}^\infty \bigg(-\sum_{j=1}^{m+3} \frac{T_j(t)}{k^j}\bigg)^l\bigg|
\leq \sum_{l = m+4}^\infty \bigg(\frac{g(t)}{|k|(1+t)^n}\bigg)^l
\leq \frac{Cg(t)}{|k|^{m+4}(1+t)^n},
\end{align}
for $t \geq 0$ and $|k| \geq K$. Now let $Q_0(t) + \frac{Q_1(t)}{k} + \frac{Q_2(t)}{k^2} + \cdots$ be the formal power series expansion of $\hat{T}(t,k)^{-1}$ as $k \to \infty$, i.e.
\begin{align*}
& Q_0(t) = I, \quad Q_1(t) = -T_1(t), \quad Q_2(t) = T_1(t)^2 - T_2(t), 
	\\
& Q_3(t) = T_1(t)T_2(t) + T_2(t) T_1(t) - T_1^3(t) - T_3(t), \quad \dots.
\end{align*}
Equation (\ref{Gjtbound}) and the inequality (\ref{ttailestimate}) imply that the function $\mathcal{E}(t,k)$ defined by
$$\mathcal{E}(t,k) = \hat{T}(t,k)^{-1} - \sum_{j=0}^{m+3} \frac{Q_j(t)}{k^j}$$
satisfies
\begin{align}\label{tEestimate}
|\mathcal{E}(t,k)| \leq \frac{Cg(t)}{|k|^{m+4}(1+t)^n}, \qquad t \geq 0, \quad |k| \geq K.
\end{align}

Let $\hat{A}(t,k)$ be given by (\ref{tAhatdef}).
Since $T_{formal}$ is a formal solution of (\ref{tpartG}), the coefficient of $k^{-j}$ in the formal expansion of $\Delta = A - \hat{A}$ as $k \to \infty$  vanishes for $j \leq m$; hence, in view of Claim 1 and (\ref{tEestimate}),
\begin{align}\nonumber
|\Delta| 
& = \bigg|A - (\hat{T}_t - 4ik^3 \hat{T} \sigma_3)\bigg(\sum_{j=0}^{m+3} \frac{Q_j}{k^j} + \mathcal{E}\bigg)\bigg|
	\\\nonumber
& \leq \frac{Cg(t)}{|k|^{m+1}(1+t)^n} + |(\hat{T}_t - 4ik^3 \hat{T} \sigma_3)\mathcal{E}|
	\\ \nonumber
& \leq \frac{Cf(t)}{|k|^{m+1}(1+t)^n}, \qquad t \geq 0, \quad |k| \geq K,
\end{align}
where $f$ is a continuous (not necessarily bounded) function in $L^1([0,\infty))$. This proves (\ref{tpreABsmall}).
\proofendcontinue

Given $K > 0$, we let $\bar{D}_\pm^K = \bar{D}_\pm \cap \{|k| \geq K\}$.

\medskip\noindent
{\bf Claim 3.} The Volterra integral equation
\begin{align} \label{tintegraleq}
\Psi(t,k) = &\; \Psi_0(t,k) - \int_t^\infty E(t,t',k) \Delta(t', k) \Psi(t', k) dt',
\end{align}
where $\Psi_0(t,k) =  [\hat{T}(t,k)]_2$ and
\begin{align*}
& E(t,t',k) =  \hat{T}(t,k) \begin{pmatrix} e^{8ik^3(t' - t)} & 0 \\ 0 & 1 \end{pmatrix} \hat{T}(t',k)^{-1},
\end{align*}
has a unique solution $\Psi(t,k)$ for each $k \in \bar{D}_+^K$. This solution satisfies $\Psi = [T]_2$.

{\it Proof of Claim 3.}
Define $\Psi_l$ for $l \geq 1$ inductively by 
\begin{align*}
\Psi_{l+1}(t,k) = - \int_t^\infty E(t,t',k) \Delta(t', k) \Psi_l(t', k) dt', \qquad t \geq 0, \quad k \in \bar{D}_+^K.
\end{align*}
Then
\begin{align}\label{tPhiliterated}
\Psi_l(t,k) = & (-1)^l \int_{t = t_{l+1} \leq t_l \leq \cdots \leq t_1 < \infty} \prod_{i = 1}^l E(t_{i+1}, t_i, k) \Delta(t_i, k) \Psi_0(t_1, k) dt_1 \cdots dt_l. 
\end{align}
Using the estimate 
\begin{align*}
  |E(t, t', k)| < C, \qquad 0 \leq t \leq t' < \infty, \quad k \in \bar{D}_+^K,
\end{align*}
as well as (\ref{tABsmall}), we find, for $t \geq 0$ and $k \in \bar{D}_+^K$,
\begin{align}\nonumber
|\Psi_l(t,k)| \leq \; & C \int_{t \leq t_l \leq \cdots \leq t_1 < \infty} \prod_{i = 1}^l  |\Delta(t_i, k)| |\Psi_0(t_1, k)| dt_1 \cdots dt_l 
	\\ \nonumber
\leq & \; \frac{C}{l!}\|\Psi_0(\cdot, k)\|_{L^\infty([t,\infty))} \|\Delta(\cdot, k)\|_{L^1([t,\infty))}^l
	\\ \label{tPhilestimate}
\leq & \; \frac{C}{l!} \|\Psi_0(\cdot, k)\|_{L^\infty([t,\infty))}\bigg(\frac{C}{|k|^{m+1}(1+t)^n}\bigg)^l.
\end{align}
Since
$$\sup_{k \in \bar{D}_+^K} \|\Psi_0(\cdot, k)\|_{L^\infty([0,\infty))} < C,$$ 
we find
$$|\Psi_l(t,k)| \leq \frac{C}{l!}\bigg(\frac{C}{|k|^{m+1}(1+t)^n}\bigg)^l, \qquad t \geq 0, \quad k \in \bar{D}_+^K.$$
Hence the series $\Psi(t,k) = \sum_{l=0}^\infty \Psi_l(t,k)$ converges absolutely and uniformly for $t \geq 0$ and $k \in \bar{D}_+^K$ to a continuous solution $\Psi(t,k)$ of (\ref{tintegraleq}).
Moreover,
\begin{align}\label{tPsiminusf}
|\Psi(t,k) - \Psi_0(t,k)| \leq \sum_{l=1}^\infty | \Psi_l(t,k)| \leq \frac{C}{|k|^{m+1}(1+t)^n}, \qquad t \geq 0, \quad k \in \bar{D}_+^K.
\end{align}
It follows from the integral equation (\ref{tintegraleq}) that $\Psi$ satisfies the second column of (\ref{tpartG}). By (\ref{tPsiminusf}), $\Psi(t,k) \sim \begin{pmatrix} 0 \\ 1 \end{pmatrix}$ as $t \to \infty$. Hence, by uniqueness of solution, $\Psi = [T]_2$ for $k \in \bar{D}_+^K$.
\proofendcontinue

\medskip\noindent
{\bf Claim 4.}  $[T]_2$ satisfies (\ref{Gasymptoticsa}).

{\it Proof of Claim 4.}
Equation (\ref{tPsiminusf}) implies that $[T]_2 = \Psi$ satisfies (\ref{Gasymptoticsa}) for $j = 0$.

A glance at the proof of Claim 2 shows that the inequality (\ref{tpreABsmall}) can be extended to derivatives of $\Delta$ with respect to $k$:
\begin{subequations}\label{tDeltaEbound}
\begin{align}\label{tDeltaEbounda}
\big|\partial_k^j\Delta(t,k)\big| \leq \frac{Cf(t)}{|k|^{m+1}(1 + t)^n}, \qquad t \geq 0, \quad |k| \geq K, \quad j = 0,1, \dots, n,
\end{align}
where $f$ is a function in $L^1([0,\infty)) \cap C([0,\infty))$. We will also need the following estimate for $j = 0,1, \dots, n$:
\begin{align}\label{tDeltaEboundb}
|\partial_k^jE(t, t', k)| < C|k|^{2j}(1+|t' - t|)^j, \qquad 0 \leq t \leq t' < \infty, \quad k \in \bar{D}_+^K.
\end{align}
\end{subequations}
Let 
\begin{align}\label{tLambda0}
\Lambda_0(t,k) = \;& [\partial_k\hat{T}(t,k)]_2 
- \int_t^\infty \frac{\partial}{\partial k}\big[E(t,t',k) \Delta(t', k)\big] \Psi(t', k) dt'.
\end{align}
Differentiating the integral equation (\ref{tintegraleq}) with respect to $k$, we find that $\Lambda := \partial_k\Psi$ satisfies
\begin{align}\label{tintegraleq3}
\Lambda(t,k) = \Lambda_0(t,k)
- \int_t^\infty E(t,t',k) \Delta(t', k) \Lambda(t', k) dt'
\end{align}
for each $k$ in the interior of $\bar{D}_+^K$; the differentiation can be justified by dominated convergence using (\ref{tDeltaEbound}) and a Cauchy estimate for $\partial_k\Psi$.
We seek a solution of (\ref{tintegraleq3}) of the form $\Lambda = \sum_{l=0}^\infty \Lambda_l$. Proceeding as in (\ref{tPhilestimate}), we find
$$|\Lambda_l(t,k)| \leq \frac{C}{l!}\|\Lambda_0(\cdot, k)\|_{L^\infty([t,\infty))}\bigg(\frac{C}{|k|^{m+1}(1+t)^n}\bigg)^l, \qquad t \geq 0, \quad k \in \bar{D}_+^K.$$
Using (\ref{tPsiminusf}) and (\ref{tDeltaEbound}) in (\ref{tLambda0}), we obtain
\begin{align}\label{HHLambda0}
|\Lambda_0(t,k) - [\partial_k\hat{T}(t,k)]_2 | \leq \frac{C}{|k|^{m-1}(1+t)^{n-1}}, \qquad t \geq 0, \quad k \in \bar{D}_+^K.
\end{align}
In particular $\|\Lambda_0(\cdot, k)\|_{L^\infty([t,\infty))}$ is bounded for $k \in \bar{D}_+^K$ and $t \geq 0$.
Thus, $\sum_{l=0}^\infty \Lambda_l$ converges uniformly on $[0,\infty) \times \bar{D}_+^K$ to a continuous solution $\Lambda$ of (\ref{tintegraleq3}), which satisfies the following analog of (\ref{tPsiminusf}):
\begin{align}\label{HHLambdaLambda0}
|\Lambda(t,k) - \Lambda_0(t,k)| \leq 
\frac{C}{|k|^{m+1}(1+t)^{n}}, \qquad t \geq 0,  \quad k \in \bar{D}_+^K.
\end{align}
Equations (\ref{HHLambda0}) and (\ref{HHLambdaLambda0}) show that $[T]_2 = \Psi$ satisfies (\ref{Gasymptoticsa}) for $j = 1$.

Proceeding inductively, we find that $\Lambda^{(j)} := \partial_k^j\Psi$ satisfies an integral equation of the form
\begin{align*}
\Lambda^{(j)}(t,k) = \Lambda^{(j)}_0(t,k)
- \int_t^\infty E(t,t',k) \Delta(t', k) \Lambda^{(j)}(t', k) dt', 
\end{align*}
where
$$\big|\Lambda^{(j)}_0(t, k) - [\partial_k^{j} \hat{T}(t, k)]_2 \big| \leq 
\frac{C}{|k|^{m+1 - 2j}(1+t)^{n-j}}, \qquad t \geq 0, \quad k \in \bar{D}_+^K.$$
If $1 \leq j \leq n$, then $\|\Lambda^{(j)}_0(\cdot, k)\|_{L^\infty([t,\infty))}$ is bounded for $k$ in compact subsets of $\bar{D}_+^K$ and $t \geq 0$; hence the associated series $ \Lambda^{(j)} = \sum_{l=0}^\infty \Lambda^{(j)}_l$ converges uniformly on compact subsets of $[0,\infty) \times \bar{D}_+^K$ to a continuous solution with the desired properties.

The above claims prove the theorem for $T$. We now consider $[U]_2$.

\medskip\noindent
{\bf Claim 5.}  $\{V_j(t), W_j(t)\}_1^{m+3}$ are $C^1$ functions of $t \geq 0$ satisfying
\begin{align*}
\begin{cases}
Z_j, W_j \in L^\infty([0,\infty)), \\
Z_j', W_j' \in L^1([0,\infty)),
\end{cases} \qquad j = 1, \dots, m+3.
\end{align*}

{\it Proof of Claim 5.}
Let $S_j$ and $\mathcal{S}_j$ refer to the statements
$$
\begin{cases}
V_j^{(o)} \in C^{[\frac{m+6-j}{3}]}([0,\infty)), \qquad
V_j^{(d)} \in C^{[\frac{m+7-j}{3}]}([0,\infty)), 
	\\
\text{$\partial^i V_j^{(o)}(t) \in L^1([0,\infty))$ for $i = 1, \dots, [\frac{m+6-j}{3}]$,}
	\\
\text{$\partial^i V_j^{(d)}(t) \in L^1([0,\infty))$ for $i = 1, \dots, [\frac{m+7-j}{3}]$,}
\end{cases}$$
and
$$
\begin{cases}
W_j^{(d)} \in C^{[\frac{m+6-j}{3}]}([0,\infty)), \qquad
W_j^{(o)} \in C^{[\frac{m+7-j}{3}]}([0,\infty)), 
	\\
\text{$\partial^i W_j^{(d)}(t) \in L^1([0,\infty))$ for $i = 1, \dots, [\frac{m+6-j}{3}]$,}
	\\
\text{$\partial^i W_j^{(o)}(t) \in L^1([0,\infty))$ for $i = 1, \dots, [\frac{m+7-j}{3}]$,}
\end{cases}$$
respectively.
Using the expressions (\ref{Fjmexplicit}) and (\ref{calFjmexplicit}) for $V_1$ and $W_1$, we conclude that $\{S_j, \mathcal{S}_j\}_1^3$ hold. The relations (\ref{xrecursive}) and (\ref{xrecursive2}) imply that if $4 \leq j \leq m+3$ and $\{S_i, \mathcal{S}_i\}_{i=j-3}^{j-1}$ hold, then $S_j$ and $\mathcal{S}_j$ also hold. Thus, by induction, $\{S_j, \mathcal{S}_j\}_{j=1}^{m+3}$. This shows that $\{V_j, W_j\}_1^{m+3}$ are $C^1$ functions satisfying $\partial V_j, \partial W_j \in L^1([0,\infty))$ for $j = 1, \dots, m+3$. Integration shows that $V_j, W_j \in L^\infty([0,\infty))$ for $j = 1, \dots, m+3$. \proofendcontinue

We define $\hat{V}$Ê and $\hat{W}$ by
$$\hat{V}(t,k) = I + \frac{V_1(t)}{k} + \cdots + \frac{V_{m+3}(t)}{k^{m+3}}, \qquad
\hat{W}(t,k) = \frac{W_1(k)}{k} + \cdots + \frac{W_{m+3}(t)}{k^{m+3}},$$
so that
$$\hat{U}(t,k) = \hat{V}(t,k) + \hat{W}(t,k) e^{8ik^3t\sigma_3},$$

\medskip\noindent
{\bf Claim 6.}  
There exists a $K> 0$ such that $\hat{V}(t,k)^{-1}$ and $\hat{W}(t,k)^{-1}$ exist for all $k \in \C$ with $|k| \geq K$. 
Moreover, letting $A = -4ik^3\sigma_3 + \mathsf{V}$ and 
\begin{align*}
\begin{cases}
\hat{A}_1(t,k) = \big(\hat{V}_t(t,k) - 4ik^3 \hat{V}(t,k) \sigma_3\big)\hat{V}(t,k)^{-1}, 
	\\
\hat{A}_2(t,k) = \big(\hat{W}_t(t,k) + 4ik^3 \hat{W}(t,k) \sigma_3\big)\hat{W}(t,k)^{-1},
\end{cases}
\qquad t \geq 0, \quad |k| \geq K,
\end{align*}
the differences 
$$\Delta_l(t,k) = A(t,k) - \hat{A}_l(t,k), \qquad l = 1,2,$$ 
satisfy
\begin{align*}
|\partial_k^j\Delta_l(t,k)| \leq \frac{C + f(t)}{|k|^{m+1}}, \qquad t \geq 0, \quad |k| \geq K, \quad j = 0, 1, \dots, n, \quad l = 1,2,
\end{align*}
where $f$ is a function in $L^1([0,\infty)) \cap C([0,\infty))$. 
In particular,
\begin{align}\label{GGDeltaest}
\|\Delta_l(\cdot, k)\|_{L^1([0,t])} \leq \frac{Ct}{|k|^{m+1}}, \qquad t \geq 0, \quad |k| \geq K, \quad l = 1,2.
\end{align}

{\it Proof of Claim 6.}
The proof uses Claim 5 and is similar to that of Claim 2.
\proofendcontinue

\medskip\noindent
{\bf Claim 7.} We have
\begin{subequations}
\begin{align}\label{VUWa}
& \bigg| \frac{\partial^j}{\partial k^j}\Big[\hat{V}(t,k) e^{-4ik^3t\hat{\sigma}_3} \hat{V}^{-1}(0,k) - \hat{U}(t,k)\Big]_2 \bigg| \leq C(1+t)^j\frac{1 + |e^{8ik^3t}|}{|k|^{m+4- 2j}}, 
	\\ \label{VUWb}
& \bigg| \frac{\partial^j}{\partial k^j}\Big[\hat{W}(t,k) e^{4ik^3t\hat{\sigma}_3} \hat{W}^{-1}(0,k) - \hat{U}(t,k)e^{-8ik^3t\sigma_3}\Big]_2 \bigg| \leq C(1+t)^j\frac{1 + |e^{8ik^3t}|}{|k|^{m+4- 2j}},
\end{align}
for all $t \geq 0$, $|k| \geq K$, and $j = 0, 1, \dots, n$.
\end{subequations}

{\it Proof of Claim 7.}
We write
$$U_{formal}(t,k) = V_{formal}(t,k) + W_{formal}(t,k)e^{8ik^3t\sigma_3},$$
where
$$V_{formal}(t,k) = I + \frac{V_1(t)}{k} + \frac{V_2(t)}{k^2} + \cdots, \qquad
W_{formal}(t,k) = \frac{W_1(k)}{k} + \frac{W_2(t)}{k^2} + \cdots.$$
Then
\begin{align}\label{VVU}
V_{formal}(t,k) e^{-4ik^3t\hat{\sigma}_3} V_{formal}^{-1}(0,k) = U_{formal}(t,k)
\end{align}
to all orders in $k$. Indeed, both sides of (\ref{VVU}) are formal solutions of (\ref{tpartG}) satisfying the same initial condition at $t = 0$. Truncating (\ref{VVU}) at order $k^{-m-3}$, it follows that
$$\hat{V}(t,k) e^{-4ik^3t\hat{\sigma}_3} \hat{V}^{-1}(0,k) = \hat{U}(t,k) + O(k^{-m-4}) + O(k^{-m-4})e^{8ik^3t\sigma_3}.$$
Using Claim 5 and estimating the inverse $\hat{V}^{-1}(0,k)$ as in Claim 2, we find (\ref{VUWa}) for $j = 0$. Using the estimate $|\partial_k e^{\pm 8ik^3t}| \leq C |t k^2 e^{\pm 8ik^3t}|$, we find (\ref{VUWa}) also for $j \geq 1$. 

Similarly, we have
$$W_{formal}(t,k) e^{4ik^3t\hat{\sigma}_3} W_{formal}^{-1}(0,k) = U_{formal}(t,k)e^{-8ik^3t\sigma_3}$$
to all orders in $k$ and truncation leads to (\ref{VUWb}).
\proofendcontinue

\medskip\noindent
{\bf Claim 8.} $[U]_2$ satisfies (\ref{Gasymptoticsb}).

{\it Proof of Claim 8.}
Using that $U(t,k)$ satisfies (\ref{tpartG}), we compute
\begin{align*}
  (\hat{V}^{-1}U)_t & = - \hat{V}^{-1}\hat{V}_t \hat{V}^{-1} U + \hat{V}^{-1} U_t
  	\\
&  = - \hat{V}^{-1}(\hat{A}_1\hat{V} + 4ik^3\hat{V}\sigma_3) \hat{V}^{-1} U + \hat{V}^{-1} (AU + 4ik^3U\sigma_3)
	\\
&  = \hat{V}^{-1}\Delta_1 U - 4ik^3[\sigma_3, \hat{V}^{-1}U].
\end{align*}
Hence
$$\Big(e^{4ik^3t\hat{\sigma}_3}\hat{V}^{-1}U\Big)_t = e^{4ik^3t\hat{\sigma}_3}\hat{V}^{-1} \Delta_1 U.$$
Integrating and using the initial condition $U(0,k) = I$, we conclude that $U$ satisfies the Volterra integral equation
\begin{align}\label{GGUVolterra}
U(t,k) = \hat{V}(t,k) e^{-4ik^3t\hat{\sigma}_3}\hat{V}^{-1}(0,k) + \int_0^t \hat{V}(t,k) e^{4ik^3(t'-t)\hat{\sigma}_3} (\hat{V}^{-1}\Delta_1 U)(t',k) dt'.
\end{align}
Letting $\Psi = [U]_2$ and $\Psi_0(t,k) = [\hat{V}(t,k) e^{-4ik^3t\hat{\sigma}_3}\hat{V}^{-1}(0,k) ]_2$, we can write the second column of (\ref{GGUVolterra}) as
$$\Psi(t,k) = \Psi_0(t,k) + \int_0^t E(t,t',k) (\Delta_1 \Psi)(t',k) dt',$$
where
$$E(t,t',k) = \hat{V}(t,k) \begin{pmatrix} e^{8ik^3(t'-t)} & 0 \\ 0 & 1 \end{pmatrix} \hat{V}^{-1}(t',k).$$
We seek a solution $\Psi(t,k) = \sum_{l=0}^\infty \Psi_l(t,k)$ where
$$\Psi_l(t,k) = (-1)^l \int_{0 \leq t_1 \leq \cdots \leq t_l \leq t < \infty} \prod_{i = 1}^l E(t_{i+1}, t_i, k) \Delta_1(t_i, k) \Psi_0(t_1, k) dt_1 \cdots dt_l.$$ 
The estimates
\begin{align}\label{GGDeltaEboundb}
|\partial_k^jE(t, t', k)| < C|k|^{2j}(1+|t' - t|)^j, \quad 0 \leq t' \leq t < \infty, \;\; k \in \bar{D}_-^K, \;\; j = 0, 1, \dots, n,
\end{align}
and
$$|\Psi_0(t,k)| \leq C, \qquad t \geq 0, \quad k \in \bar{D}_-^K,$$
together with  (\ref{GGDeltaest}) yield
\begin{align*}\nonumber
|\Psi_l(t,k)| \leq \; & C \int_{0 \leq t_1 \leq \cdots \leq t_l \leq t < \infty} \prod_{i = 1}^l  |\Delta_1(t_i, k)| |\Psi_0(t_1, k)| dt_1 \cdots dt_l 
	\\ \nonumber
\leq & \; \frac{C}{l!}\|\Psi_0(\cdot, k)\|_{L^\infty([0,t])} \|\Delta_1(\cdot, k)\|_{L^1([0,t])}^l
	\\ 
\leq &\; \frac{C}{l!} \bigg(\frac{Ct}{|k|^{m+1}}\bigg)^l, \qquad t \geq 0, \quad k \in \bar{D}_-^K.
\end{align*}
Hence
\begin{align}\label{GGPsiPsi0}
|\Psi(t,k) - \Psi_0(t,k)| \leq \sum_{l=1}^\infty | \Psi_l(t,k)| \leq \frac{Ct e^{\frac{Ct}{|k|^{m+1}}}}{|k|^{m+1}}, \qquad t \geq 0, \quad k \in \bar{D}_-^K.
\end{align}
Equations (\ref{VUWa}) and (\ref{GGPsiPsi0}) prove the second column of (\ref{Gasymptoticsb}) for $j = 0$.

Differentiating the integral equation (\ref{GGUVolterra}) with respect to $k$, we find that $\Lambda := \partial_k\Psi$ satisfies
\begin{align}\label{GGintegraleq3}
\Lambda(t,k) = \Lambda_0(t,k)
+ \int_0^t E(t,t',k) \Delta_1(t', k) \Lambda(t', k) dt'
\end{align}
for each $k$ in the interior of $\bar{D}_-^K$, where
\begin{align}\label{GGLambda0}
\Lambda_0(t,k) = \;& [\partial_k\Psi_0(t,k)]_2 
+ \int_0^t \frac{\partial}{\partial k}\big[E(t,t',k) \Delta_1(t', k)\big] \Psi(t', k) dt'.
\end{align}
We seek a solution of (\ref{GGintegraleq3}) of the form $\Lambda = \sum_{l=0}^\infty \Lambda_l$. Proceeding as above, we find
$$|\Lambda_l(t,k)| \leq \frac{C}{l!}\|\Lambda_0(\cdot, k)\|_{L^\infty([0, t])}\bigg(\frac{Ct}{|k|^{m+1}}\bigg)^l, \qquad t \geq 0, \quad k \in \bar{D}_-^K.$$
Using (\ref{Gestb}) and (\ref{GGDeltaEboundb}) in (\ref{GGLambda0}), we obtain
\begin{align}\label{GGFF1}
|\Lambda_0(t,k) - [\partial_k\Psi_0(t,k)]_2| \leq \frac{C(1+t)^2}{|k|^{m-1}}, \qquad t \geq 0, \quad k \in \bar{D}_-^K.
\end{align}
Thus $\sum_{l=0}^\infty \Lambda_l$ converges uniformly on compact subsets of $[0,\infty) \times \bar{D}_-^K$ to a continuous solution $\Lambda$ of (\ref{GGintegraleq3}), which satisfies
\begin{align}\label{GGFF2}
|\Lambda(t,k) - \Lambda_0(t,k)| 
\leq \frac{C(1+t)^2te^{\frac{Ct}{|k|^{m+1}}}}{|k|^{m+1}}, \qquad t \geq 0,  \quad k \in \bar{D}_-^K.
\end{align}
Equations (\ref{VUWa}), (\ref{GGFF1}), and (\ref{GGFF2}) show that $[U]_2 = \Psi$ satisfies (\ref{Gasymptoticsb}) for $j = 1$. 
Extending the above argument, we find that (\ref{Gasymptoticsb}) holds also for $j = 2, \dots, n$.
\proofendcontinue

\medskip\noindent
{\bf Claim 9.}  $[U]_2$ satisfies (\ref{Gasymptoticsc}). 

{\it Proof of Claim 9.}
Let $w(t,k) = U(t,k) e^{-8ik^3t\sigma_3}$. Then $w$ satisfies $w_t = Aw - 4ik^3w\sigma_3$. 
Thus
\begin{align*}
  (\hat{W}^{-1}w)_t & = - \hat{W}^{-1}\hat{W}_t \hat{W}^{-1} w + \hat{W}^{-1} w_t
  	\\
&  = - \hat{W}^{-1}(\hat{A}_2\hat{W} - 4ik^3\hat{W}\sigma_3) \hat{W}^{-1} w + \hat{W}^{-1} (Aw - 4ik^3w\sigma_3)
	\\
&  = \hat{W}^{-1}\Delta_2 w + 4ik^3[\sigma_3, \hat{W}^{-1}w].
\end{align*}
Hence
$$\Big(e^{-4ik^3t\hat{\sigma}_3}\hat{W}^{-1}w\Big)_t = e^{-4ik^3t\hat{\sigma}_3}\hat{W}^{-1} \Delta_2 w.$$
Integrating and using the initial condition $w(0,k) = I$, we conclude that $w$ satisfies the Volterra integral equation
\begin{align}\label{GGWhatVolterra}
w(t,k) = \hat{W}(t,k) e^{4ik^3t\hat{\sigma}_3}\hat{W}^{-1}(0,k) + \int_0^t \hat{W}(t,k) e^{4ik^3(t-t')\hat{\sigma}_3} (\hat{W}^{-1}\Delta_2 w)(t',k) dt'.
\end{align}
Letting $\Psi = [w]_2$ and $\Psi_0(t,k) = [\hat{W}(t,k) e^{4ik^3t\hat{\sigma}_3}\hat{W}^{-1}(0,k) ]_2$, we can write the second column of (\ref{GGWhatVolterra}) as
$$\Psi(t,k) = \Psi_0(t,k) + \int_0^t E(t,t',k) (\Delta_2 \Psi)(t',k) dt',$$
where
$$E(t,t',k) = \hat{W}(t,k) \begin{pmatrix} e^{8ik^3(t-t')} & 0 \\ 0 & 1 \end{pmatrix} \hat{W}^{-1}(t',k).$$
As in the proof of Claim 8, the estimates
$$|\partial_k^jE(t, t', k)| < C|k|^{2j}(1+|t' - t|)^j, \quad 0 \leq t' \leq t < \infty, \;\; k \in \bar{D}_+^K, \;\; j = 0, 1, \dots, n,$$
and
$$|\Psi_0(t,k)| \leq C, \qquad t \geq 0, \quad k \in \bar{D}_+^K,$$
together with  (\ref{GGDeltaest}) yield
\begin{align}\label{FFPsiPsi0}
|\Psi(t,k) - \Psi_0(t,k)| \leq \frac{Ct e^{\frac{Ct}{|k|^{m+1}}}}{|k|^{m+1}}, \qquad t \geq 0, \quad k \in \bar{D}_+^K.
\end{align}
Equations (\ref{VUWb}) and (\ref{FFPsiPsi0}) prove the second column of (\ref{Gasymptoticsc}) for $j = 0$. Proceeding as in the proof of Claim 8, (\ref{Gasymptoticsc}) follows also for $j = 1, \dots, n$. 
\proofend

\subsection{The spectral functions $\{A(k), B(k)\}$}
We let $S(k) = T(0,k)$ and define the spectral functions $A(k)$ and $B(k)$ for $k \in \bar{D}_+$ by 
\begin{align}\label{ABdef}
S(k) = \begin{pmatrix} 
\overline{A(\bar{k})} 	&	B(k)	\\
\lambda \overline{B(\bar{k})}	&	A(k)
\end{pmatrix}, \qquad k \in (\bar{D}_-, \bar{D}_+).
\end{align}

\begin{theorem}\label{ABth}
Suppose $\{g_j(t)\}_0^2$ satisfy (\ref{gjassump}) for some integers $m \geq 1$ and $n \geq 1$.
Then the spectral functions  $A(k)$ and $B(k)$ have the following properties:
\begin{enumerate}[$(a)$]
\item $A(k)$ and $B(k)$ are continuous for $k \in \bar{D}_+$ and analytic for $k \in D_+$.

\item There exist complex constants $\{A_j, B_j\}_1^m$ such that
\begin{subequations}\label{ABexpansions}
\begin{align}\nonumber
& A(k) = 1 + \frac{A_1}{k} + \cdots + \frac{A_m}{k^m} + O\bigg(\frac{1}{k^{m+1}}\bigg), 
	\\
& B(k) = \frac{B_1}{k} + \cdots + \frac{B_m}{k^m} + O\bigg(\frac{1}{k^{m+1}}\bigg),	
\end{align}
uniformly as $k \to \infty$ with $k \in \bar{D}_+$.

\item For $j = 1, \dots n$, the derivatives $A^{(j)}(k)$ and $B^{(j)}(k)$ have continuous extensions to $k \in \bar{D}_+$ and 
\begin{align}\nonumber
& A^{(j)}(k) = \frac{d^j}{dk^k}\bigg(1 + \frac{A_1}{k} + \cdots + \frac{A_m}{k^m}\bigg) 
+ O\bigg(\frac{1}{k^{m+1-2j}}\bigg),
	\\
& B^{(j)}(k) = \frac{d^j}{dk^j}\bigg(\frac{B_1}{k} + \cdots + \frac{B_m}{k^m}\bigg) + O\bigg(\frac{1}{k^{m+1-2j}}\bigg),
\end{align}
\end{subequations}
uniformly as $k \to \infty$ with $k \in \bar{D}_+$.

\item $A$ and $B$ obey the symmetries
\begin{align}\label{ABsymm}
\begin{cases}
A(k) = \overline{A(-\bar{k})}, \\ 
B(k) = \overline{B(-\bar{k})},
\end{cases} \qquad k \in \bar{D}_+.
\end{align}

\item $A$ and $B$ satisfy the relation
\begin{align}\label{ABdetone}
A(k)\overline{A(\bar{k})} - \lambda B(k) \overline{B(\bar{k})} = 1, \qquad k \in \Gamma.
\end{align}

\end{enumerate}
\end{theorem}
\proofbegin 
Letting $A_j = (T_j(0))_{22}$ and $B_j = (T_j(0))_{12}$, properties $(a)$-$(c)$ follow immediately from Theorems \ref{tth1} and \ref{tth2}. 
Property $(d)$ is a consequence of the symmetry $T(x,k) = \overline{T(x,-\bar{k})}$.
Property $(e)$ follows since $\det T = 1$. 
\proofend

\subsection{More spectral functions}
Suppose $u_0, g_0, g_1, g_2$ satisfy (\ref{uassump}) and (\ref{gjassump}) for some integers $m \geq 1$ and $n \geq 1$. 
Let $\{a(k), b(k), A(k), B(k)\}$ be given by (\ref{abdef}) and (\ref{ABdef}). 

We define the spectral functions $c(k)$ and $d(k)$ by
\begin{align*}
& c(k) = A(k)b(k) - B(k)a(k), \qquad k \in \bar{D}_1 \cup \R,
	\\
& d(k) = a(k)\overline{A(\bar{k})} -  \lambda b(k) \overline{B(\bar{k})}, \qquad k \in \bar{D}_2.
\end{align*}
Then
$$S(k)^{-1}s(k) = \begin{pmatrix} \overline{d(\bar{k})} & c(k) \\ \lambda \overline{c(\bar{k})} & d(k) \end{pmatrix}, \qquad k \in (\partial D_3, \partial D_2).$$

We also define spectral functions $h(k)$ and $r(k)$ by
\begin{subequations}\label{hrdef}
\begin{align}\label{hdef}
& h(k) = -\frac{ \overline{B(\bar{k})}}{a(k) d(k)}, \qquad k \in \bar{D}_2,
	\\ \label{rdef}
& r(k) = \frac{\overline{c(\bar{k})}}{d(k)}
= \frac{\overline{b(\bar{k})}}{a(k)} + h(k), \qquad k \in \R.
\end{align}
 \end{subequations}

\section{The mKdV equation in the quarter plane}\label{mkdvsec}\nequation
In this section, we apply the results from the preceding sections to express the solution of the mKdV equation in the quarter plane in terms of the solution of a RH problem. 

Before stating the main result, we need to recall some definitions related to $L^p$-RH problems. We use the notation of \cite{LCarleson}. Further details can be found in \cite{LCarleson} (see also the appendix). Let $\mathcal{J}$ denote the collection of all subsets $\gamma$ of the Riemann sphere $\hat{\C} = \C \cup \{\infty\}$ such that $\gamma$ is homeomorphic to the unit circle and
\begin{align}\label{carlesondef}
 \sup_{z \in \gamma \cap \C} \sup_{r > 0} \frac{|\gamma \cap D(z, r)|}{r} < \infty,
\end{align}
where $D(z, r)$ denotes the disk of radius $r$ centered at $z$. Curves satisfying (\ref{carlesondef}) are called Carleson curves. Let $1 \leq p < \infty$. If $D$ is the bounded component of $\hat{\C} \setminus \gamma$ where $\gamma \in \mathcal{J}$ and $\infty \notin \gamma$, then a function $f$ analytic in $D$ belongs to the Smirnoff class $E^p(D)$ if there exists a sequence of rectifiable Jordan curves $\{C_n\}_1^\infty$ in $D$, tending to the boundary in the sense that $C_n$ eventually surrounds each compact subdomain of $D$, such that $\sup_{n \geq 1} \int_{C_n} |f(z)|^p |dz| < \infty$. 
If $D \subset \hat{\C}$ is bounded by an arbitrary curve in $\mathcal{J}$, $E^p(D)$ is defined as the set of functions $f$ analytic in $D$ for which $f \circ \varphi^{-1} \in E^p(\varphi(D))$, where $\varphi(z) = \frac{1}{z - z_0}$ and $z_0$ is any point in $\C \setminus \bar{D}$. The subspace of $E^p(D)$ consisting of all functions $f \in E^p(D)$ such that $z f(z) \in E^p(D)$ is denoted by $\dot{E}^p(D)$.
We let $E^\infty(D)$ denote the space of bounded analytic functions in $D$. 

Let $n \geq 1$ be an integer and let $\Sigma$ be a Carleson jump contour. Given an $n \times n$-matrix valued function $v: \Sigma \to GL(n, \C)$, a {\it solution of the $L^p$-RH problem determined by $(\Sigma, v)$} is an $n \times n$-matrix valued function $m \in I + \dot{E}^p(\hat{\C} \setminus \Sigma)$ such that the nontangential boundary values $m_\pm$ satisfy $m_+ = m_- v$ a.e. on $\Sigma$. 

Given a Carleson jump contour $\Sigma$ and $a,b \in \R$ with $a < b$, we call $W_{a,b} = \{a \leq \arg k \leq b\}$ a nontangential sector at $\infty$ if there exists a $\delta > 0$ such that $W_{a - \delta, b + \delta}$ does not intersect $\Sigma \cap \{|z| > R\}$ whenever $R>0$ is large enough. If $f(k)$ is a function of $k \in \C \setminus \Sigma$, we say that $f$ has nontangential limit $L$ at $\infty$, written 
$$\ntlim_{k\to\infty} f(k) = L,$$ 
if $\lim_{\substack{k \to \infty \\ k \in W_{a,b}}}f(k) = L$ for every nontangential sector $W_{a,b}$ at $\infty$.

The following theorem expresses the solution of (\ref{mkdv}) in the quarter plane $\{x\geq 0, t \geq 0\}$ in terms of the solution of an $L^2$-RH problem. 

\begin{theorem}\label{existenceth}
Suppose $u_0, g_0, g_1, g_2$ satisfy (\ref{uassump}) and (\ref{gjassump}) with $n = 1$ and $m = 4$, i.e., suppose
\begin{align*}
\begin{cases}
 u_0 \in C^5([0,\infty)), \qquad g_0 \in C^3([0,\infty)), \qquad g_1, g_2 \in C^2([0,\infty)),
	\\
(1 + x) u_0^{(i)}(x) \in L^1([0,\infty)), \qquad i = 0, 1, \dots, 5,
	\\
 (1 + t) g_0^{(i)}(t) \in L^1([0,\infty)), \qquad i = 0, 1, 2, 3,
	\\
 (1 + t) g_1^{(i)}(t), (1 + t) g_2^{(i)}(t) \in L^1([0,\infty)), \qquad i = 0, 1, 2.
\end{cases}
\end{align*}
Define the spectral functions $h(k)$ and $r(k)$ by (\ref{hrdef}). Define the jump matrix $J(x,t,k)$ by 
\begin{align}\label{Jdef}
&J(x,t,k) = \begin{cases} 
 \begin{pmatrix} 1 & 0 \\ \lambda h(k) e^{-2ikx + 8ik^3t} & 1 \end{pmatrix}, & k \in \partial D_1,
	\\
 \begin{pmatrix} 1 & - \overline{r(\bar{k})} e^{2ikx - 8ik^3t} \\
\lambda r(k)e^{-2ikx + 8ik^3t}& 1 -  \lambda|r(k)|^2\end{pmatrix}, & k \in \R,
	\\ 
 \begin{pmatrix} 1 & - \overline{h(\bar{k})} e^{2ikx - 8ik^3t} \\ 0 & 1 \end{pmatrix}, & k \in \partial D_4.
\end{cases}
\end{align}
If $\lambda = -1$, assume that $a(k)$ is nonzero for $k \in \bar{\C}_-$ and that $d(k)$ is nonzero for $k \in \bar{D}_2$. If $\lambda = 1$, suppose the homogeneous RH problem determined by $(\Gamma, J(x,t,\cdot))$ (see equation (\ref{homogenousRH})) has only the trivial solution for each $(x,t) \in [0, \infty) \times [0,\infty)$.

Suppose the spectral functions satisfy
\begin{align}\label{GR}
& A(k)b(k) - B(k)a(k) = 0, \qquad k \in \bar{D}_1.
\end{align}

Then the $L^2$-RH problem
\begin{align}\label{RHM}
\begin{cases}
M(x, t, \cdot) \in I + \dot{E}^2(\C \setminus \Gamma),\\
M_+(x,t,k) = M_-(x, t, k) J(x, t, k) \quad \text{for a.e.} \ k \in \Gamma,
\end{cases}
\end{align}
has a unique solution for each $(x,t) \in [0,\infty) \times [0, \infty)$. Moreover, the nontangential limit
\begin{align}\label{ulim}
u(x,t) = -2i\ntlim_{k\to\infty} (kM(x,t,k))_{12}
\end{align}
exists for each $(x,t) \in [0,\infty) \times [0, \infty)$ and the function $u(x,t)$ defined by (\ref{ulim}) has the following properties:
\begin{enumerate}[$(a)$]
\item $u:[0,\infty) \times [0, \infty) \to \R$ is $C^3$ in $x$ and $C^1$ in $t$.

\item $u(x,t)$ satisfies (\ref{mkdv}) for $x \geq 0$ and $t \geq 0$.

\item $u(0,t) = g_0(t)$, $u_x(0,t) = g_1(t)$, and $u_{xx}(0,t) = g_2(t)$ for $t \geq 0$.

\item $u(x,0) = u_0(x)$ for $x \geq 0$.
\end{enumerate}
\end{theorem}
\proofbegin
The proof proceeds through a series of claims. 

\medskip\noindent
{\bf Claim 1.} $a(k)$ is nonzero for $k \in \bar{\C}_-$ and $d(k)$ is nonzero for $k \in \bar{D}_2$.

{\it Proof of Claim 1.}
For $\lambda = -1$ this holds by assumption. If $\lambda = 1$, Theorem \ref{abth} shows that $a(k)$ is nonzero for $k \in \bar{\C}_-$, while the arguments of \cite{Lsolitonfree} show that $d(k)$ is nonzero in $\bar{D}_2$.
\proofendcontinue

In view of Claim 1, Theorems \ref{abth} and \ref{ABth} imply that $r(k)$ and $h(k)$ have the following properties:
\begin{itemize}
\item $r \in C^1(\R)$.

\item $h$ is analytic in $D_2$ and $h, h'$ have continuous extensions to $\bar{D}_2$.

\item There exist complex constants $\{r_j, h_j\}_1^4$ such that
\begin{align}
r^{(j)}(k) = \frac{d^j}{dk^j}\bigg(\frac{r_1}{k} + \cdots + \frac{r_4}{k^4}\bigg) + O(k^{-5}), \qquad |k| \to \infty, \quad k \in \R, \quad j = 0,1,
	\\
h^{(j)}(k) = \frac{d^j}{dk^j}\bigg(\frac{h_1}{k} + \cdots + \frac{h_4}{k^4}\bigg) + O(k^{-5}), \qquad k \to \infty, \quad k \in \bar{D}_2, \quad j = 0,1.
\end{align}

\item $r(k) = \overline{r(-\bar{k})}$ for $k \in \R$ and $h(k) = \overline{h(-\bar{k})}$ for $k \in \bar{D}_2$.
\end{itemize}
Relation (\ref{GR}) implies that $r_j = 0$ for $j = 1, \dots, 4$, so that in fact
\begin{align}\label{rkm4}
r^{(j)}(k) = O(k^{-5}), \qquad |k| \to \infty, \quad k \in \R, \quad j = 0,1.
\end{align}
Indeed, the expansions in (\ref{abexpansions}) of $\{a(k), b(k)\}$ are valid as $k \to \infty$ in $\bar{D}_1 \cup \bar{D}_2$ and the expansions in (\ref{ABexpansions}) of $\{A(k), B(k)\}$ are valid as $k \to \infty$ in $\bar{D}_1 \cup \bar{D}_3$.
Hence $c(k)$  has an expansion
$$c(k) = \frac{c_1}{k} + \cdots + \frac{c_4}{k^4} + O(k^{-5}), \qquad k \to \infty, \qquad k \in \bar{D}_1 \cup \R,$$
where the coefficients $\{c_j\}_1^4$ are the {\it same} as $k \to \infty$ in $\bar{D}_1$ and in $\R$. Since $c(k)$ vanishes identically in $\bar{D}_1$ by (\ref{GR}), we have $c_j = 0$ for $j = 1, \dots, 4$. Since $r(k) = \frac{\overline{c(\bar{k})}}{d(k)}$, this proves (\ref{rkm4}).

\medskip\noindent
{\bf Claim 2.} (Vanishing Lemma) Let $x \geq 0$ and $t \geq 0$. Suppose $N(x,t,k)$ is a solution of the homogeneous $L^2$-RH problem determined by $(\Gamma, J(x,t,\cdot))$, i.e. 
\begin{align}\label{homogenousRH}
\begin{cases}
N(x,t,\cdot) \in \dot{E}^2(\C \setminus \Gamma),\\
N_+(x,t,k) = N_-(x,t,k) J(x, t, k) \quad \text{for a.e.} \ k \in \Gamma.
\end{cases}
\end{align}
Then $N$ vanishes identically.

{\it Proof of Claim 2.}
For $\lambda = 1$, this holds by assumption. Thus suppose $\lambda = -1$. 
We write $N(k) := N(x,t,k)$ and let $G(k) = N(k) \overline{N(\bar{k})}^T$. Applying Lemma \ref{intfmlemma} with  $m=2$ and $n =0$, we obtain $G \in \dot{E}^1(\C \setminus \Gamma)$ and
$$\int_{\partial D_3} G_+(k) dk = 0, \qquad \int_{\partial D_4} G_-(k) dk = 0.$$
Adding these equations, we find
$$\int_\R G_+ dk - \int_{\partial D_4} (G_+ - G_-)dk = 0.$$
Since
\begin{align*}
G_+(k) & = N_+(k) \overline{N_-(\bar{k})}^T
= N_-(k) J(x,t,k) \overline{N_-(\bar{k})}^T
	\\
&= N_-(k) \overline{J(x,t,\bar{k})}^T \overline{N_-(\bar{k})}^T
 = G_-(k), \qquad k \in \Gamma,
\end{align*}
this gives
$$0 = \int_\R G_+(k) dk = \int_\R N_-(k)J(x,t,k) \overline{N_-(k)}^T dk.$$
For $k \in \R$, $J(x,t,k)$ is a Hermitian matrix with $(11)$ entry $1$ and determinant one; by Sylvester's criterion it is positive definite. 
It follows that $N_- = 0$ a.e. on $\R$. But then $N_+ = N_- J$ also vanishes a.e. on $\R$. 

Let $B_r$ be a small open ball contained in $D_2 \cup \R \cup D_3$ centered at some nonzero real number. The function $\tilde{N}(k)$ defined by
$$\tilde{N}(k) = \frac{1}{2\pi i} \int_{\partial B_r} \frac{N(s)}{s-k}ds$$
is analytic in $B_r$ and equals $N(k)$ for $k \in B_r\setminus \R$. Indeed, if $k \in B_r \cap \C_+$, then
$$\tilde{N}(k) = \frac{1}{2\pi i} \int_{\partial (B_r \cap \C_+)} \frac{N_+(s)}{s-k}ds + \frac{1}{2\pi i} \int_{\partial (B_r \cap \C_-)} \frac{N_-(s)}{s-k}ds
= N(k) + 0,$$
and a similar argument applies if $k \in B_r \cap \C_-$. We infer that $N(k)$ Êis analytic in $B_r$. Since $N = 0$ on $B_r \cap \R$, it follows by analytic continuation that $N$ vanishes identically for $k \in D_2 \cup \R \cup D_3$. But then $N_\pm = 0$ on $\partial D_1 \cup \partial D_4$ so by a similar argument we find that $N$ vanishes for all $k \in \C \setminus \Gamma$.
\proofendcontinue

The global relation (\ref{GR}) implies that $h(0) = - \frac{\overline{b(0)}}{a(0)}$ and $r(0) = 0$.\footnote{The symmetries (\ref{absymm}) and (\ref{ABsymm}) imply that the values of $a,b,A,B$ at $k = 0$ are real.}
It follows that unless $b(0) = 0$, the matrix $J$ is not continuous at $k = 0$. 
In order to obtain a jump matrix which is continuous and which approaches the identity matrix sufficiently fast as $k \to \infty$, we introduce $m(x,t,k)$ by
$$m(x,t,k) 
= \begin{cases} M(x,t,k)  \begin{pmatrix} 1 & 0 \\ -\lambda h_a(k) e^{-2ikx + 8ik^3t} & 1 \end{pmatrix}, & k \in D_1,
	\\
M(x,t,k)  \begin{pmatrix} 1 & -\overline{h_a(\bar{k})} e^{2ikx - 8ik^3t} \\ 0 & 1 \end{pmatrix}, & k \in D_4,
	\\
M(x,t,k), & \text{otherwise},		
\end{cases}$$
where $h_a(k)$ is a rational function such that $h_a$ has no poles in $\bar{D}_1$, $h_a(0) = h(0)$, $h_a(k) = \overline{h_a(-\bar{k})}$, and 
$$h_a(k) = \frac{h_1}{k} + \cdots + \frac{h_4}{k^4} + O(k^{-5}), \qquad k \to \infty, \quad k \in \C.$$
Then
\begin{align}\label{hhadecay}
h(k) - h_a(k) = O(k^{-5}), \qquad k \to \infty, \quad k \in \partial D_1.
\end{align}
It is easy to see that such a function $h_a$ exists.

Lemma \ref{deformationlemma} implies that the $L^2$-RH problem for $M$ is equivalent to the $L^2$-RH problem
\begin{align}\label{mRH}
\begin{cases}
m(x, t, \cdot) \in I + \dot{E}^2(\C \setminus \Gamma),\\
m_+(x,t,k) = m_-(x, t, k) v(x, t, k) \quad \text{for a.e.} \ k \in \Gamma,
\end{cases}
\end{align}
where
\begin{align}\label{vdef}
&v(x,t,k) = \begin{cases} 
\begin{pmatrix} 1 & 0 \\ \lambda (h(k) - h_a(k)) e^{-2ikx + 8ik^3t} & 1 \end{pmatrix}, \qquad k \in \partial D_1,
	\\
\begin{pmatrix} 1 & - \overline{r(\bar{k})} e^{2ikx - 8ik^3t} \\
\lambda r(k)e^{-2ikx + 8ik^3t} & 1 -  \lambda |r(k)|^2\end{pmatrix}, \qquad k \in \R,
	\\ 
\begin{pmatrix} 1 & - (\overline{h(\bar{k})}- \overline{h_a(\bar{k})}) e^{2ikx - 8ik^3t} \\ 0 & 1 \end{pmatrix}, \qquad k \in \partial D_4,
\end{cases}
\end{align}

Let $\mathcal{B}(L^2(\Gamma))$ denote the space of bounded linear operators on $L^2(\Gamma)$. 
Defining the nilpotent matrices $w^\pm(x,t,k)$ by
$$w^- = \begin{cases}
\begin{pmatrix} 0 & 0 \\ \lambda(h(k) - h_a(k)) e^{-2ikx + 8ik^3t} & 0 \end{pmatrix}, & k \in \partial D_1,
	\\
\begin{pmatrix} 0 & 0 \\ \lambda r(k) e^{-2ikx + 8ik^3 t} & 0 \end{pmatrix}, & k \in \R,
	\\
0, & k \in \partial D_4,
\end{cases}
$$
and
$$w^+ = \begin{cases}
0, & k \in \partial D_1,
	\\
\begin{pmatrix} 0 & -\overline{r(\bar{k})} e^{2ikx - 8ik^3t} \\ 0 & 0 \end{pmatrix}, & k \in \R,
	\\
\begin{pmatrix} 0 & -(\overline{h(\bar{k})}-  \overline{h_a(\bar{k})}) e^{2ikx - 8ik^3t} \\ 0 & 0 \end{pmatrix}, & k \in \partial D_4,
\end{cases}$$
we can write $v = (v^-)^{-1}v^+$, where $v^+ = I + w^+$ and $v^- = I - w^-$. Let $\mathcal{C}_w$ be the operator defined in (\ref{Cwdef}). For each $(x,t) \in [0,\infty) \times [0,\infty)$, we have $v^\pm \in C(\Gamma)$ and $v^{\pm}, (v^{\pm})^{-1} \in I + L^2(\Gamma) \cap L^\infty(\Gamma)$. Therefore, Claim 2 and Lemma \ref{Fredholmzerolemma} imply that $I - \mathcal{C}_w \in \mathcal{B}(L^2(\Gamma))$ is bijective and that the $L^2$-RH problem determined by $(\Gamma, v)$ has a unique solution $m(x,t,k)$ for each  $(x,t) \in [0,\infty) \times [0,\infty)$. By Lemma \ref{mulemma}, this solution is given by
$$m = I + \mathcal{C}(\mu (w^+ + w^-)) \in I + \dot{E}^2(\C \setminus \Gamma),$$
where
$$\mu = I + (I - \mathcal{C}_w)^{-1}\mathcal{C}_w I \in I + L^2(\Gamma).$$
By the open mapping theorem, $(I - \mathcal{C}_w)^{-1} \in \mathcal{B}(L^2(\Gamma))$ for each  $(x,t)$. 

\medskip\noindent
{\bf Claim 3.} The map
\begin{align}\label{xttomuI}
(x,t) \mapsto \mu(x,t, \cdot) - I: [0,\infty) \times [0, \infty) \to L^2(\Gamma)
\end{align}
is $C^3$ in $x$  and $C^1$ in $t$.

{\it Proof of Claim 3.}
In view of (\ref{rkm4}) and (\ref{hhadecay}), the maps
\begin{subequations}\label{xtww}
\begin{align}\label{xtwwa}
& (x,t) \mapsto (w^+(x,t,\cdot), w^-(x,t,\cdot)):[0, \infty) \times [0, \infty) \to L^2(\Gamma) \times L^2(\Gamma),
	\\
& (x,t) \mapsto (w^+(x,t,\cdot), w^-(x,t,\cdot)):[0, \infty) \times [0, \infty) \to L^\infty(\Gamma) \times L^\infty(\Gamma),
\end{align}
\end{subequations}
are $C^3$ in $x$ and $C^1$ in $t$.
On the other hand, the map
\begin{align}
(w^+, w^-) \mapsto I - \mathcal{C}_w: L^\infty(\Gamma) \times L^\infty(\Gamma) \to \mathcal{B}(L^2(\Gamma))
\end{align}
is smooth by the estimate
\begin{align*}
\|\mathcal{C}_w\|_{\mathcal{B}(L^2(\Gamma))} \leq C \max\big\{\|w^+\|_{L^\infty(\Gamma)}, \|w^-\|_{L^\infty(\Gamma)} \big\}.
\end{align*}
Moreover, the bilinear map
\begin{align}\label{wwCwI}
(w^+, w^-) \mapsto \mathcal{C}_wI:  L^2(\Gamma) \times L^2(\Gamma) \to L^2(\Gamma)
\end{align}
is smooth by the estimate
$$\|\mathcal{C}_wI\|_{L^2(\Gamma)} \leq C \max\big\{\|w^+\|_{L^2(\Gamma)}, \|w^-\|_{L^2(\Gamma)} \big\}.$$
Since (\ref{xttomuI}) can be viewed as a composition of maps of the form (\ref{xtww})-(\ref{wwCwI}) together with the smooth inversion map $I - \mathcal{C}_w \mapsto (I - \mathcal{C}_w)^{-1}$, it follows that (\ref{xttomuI}) is $C^3$ in $x$  and $C^1$ in $t$.
\proofendcontinue

\medskip\noindent
{\bf Claim 4.} The nontangential limit in (\ref{ulim}) exists for each $(x,t) \in [0,\infty) \times [0, \infty)$. Moreover, defining $u(x,t)$ by (\ref{ulim}), it holds that
\begin{align}\label{Mlax}
\begin{cases}
m_x - ik[\sigma_3, m] = \mathsf{U}m,
	\\
m_t + 4ik^3[\sigma_3, m] = \mathsf{V}m,
\end{cases} \qquad (x,t) \in [0,\infty) \times [0, \infty),
\end{align}
where $\mathsf{U}(x,t)$  and $\mathsf{V}(x,t,k)$ are defined by (\ref{mathsfUVdef}).

{\it Proof of Claim 4.}
For each $k \in \C \setminus \Gamma$, the linear map
$$f \mapsto (\mathcal{C}f)(k) = \frac{1}{2\pi i} \int_\Gamma \frac{f(s)ds}{s - k}: L^2(\Gamma) \to \C$$
is bounded. Also, by (\ref{xttomuI}) and (\ref{xtww}),
$$(x,t) \mapsto \mu(w^+ + w^-) = (\mu- I)(w^+ + w^-) + (w^+ + w^-): [0,\infty) \times [0,\infty) \to L^2(\Gamma)$$
is $C^3$ in $x$ and $C^1$ in $t$.
Hence, for each $k \in \C \setminus \Gamma$, 
\begin{align*}
(x,t) \mapsto m(x,t,k) = &\; I + \frac{1}{2\pi i}\int_\Gamma \frac{\mu(x,t,s)(w^+ + w^-)(x,t,s)}{s-k} ds
\end{align*}
is $C^3$ in $x$ and $C^1$ in $t$. Using that $\{\partial_x^j(\mu(w^+ + w^-))\}_0^3$ and $\partial_t(\mu(w^+ + w^-))$ exist in $L^2(\Gamma)$, we find
\begin{subequations}\label{MxMt}
\begin{align}\label{MxMta}
& (\partial_x^jm)(x,t,\cdot) = \mathcal{C}(\partial_x^j(\mu(w^+ + w^-))) \in \dot{E}^2(\C \setminus \Gamma), \qquad j = 1,2,3,
	\\
& (\partial_tm)(x,t,\cdot) = \mathcal{C}(\partial_t(\mu(w^+ + w^-))) \in \dot{E}^2(\C \setminus \Gamma),
\end{align}
\end{subequations}
for each  $(x,t)$.

Let $\psi = m e^{-i(-k x + 4k^3t)\sigma_3}$. Then $\psi_x \psi^{-1} = (m_x + ikm\sigma_3)m^{-1}$.
Since $m = I + \dot{E}^2(\C \setminus \Gamma)$ and $m_x \in \dot{E}^2(\C \setminus \Gamma)$, Lemmas \ref{EpCnlemma} and \ref{intfmlemma} show that there exist functions $f_0(x,t)$ and $f_1(x,t)$ such that the function $f$ defined by
\begin{align*}
f(x,t,k)  =\frac{1}{(k + i)^2}(m_x + ikm\sigma_3)m^{-1} 
- \frac{f_0(x,t)}{(k+i)^2} - \frac{f_1(x,t)}{k+i}
\end{align*}
lies in $\dot{E}^1(\C \setminus \Gamma)$. The jump condition (\ref{mRH}) for $m$ implies that $f_+ = f_-$ a.e. on $\Gamma$ and hence that $f$ vanishes identically:
$$f(x,t,\cdot) = \mathcal{C}(f_+ - f_-) = 0.$$
We conclude that there exist functions $F_0(x,t)$ and $F_1(x,t)$ such that 
\begin{align}\label{MxikM}
m_x + ik m\sigma_3 = (F_0(x,t) + k F_1(x,t))m, \qquad k \in \C \setminus \Gamma.
\end{align}

Let $W$ be a nontangential sector at $\infty$ with respect to $\Gamma$, that is, $W = \{k \in \C \setminus\{0\}\, | \, \alpha \leq \arg k \leq \beta\}$ where $\alpha, \beta$ are such that $W \subset \C \setminus \Gamma$.
We write
\begin{align*}
m(x,t,k) 
& = I + \frac{1}{2\pi i}\int_\Gamma \frac{(\mu(w^+ + w^-))(x,t,s)}{s-k} ds
	\\
& = I - \frac{1}{2\pi i}\int_\Gamma (\mu(w^+ + w^-))(x,t,s) \bigg(\frac{1}{k} + \frac{s}{k^2} + \frac{s^2}{k^3} + \frac{s^3}{k^4} + \frac{s^4}{k^4(k - s)}\bigg) ds.
\end{align*}
The $L^2$ norm of $\mu(x,t, \cdot) - I$ as well as the $L^1$ and $L^2$ norms of $s^3(w^+ + w^-)(x,t,s)$ and $\frac{s^4}{k - s} (w^+ + w^-)(x,t,s)$ are bounded (the latter uniformly with respect to all $k \in W$ with $|k| > 1$). Hence, 
\begin{align}\label{Mexpansiona}
m(x,t,k) = I + \sum_{j=1}^3 \frac{m_j(x,t)}{k^j} + O(k^{-4}), \qquad k \to \infty, \quad k \in W,
\end{align}
where the error term is uniform with respect to $\arg k \in [\alpha, \beta]$ and
$$m_j(x,t) = - \frac{1}{2\pi i}\int_\Gamma (\mu(w^+ + w^-))(x,t,s) s^{j-1} ds, \qquad j = 1,2,3.$$
We infer that the nontangential limit in (\ref{ulim}) exists and that
\begin{align}\label{ulim2}
u(x,t) = - 2i \ntlim_{k\to \infty} (km(x,t,k))_{12} = \frac{1}{\pi}\int_\Gamma (\mu(w^+ + w^-))_{12}(x,t,s) ds.
\end{align}

Similarly, since the $L^2$-norms of $\mu_x$, $\mu_t$ as well as the $L^1$ and $L^2$ norms of $(s^2 + \frac{s^3}{k-s})(w^+ + w^-)_x$ and $(1 + \frac{s}{k-s})(w^+ + w^-)_t$ are bounded,
\begin{align}\label{Mexpansionb}
\begin{cases}
 m_x = \sum_{j=1}^2 \frac{\partial_xm_{j}(x,t)}{k^j} + O(k^{-3}), 
	\\
m_t = O(k^{-1}), 
\end{cases}	
\qquad k \to \infty, \quad k \in W.
\end{align}
Substituting (\ref{Mexpansiona}) and (\ref{Mexpansionb}) into (\ref{MxikM}) the terms of $O(k)$ and $O(1)$ yield
$$F_1 = i\sigma_3, \qquad F_0 = -i[\sigma_3, m_1].$$
The symmetry $J(x,t,k) = \sigma_j J(x,t,-k)^{-1}\sigma_j$, where $j=1$ if $\lambda = 1$ and $j = 2$ if $\lambda = -1$, implies that $\sigma_j m(x,t,-k) \sigma_j$ satisfies the same $L^2$-RH problem as $m(x,t,k)$; so by uniqueness $m(x,t,k) = \sigma_j m(x,t,-k) \sigma_j$.
Hence
$$m_j(x,t) =  \begin{cases} (-1)^j \sigma_1 m_j(x,t) \sigma_1, \quad & \lambda = 1, 
	\\
(-1)^j \sigma_2 m_j(x,t,) \sigma_2, & \lambda = -1,
\end{cases} \quad j = 1,2,3.$$
We give the remainder of the proof of Claim 4  in the case of $\lambda = 1$; the case of $\lambda = -1$ is similar.

Assume $\lambda = 1$. Writing $m_j = a_j \sigma_1 + b_j \sigma_2 + c_j \sigma_3 + d_jI$, where $a_j,b_j,c_j,d_j$ are scalar-valued functions of $(x,t)$, we find $a_j = d_j = 0$ for $j$ odd and $b_j = c_j = 0$ for $j$ even. We infer that
\begin{align} \label{Mjs}
\begin{cases}
m_1(x,t) = b_1(x,t) \sigma_2 + c_1(x,t) \sigma_3,
  	\\
  m_2(x,t) = a_2(x,t) \sigma_1 + d_2(x,t) I,
  	\\
  m_3(x,t) = b_3(x,t) \sigma_2 + c_3(x,t) \sigma_3.
\end{cases}
\end{align}
Hence $b_1 = -u/2$ and $F_0 = -i[\sigma_3, m_1] = -2b_1 \sigma_1$; thus
$$F_0 + kF_1 = ik\sigma_3 + u(x,t) \sigma_1.$$
Recalling (\ref{MxikM}), this proves that 
\begin{align}\label{mxiksigma}
m_x - ik[\sigma_3, m] = \mathsf{U}m.
\end{align}

Substituting (\ref{Mexpansiona}) and (\ref{Mexpansionb}) into (\ref{mxiksigma}), the terms of $O(k^{-1})$ and $O(k^{-2})$ yield
\begin{align}\label{cadb}
& a_2 = \frac{u_x}{4} - \frac{iu}{2} (m_1)_{11}, \qquad
b_3 = -\frac{1}{8}\big(u^3 + 4u (m_2)_{11} + 2 i u_x (m_1)_{11} - u_{xx}\big).
\end{align}
Since $m-I$ and $m_t$ belong to Ê$\dot{E}^2(\C \setminus \Gamma)$, there exist functions $\{g_j(x,t)\}_1^4$ such that the function $g$ defined by
\begin{align*}
g(x,t,k) & = \frac{1}{(k + i)^4} \psi_t \psi^{-1} 
- \sum_{j=1}^4 \frac{g_j(x,t)}{(k+i)^j} 
	\\
& =\frac{1}{(k + i)^4}[m_t -4ik^3 m\sigma_3]m^{-1} 
- \sum_{j=1}^4 \frac{g_j(x,t)}{(k+i)^j}
\end{align*}
lies in $\dot{E}^1(\C \setminus \Gamma)$. The jump condition (\ref{mRH}) for $m$ implies that $g_+ = g_-$ a.e. on $\Gamma$ and hence that $g$ vanishes identically.
We conclude that there exist functions $\{G_j(x,t)\}_0^3$ such that
\begin{align}\label{MxikM2}
m_t - 4ik^3 m\sigma_3 = \bigg(\sum_{j=0}^3 k^j G_j(x,t) \bigg) m, \qquad k \in \C \setminus \Gamma.
\end{align}
Substituting (\ref{Mexpansiona}) and (\ref{Mexpansionb}) into (\ref{MxikM2}) the terms of $O(k^j)$, $j = 0,\dots, 3$, yield
\begin{align*}
& G_3 = -4i\sigma_3, \qquad 
G_2 = -4im_1 \sigma_3 - G_3 m_1, \qquad
G_1 = -4im_2 \sigma_3 - G_3 m_2 - G_2 m_1, 
	\\
& G_0 = -4im_3 \sigma_3 - G_3 m_3 - G_2 m_2- G_1 m_1.
\end{align*}
Equations (\ref{Mjs}) and (\ref{cadb}) now show that $\sum_{j=0}^3 k^j G_j = - 4ik^3 \sigma_3 + \mathsf{V}$.
\proofendcontinue

\medskip\noindent
{\bf Claim 5.} $u(x,t)$ is $C^3$ in $x \geq 0$ and $C^1$  in $t \geq 0$.

{\it Proof of Claim 5.}
By (\ref{ulim2}), we have
$$u(x,t) = \frac{1}{\pi}\int_\Gamma ((\mu - I)w^+)_{12}(x,t,s) ds
 + \frac{1}{\pi}\int_\Gamma (w^+)_{12}(x,t,s) ds.$$
The claim follows from the differentiability properties of the maps in (\ref{xttomuI}) and (\ref{xtwwa}) as well as the fact that the map
$$(x,t) \mapsto w^+(x,t,\cdot):[0, \infty) \times [0, \infty) \to L^1(\Gamma)$$
is $C^3$ in $x$ and $C^1$ in $t$.
\proofendcontinue

\medskip\noindent
{\bf Claim 6.} $u(x,t)$ satisfies the mKdV equation (\ref{mkdv}) for $x \geq 0$ and $t \geq 0$.

{\it Proof of Claim 6.}
Equations (\ref{Mlax}) and Claim 5 imply that $m_x$ is $C^1$ in $t$ and that $m_t$ is $C^1$ in $x$.
Hence the mixed partials $m_{xt}$ and $m_{tx}$ exist and are equal for $x \geq0$ and $t \geq 0$. The compatibility of (\ref{Mlax}) implies that $u$ satisfies (\ref{mkdv}).
\proofendcontinue

\medskip\noindent
{\bf Claim 7.} $u(x,0) = u_0(x)$ for $x \geq 0$.

{\it Proof of Claim 7.}
Consider the $x$-part (\ref{xpartF}) with potential $u(x)$ given by  $u_0(x)$. Let $X(x,k)$ and $Y(x,k)$ denote the associated eigenfunctions defined in Theorem \ref{xth1}. The results of Section \ref{xsec} imply that the function $m^{(x)}(x,k)$ defined by
\begin{align}
m^{(x)}(x,k) =
\begin{cases}
 \left(\frac{[Y(x,k)]_1}{a(k)}, [X(x,k)]_2\right), \qquad \im k < 0,
	\\
 \left([X(x,k)]_1, \frac{[Y(x,k)]_2}{\overline{a(\bar{k})}}\right), \qquad \im k > 0,
\end{cases}
\end{align}
satisfies the $L^2$-RH problem
\begin{align}\label{MxRH}
\begin{cases}
m^{(x)}(x, \cdot) \in I + \dot{E}^2(\C \setminus \R),\\
m_+^{(x)}(x,k) = m_-^{(x)}(x, k) \begin{pmatrix} 1 & - \frac{b(k)}{\overline{a(\bar{k})}} e^{2ikx} \\
\frac{\lambda \overline{b(\bar{k})}}{a(k)} e^{-2ikx} & \frac{1}{|a(k)|^2} \end{pmatrix}  \quad \text{for a.e.} \ k \in \R.
\end{cases}
\end{align}
Moreover, by (\ref{varphiasymptoticsa}) and the explicit expression (\ref{Fjexplicit}) for $X_1$, we have
\begin{align}\label{u0lim}
u_0(x) = -2i \lim_{k\to\infty} (km^{(x)}(x,k))_{12},
\end{align}
where the limit is taken along any direction in $\im k \leq 0$.

On the other hand, in view of Claim 1,
$$h(k) e^{-2ikx} \in \dot{E}^2(D_2) \cap E^\infty(D_2), \qquad \overline{h(\bar{k})} e^{2ikx} \in \dot{E}^2(D_3) \cap E^\infty(D_3).$$ 
Hence, Lemma \ref{deformationlemma}  together with the expression (\ref{Jdef}) for $J$ show that the function $M^{(x)}(x,k)$ defined by
$$M^{(x)}(x,k) = \begin{cases} M(x, 0, k), & k \in D_1 \cup D_4,
	\\
M(x,0,k) \begin{pmatrix} 1 & 0 \\ \lambda h(k) e^{-2ikx} & 1 \end{pmatrix}, & k \in D_2,
	\\
M(x,0,k) \begin{pmatrix} 1 & \overline{h(\bar{k})} e^{2ikx} \\ 0 & 1 \end{pmatrix}, & k \in D_3,
\end{cases}$$
also satisfies (\ref{MxRH}). By uniqueness, $M^{(x)} = m^{(x)}$. Comparing the definition (\ref{ulim}) of $u(x,t)$ with (\ref{u0lim}), we obtain $u(x,0) = u_0(x)$ for $x \geq 0$.
\proofendcontinue

\medskip\noindent
{\bf Claim 8.} $u(0,t) = g_0(t)$ and $u_x(0,t) = g_1(t)$ for $t \geq 0$.

\begin{figure}
\begin{center}
\bigskip \bigskip
\begin{overpic}[width=.45\textwidth]{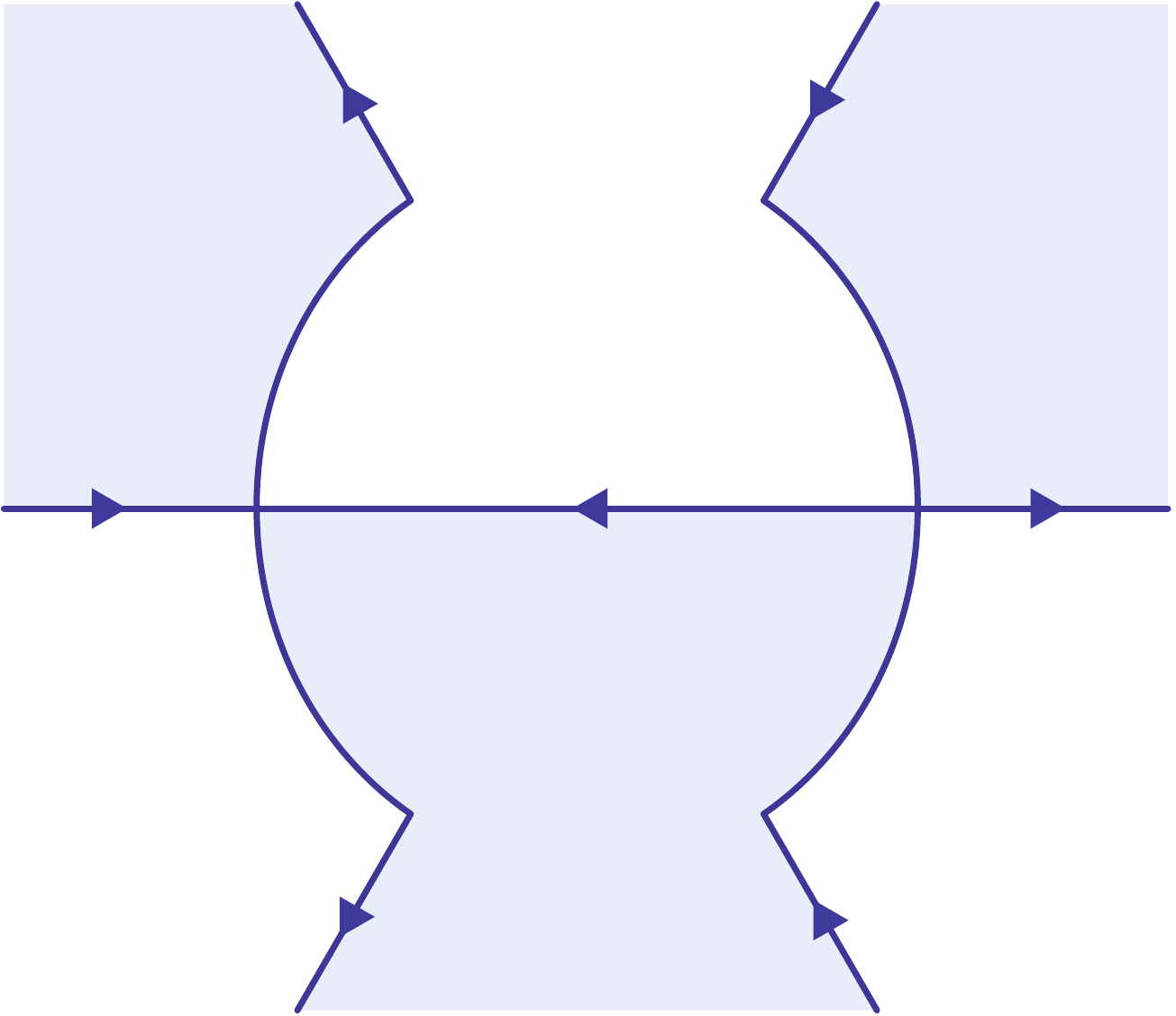}
      \put(82,67){$\mathcal{D}_3$}
      \put(46,60){$\mathcal{D}_4$}
      \put(10,67){$\mathcal{D}_3$}
      \put(10,15){$\mathcal{D}_2$}
      \put(46,22){$\mathcal{D}_1$}
      \put(82,15){$\mathcal{D}_2$}
      \put(102,41.2){$\Sigma$}
      \end{overpic}
     \begin{figuretext}\label{deformedDjs.pdf}
       The contour $\Sigma$ and the deformed domains $\{\mathcal{D}_j\}_1^4$ in the complex $k$-plane.
     \end{figuretext}
     \end{center}
\end{figure}

{\it Proof of Claim 8.}
Since $A(k) \to 1$ uniformly as $k \to \infty$, $k \in \bar{D}_3$, we can define deformed domains $\{\mathcal{D}_j\}_1^4$ so that $A(k)$ is nonzero in $\bar{\mathcal{D}}_3$, see Figure \ref{deformedDjs.pdf}.
We let $\Sigma = \R \cup (\bar{\mathcal{D}}_1 \cap \bar{\mathcal{D}}_2) \cup (\bar{\mathcal{D}}_3 \cap \bar{\mathcal{D}}_4)$ denote the contour separating the $\mathcal{D}_j$'s oriented as in Figure \ref{deformedDjs.pdf}.
We choose the $\mathcal{D}_j$'s so that $\Sigma$ is invariant under the involution $k \to \bar{k}$. 

Consider the $t$-part (\ref{tpartG}) defined in terms of $\{g_j(t)\}_0^2$. Let $T(t,k)$ and $U(t,k)$ denote the associated eigenfunctions defined in Theorem \ref{tth1}. 
The relation (\ref{GR}) and the condition $\det S(k) = 1$ imply
$$d(k) = a(k)\overline{A(\bar{k})} -  \lambda \frac{B(k) a(k)}{\overline{A(\bar{k})}} \overline{B(\bar{k})}
= \frac{a(k)}{A(k)}, \qquad k \in \partial D_1.$$
This shows that $A(k)$ admits an analytic continuation to $D_2$. Since $B = bA/a$ in $\bar{D}_1$ by (\ref{GR}), it follows that $B(k)$Ê also admits an analytic continuation to $D_2$.
The relation $T = U e^{-4ik^3t\hat{\sigma}_3}S(k)$ together with the assumption (\ref{GR}) yield the identity
$$\frac{[T(t,k)]_2}{A(k)} = \frac{b(k)}{a(k)} e^{-8ik^3t} [U(t,k)]_1 + [U(t,k)]_2, \qquad t \geq 0, \quad k \in \bar{D}_1.$$
Since $U(t, \cdot)$ is entire, this shows that $[T]_2/A$ also admits an analytic continuation to $D_2$. 
Hence we may define $m^{(t)}$ by
\begin{align}
m^{(t)}(t,k) =
\begin{cases}
 \left([U(t,k)]_1, \frac{[T(t,k)]_2}{A(k)} \right), \qquad k \in \mathcal{D}_1 \cup \mathcal{D}_3,
	\\
 \left(\frac{[T(t,k)]_1}{\overline{A(\bar{k})}}, [U(t,k)]_2 \right), \qquad k \in \mathcal{D}_2 \cup \mathcal{D}_4.
\end{cases}
\end{align}
The results of Section \ref{tsec} imply that $m^{(t)}(t,k)$ satisfies the $L^2$-RH problem
\begin{align}\label{MtRH}
\begin{cases}
m^{(t)}(t, \cdot) \in I + \dot{E}^2(\C \setminus \Sigma),\\
m^{(t)}_+(t,k) = m^{(t)}_-(t, k) \begin{pmatrix} 1 & \frac{B(k)}{A(k)} e^{-8ik^3t} \\
- \frac{\lambda\overline{B(\bar{k})}}{\overline{A(\bar{k})}} e^{8ik^3t} & \frac{1}{A(k)\overline{A(\bar{k})}} \end{pmatrix}  \quad \text{for a.e.} \ k \in \Gamma.
\end{cases}
\end{align}
Moreover, by (\ref{Gasymptoticsa}) and the explicit expression (\ref{Vjexplicit}) for $V_j$, $j = 1,2,3$, we have
\begin{subequations}\label{gjlim}
\begin{align}
& g_0(t) = -2i (m_1)_{12}(t),
	\\
& g_1(t) = 4(m_2)_{12}(t) - 2ig_0(t) (m_1)_{22}(t),
	\\
& g_2(t) = \lambda g_0^3(t) + 8i(m_3)_{12}(t) + 4g_0(t)(m_2)_{22}(t) - 2ig_1(t)(m_1)_{22}(t),
\end{align}
\end{subequations}
where 
$$m_j(t) = \ntlim_{k \to \infty} k^jm^{(t)}(t,k), \qquad j = 1,2,3,$$
and the limit is taken in a nontangential sector in $D_2 \cup D_4$.

On the other hand, deforming the contour from $\Gamma$ to $\Sigma$, we find that the function $\mathcal{M}$ defined by
$$\mathcal{M} = \begin{cases} M \begin{pmatrix} 1 & 0 \\ \lambda h(k) e^{-2ikx + 8ik^3t} & 1 \end{pmatrix}, & k \in \mathcal{D}_1 \cap D_2,
	\\
M \begin{pmatrix} 1 & \overline{h(\bar{k})} e^{2ikx - 8ik^3t} \\ 0 & 1 \end{pmatrix}, & k \in \mathcal{D}_4 \cap D_3,
 	\\
M, & \text{otherwise},
\end{cases}$$
satisfies 
\begin{align}\label{RHMcal}
\begin{cases}
\mathcal{M}(x, t, \cdot) \in I + \dot{E}^2(\C \setminus \Sigma),\\
\mathcal{M}_+(x,t,k) = \mathcal{M}_-(x, t, k) \mathcal{J}(x, t, k) \quad \text{for a.e.} \ k \in \Sigma,
\end{cases}
\end{align}
where $\mathcal{J}$ is defined by
\begin{align*}
&\mathcal{J}(x,t,k) = \begin{cases} 
 \begin{pmatrix} 1 & 0 \\ \lambda h(k) e^{-2ikx + 8ik^3t} & 1 \end{pmatrix}, & k \in \bar{\mathcal{D}}_1 \cap \bar{\mathcal{D}}_2,
	\\
 \begin{pmatrix} 1 & - \overline{r(\bar{k})} e^{2ikx - 8ik^3t} \\
\lambda r(k)e^{-2ikx + 8ik^3t}& 1 -  \lambda|r(k)|^2\end{pmatrix}, & k \in \bar{\mathcal{D}}_2 \cap \bar{\mathcal{D}}_3,
	\\ 
 \begin{pmatrix} 1 & - \overline{h(\bar{k})} e^{2ikx - 8ik^3t} \\ 0 & 1 \end{pmatrix}, & k \in \bar{\mathcal{D}}_3 \cap \bar{\mathcal{D}}_4,
 	\\
\begin{pmatrix} \frac{1}{|a(k)|^2} & \frac{b(k)}{\overline{a(\bar{k})}} e^{2ikx-8ik^3t} \\
- \frac{\lambda \overline{b(\bar{k})}}{a(k)} e^{-2ikx+8ik^3t} & 1 \end{pmatrix},
& k \in \bar{\mathcal{D}}_4 \cap \bar{\mathcal{D}}_1.
\end{cases}
\end{align*}
Furthermore, since $a(k)$, $d(k)$, and $A(k)$Ê are nonzero in $\bar{\mathcal{D}}_1 \cup \bar{\mathcal{D}}_2$, $\bar{\mathcal{D}}_2$, and $\bar{\mathcal{D}}_3$, respectively, we have
$$\bigg(\frac{d(k)}{\overline{A(\bar{k})}}\bigg)^{\pm 1}  \in 1 + (\dot{E}^2 \cap E^\infty)(\mathcal{D}_2), \qquad
a(k)^{\pm 1} \in 1 + (\dot{E}^2 \cap E^\infty)(\mathcal{D}_1 \cap \mathcal{D}_2).$$
It follows that the functions
\begin{align}\nonumber
& G_1(t,k) =  \begin{pmatrix} a(k) & 0 \\ 0 & \frac{1}{a(k)} \end{pmatrix}, 
\qquad G_2(t,k) = \begin{pmatrix} \frac{d(k)}{\overline{A(\bar{k})}} & - b(k) e^{-8ik^3t} \\ 0 & \frac{\overline{A(\bar{k})}}{d(k)} \end{pmatrix}, 
	\\ \label{Gjdef}
& G_3(t,k) = \begin{pmatrix} \frac{A(k)}{\overline{d(\bar{k})}} & 0 \\ - \lambda \overline{b(\bar{k})} e^{8ik^3t} & \frac{\overline{d(\bar{k})}}{A(k)} \end{pmatrix},
\qquad
G_4(t,k) = \begin{pmatrix} \frac{1}{\overline{a(\bar{k})}} & 0 \\ 0 & \overline{a(\bar{k})} \end{pmatrix},
\end{align}
satisfy $G_j(t,\cdot) \in \dot{E}^2(\mathcal{D}_j) \cap E^\infty(\mathcal{D}_j)$ for $j = 1,\dots, 4$.
Thus a contour deformation argument together with the expression (\ref{Jdef}) for $J$ show that the function $M^{(t)}(t,k)$ defined by
$$M^{(t)}(t,k) = \mathcal{M}(0,t,k)G_j(t,k), \qquad k \in \mathcal{D}_j, \quad j = 1,\dots,4,$$
also satisfies (\ref{MtRH}). 
By uniqueness, $M^{(t)} = m^{(t)}$. Comparing the definition (\ref{ulim}) of $u(x,t)$ with (\ref{gjlim}), we obtain $u(0,t) = g_0(t)$ for $t \geq 0$. Comparing (\ref{gjlim}) with equation (\ref{cadb}) in the defocusing case, or with an equation analogous to (\ref{cadb}) in the focusing case, we obtain also $u_x(0,t) = g_1(t)$ and $u_{xx}(0,t) = g_2(t)$ for $t \geq 0$.
\proofend

\begin{remark}\upshape
In the focusing (i.e. $\lambda = -1$) case, the function $d(k)$ could have zeros. The assumption in Theorem \ref{existenceth} that $d(k)$ is nonzero in $\bar{D}_2$ is made purely for convenience. Indeed, suppose $d(k)$ has zeros in $\bar{D}_2$. Since $d(k) \to 1$ as $k \to \infty$, we can define deformed domains $\{\mathcal{D}_j\}_1^4$ as in Figure \ref{deformedDjs.pdf} such that $d(k)$ is nonzero in $\mathcal{D}_2$. By formulating the RH problem (\ref{RHM}) in terms of the deformed domains $\{\mathcal{D}_j\}$ instead of $\{D_j\}$, the same proof goes through with obvious modifications.
\end{remark}

\begin{remark}\upshape
In the focusing (i.e. $\lambda = -1$) case, the function $a(k)$ could have zeros. The assumption in Theorem \ref{existenceth} that $a(k)$ is nonzero inÊ $\bar{\C}_-$ can be weakened. In particular, if $a(k)$Ê has a finite number of simple poles, then these poles can be easily treated as in \cite{BFS2004}.  
\end{remark}

\begin{remark}\upshape
Let us comment on the existence of a vanishing lemma for the RH problem (\ref{RHM}). Zhou showed in \cite{Z1989} (see Theorem 9.3 of \cite{Z1989}) that a homogeneous RH problem with a Schwarz reflection invariant contour has only the trivial solution provided that the jump matrix $J(k)$ Êsatisfies (a) $J(k) = \overline{J(\bar{k})}^T$ for $k \in \Gamma \setminus \R$ and (b) $\re J(k)$ is positive definite on $\R$.\footnote{The real and imaginary parts of a square matrix $A$ are defined by $\re A = \frac{1}{2}(A + \bar{A}^T)$ and $\im A = \frac{1}{2i}(A - \bar{A}^T)$.} In the focusing (i.e. $\lambda = -1$) case, the jump matrix (\ref{Jdef}) satisfies these conditions, so there exists a vanishing lemma. 
In the defocusing (i.e. $\lambda = 1$) case, we have not been able to establish a vanishing lemma. Let us explain the main difficulties. Let $\lambda = 1$. Then the condition (a) is not satisfied; instead the jump matrix $J$ satisfies $J(k) = \sigma_3 \overline{J(\bar{k})}^T \sigma_3$. This implies that if $M(k)$ is a solution of the homogeneous RH problem, then the function $M(k)\sigma_3 \overline{M(\bar{k})}^T$ has no jump across $\Gamma$; hence it vanishes identically. In particular,
$$M_-(k)J(k) \sigma_3 \overline{M_-(k)}^T = M_-(k) \sigma_3 \overline{J(k)}^T \overline{M_-(k)}^T = 0, \qquad k \in \R.$$
However, since no linear combination of $J \sigma_3$ and $\sigma_3 \bar{J}^T$ is positive definite, we cannot conclude that $M_- = 0$. 
More generally, proceeding as in the appendix of \cite{FZ1992}, we may consider a sectionally analytic function $G(k)$ defined by 
$$G(k) = M(k) H_j(k) \overline{M(\bar{k})}^T, \qquad k \in D_j, \quad j = 1,\dots, 4,$$ 
where $H_j(k)$ is any $2 \times 2$-matrix valued function which is analytic and bounded in $D_j$.
The requirement that $G$ be continuous across $\partial D_1$ enforces the condition $J(k) H_1(k) = H_2(k) \overline{J(\bar{k})}^T$ on $\partial D_1$, i.e.
$$H_2(k) = \begin{pmatrix} 1 & 0 \\ h(k) e^{-2ikx + 8ik^3t} & 1 \end{pmatrix} H_1(k) \begin{pmatrix} 1 & 0 \\ h(k) e^{-2ikx + 8ik^3t} & 1 \end{pmatrix}.$$
Since we want $H_2$ to be analytic and bounded in $D_2$ we assume that $H_1$ is traceless and lower triangular. 
Similarly, 
continuity across $\partial D_4$ enforces the condition $J(k) H_3(k) = H_4(k) \overline{J(\bar{k})}^T$ on $\partial D_4$, so we need to assume that $H_4$ is traceless and upper triangular. 
Then $H_2$ and $H_3$ are also traceless and lower resp. upper triangular.
The equations
\begin{align*}
& 0 = \int_\R G_+(k) dk = \int_\R M_-(k)J(k)H_3(k)\overline{M_-(k)}^T dk, \qquad 
	\\
& 0 = \int_\R G_-(k) dk = \int_\R M_-(k)H_2(k)\overline{J(k)}^T\overline{M_-(k)}^Tdk,
\end{align*}
imply that
$$ \int_\R M_- X \overline{M_-}^T dk = 0, \qquad X := c_1 JH_3 + c_2 \overline{H_3}^T\overline{J}^T
+ c_3 H_2\overline{J}^T + c_4 J \overline{H_2}^T,$$
for any complex constants $\{c_j\}_1^4$.
We seek $\{c_j\}_1^4$ such that $X$ is positive definite for $k \in \R$. Since $X$ is positive definite iff $\re X$ is positive definite, it is enough to consider positive definiteness of
$$Y := d_1 \re(JH_3) + d_2 \im(JH_3) + d_3 \re(H_2\overline{J}^T)
+ d_4 \im(H_2\overline{J}^T),$$
where $\{d_j\}_1^4$ are real coefficients. The matrix $Y$ is positive definite iff the inequalities $Y_{11}> 0$ and $\det Y > 0$ are satisfied. However, it appears that no real numbers $\{d_j\}_1^4$ and traceless lower resp. upper triangular matrices $H_2(k)$ and $H_3(k)$ fulfill these inequalities. 
\end{remark}

\appendix
\section{$L^2$-Riemann-Hilbert problems} \label{RHapp}
\renewcommand{\theequation}{A.\arabic{equation}}\nequation
We collect some results on Smirnoff classes and $L^2$-RH problems; detailed proofs can be found in \cite{LCarleson}. In the context of smooth contours, more information on $L^2$-RH problems can be found in \cite{D1999, DZ2002b, FIKN2006, Z1989}.

Assume $\Gamma = \partial D_+ = -\partial D_-$ is a Carleson jump contour. If $f \in \dot{E}^2(D_+)$ or  $f \in \dot{E}^2(D_-)$, then the nontangential limits of $f(z)$ as $z$ approaches the boundary exist a.e. on $\Gamma$ and the boundary function belongs to $L^2(\Gamma)$.
If $h \in L^2(\Gamma)$, then the Cauchy transform $\mathcal{C}h$ defined by
\begin{align}\label{Cauchytransform}
(\mathcal{C}h)(z) = \frac{1}{2\pi i} \int_\Gamma \frac{h(s)}{s - z} ds, \qquad z \in \C \setminus \Gamma,
\end{align}
satisfies $\mathcal{C}h \in \dot{E}^2(D_+ \cup D_-)$. 
We denote the nontangential boundary values of $\mathcal{C}f$ from the left and right sides of $\Gamma$ by $\mathcal{C}_+ f$ and $\mathcal{C}_-f$ respectively. Then $\mathcal{C}_+$ and $\mathcal{C}_-$ are bounded operators on $L^2(\Gamma)$ and $\mathcal{C}_+ - \mathcal{C}_- = I$.
Given two functions $w^\pm \in L^2(\Gamma) \cap L^\infty(\Gamma)$, we define the operator $\mathcal{C}_{w}: L^2(\Gamma) + L^\infty(\Gamma) \to L^2(\Gamma)$ by 
\begin{align}\label{Cwdef}
\mathcal{C}_{w}(f) = \mathcal{C}_+(f w^-) + \mathcal{C}_-(f w^+).
\end{align}
Then
\begin{align}\label{Cwnorm}
\|\mathcal{C}_w\|_{\mathcal{B}(L^2(\Gamma))} \leq C \max\big\{\|w^+\|_{L^\infty(\Gamma)}, \|w^-\|_{L^\infty(\Gamma)} \big\}.
\end{align}
where $C = \max\{\|\mathcal{C}_+\|_{\mathcal{B}(L^2(\Gamma))}, \|\mathcal{C}_-\|_{\mathcal{B}(L^2(\Gamma))}\} < \infty$ and $\mathcal{B}(L^2(\Gamma))$ denotes the Banach space of bounded linear maps $L^2(\Gamma) \to L^2(\Gamma)$.

The next lemma shows that if $v = (v^-)^{-1}v^+$ and $w^\pm = \pm v^\pm \mp I$ then the $L^2$-RH problem determined by $(\Gamma, v)$ is equivalent to the following singular integral equation for $\mu \in I + L^2(\Gamma)$:
\begin{align}\label{rhoeq}
\mu - I = \mathcal{C}_w(\mu)  \quad \text{in}\quad L^2(\Gamma).
\end{align}

\begin{lemma}\label{mulemma}
Given $v^\pm: \Gamma \to GL(n, \C)$, let $v = (v^-)^{-1}v^+$, $w^+ = v^+ - I$, and $w^- = I - v^-$. Suppose $v^\pm, (v^\pm)^{-1} \in I +  L^2(\Gamma) \cap L^\infty(\Gamma)$.
If $m \in I + \dot{E}^2(\C \setminus \Gamma)$ satisfies the $L^2$-RH problem determined by $(\Gamma, v)$, then $\mu = m_+ (v^+)^{-1} = m_- (v^-)^{-1} \in I + L^2(\Gamma)$ satisfies (\ref{rhoeq}). 
Conversely, if $\mu \in I + L^2(\Gamma)$ satisfies (\ref{rhoeq}), then
$m = I + \mathcal{C}(\mu(w^+ + w^-)) \in I + \dot{E}^2(\C \setminus \Gamma)$ satisfies the $L^2$-RH problem determined by $(\Gamma, v)$. 
\end{lemma}

\begin{lemma}\label{Fredholmzerolemma}
Given $v^\pm: \Gamma \to GL(n, \C)$, let $v = (v^-)^{-1}v^+$, $w^+ = v^+ - I$, and $w^- = I - v^-$. Suppose $v^\pm, (v^\pm)^{-1} \in  I + L^2(\Gamma) \cap L^\infty(\Gamma)$ and $v^\pm \in C(\Gamma)$.
If $w^\pm$ are nilpotent matrices, then each of the following four statements implies the other three:
\begin{enumerate}[$(a)$]
\item The map $I - \mathcal{C}_w:L^2(\Gamma) \to L^2(\Gamma)$ is bijective.

\item The $L^2$-RH problem determined by $(\Gamma, v)$ has a unique solution. 

\item The homogeneous $L^2$-RH problem determined by $(\Gamma, v)$ has only the zero solution.

\item  The map $I - \mathcal{C}_w: L^2(\Gamma) \to L^2(\Gamma)$ is injective.
\end{enumerate}
\end{lemma}

\begin{lemma}[Uniqueness]\label{uniquelemma}
Suppose $v: \Gamma \to GL(2, \C)$ satisfies $\det v = 1$ a.e. on $\Gamma$.
If the solution of the $L^2$-RH problem determined by $(\Gamma, v)$ exists, then it is unique and has unit determinant.
\end{lemma}

\begin{lemma}\label{EpCnlemma}
Let $D$ be a subset of $\hat{\C}$ bounded by a curve $\Gamma \in \mathcal{J}$ and let $f:D \to \C$ be an analytic function.
Suppose there exist curves $\{C_n\}_1^\infty \subset \mathcal{J}$ in $D$, tending to $\Gamma$ in the sense that $C_n$ eventually surrounds each compact subset of $D \subset \hat{\C}$, such that $\sup_{n \geq 1} \int_{C_n} |z - z_0|^{p-2} |f(z)|^p |dz| < \infty$. Then $f \in \dot{E}^p(D)$. 
\end{lemma}
 
 \begin{lemma}[Contour deformation]\label{deformationlemma}
Let $\gamma \in \mathcal{J}$. Suppose that, reversing the orientation on a subcontour if necessary, $\hat{\Gamma} = \Gamma \cup \gamma$ is a Carleson jump contour. 
Let $B_+$ and $B_-$ be the two components of $\hat{\C} \setminus \gamma$. Let $\hat{D}_\pm$ be the open sets such that $\hat{\C} \setminus \hat{\Gamma} = \hat{D}_+ \cup \hat{D}_-$ and $\partial \hat{D}_+ = - \partial \hat{D}_- = \hat{\Gamma}$. 
Let $\hat{D} = \hat{D}_+ \cup \hat{D}_-$.
Let $\gamma_+$ and $\gamma_-$ be the parts of $\gamma$ that belong to the boundary of $\hat{D}_+ \cap B_+$ and $\hat{D}_- \cap B_+$, respectively. 
Suppose $v: \Gamma \to GL(n, \C)$. Suppose $m_0:\hat{D} \cap B_+ \to GL(n,\C)$ satisfies
\begin{align}\label{m0m0inv}
m_0, m_0^{-1} \in I + \dot{E}^2(\hat{D} \cap B_+) \cap E^\infty(\hat{D} \cap B_+).
\end{align}
Define $\hat{v}:\hat{\Gamma} \to GL(n, \C)$ by
\begin{align*}
\hat{v} 
=  \begin{cases}
 m_{0-} v m_{0+}^{-1} & \text{on} \quad  \Gamma \cap B_+, \\
m_{0+}^{-1} & \text{on} \quad \gamma_+, \\
m_{0-} & \text{on} \quad \gamma_-, \\
v & \text{on} \quad \Gamma \cap B_-.
\end{cases}
\end{align*}
Let $m$ and $\hat{m}$ be related by
\begin{align}\label{hatmdefmm0}
\hat{m} = \begin{cases}
mm_0^{-1} & \text{on} \quad \hat{D} \cap B_+,\\
m & \text{on} \quad \hat{D} \cap B_-.
\end{cases}
\end{align}
Then $m(z)$ satisfies the $L^2$-RH problem determined by $(\Gamma,v)$ if and only if $\hat{m}(z)$ satisfies the $L^2$-RH problem determined by $(\hat{\Gamma}, \hat{v})$.
\end{lemma}

\begin{lemma}\label{intfmlemma}
Let $D$ be a subset of $\hat{\C}$ bounded by a curve $\Gamma \in \mathcal{J}$. Let $m \geq 1$ be an integer and let $\{f_j\}_1^m$ be functions in $\dot{E}^m(D)$. Let $f = \prod_{j = 1}^m f_j$. Then $f \in \dot{E}^1(D)$ and
\begin{align}\label{intfm}
\int_\Gamma \frac{f_+(z)}{(z-z_0)^{n-m+2}} dz = 0, \qquad n = 0,1, \dots.
\end{align}
\end{lemma}
\proofbegin
Let $z_0 \in \C \setminus \bar{D}$ and let $\varphi(z) = \frac{1}{z-z_0}$.
By Proposition 3.6 in \cite{LCarleson}, the assumption $f_j \in \dot{E}^m(D)$ is equivalent to $w^{-1} f_j(\varphi^{-1}(w)) \in \dot{E}^m(\varphi(D))$. Thus, by Lemma 3.8 in \cite{LCarleson}, $w^{-m} f(\varphi^{-1}(w))^m \in \dot{E}^1(\varphi(D)) = E^1(\varphi(D))$.
According to Theorem 10.4 of \cite{D1970}, 
$$\int_{\varphi(\Gamma)} w^{n-m} f_+(\varphi^{-1}(w))^m dw = 0, \qquad n = 0, 1, \dots.$$
The change of variables $w = \varphi(z)$ yields (\ref{intfm}).
\proofend

\bigskip
\noindent
{\bf Acknowledgement} {\it The author acknowledges support from the EPSRC, UK.}

\bibliographystyle{plain}
\bibliography{is}

\end{document}